\documentclass[a4paper,10pt]{article}
\usepackage{amssymb}
\usepackage{latexsym}
\usepackage{amsmath}
\usepackage{amsthm}
\usepackage{amsfonts}
\usepackage{amssymb}
\usepackage{amsmath,amscd}
\usepackage{graphicx}

\newtheorem{Th}{Theorem}
\newtheorem{Prop}{Proposition}

\newtheorem{Lm}{Lemma}
\newtheorem{Lma}{Lemma}[section]
\newtheorem{Dfi}{Definition}
\newtheorem{Rm}{Remark}
\newtheorem{Con}{Open Problem}
\newcommand{\be}{\begin{equation}}
\newcommand{\ee}{\end{equation}}
\newcommand{\bes}{\begin{equation*}}
\newcommand{\ees}{\end{equation*}}

\newcommand{\R}{\mathbb{R}}
\newcommand{\N}{\mathbb{N}}
\newcommand{\C}{\mathbb{C}}
\newcommand{\Z}{\mathbb{Z}}

\newcommand\res{\mathop{\hbox{\vrule height 7pt width .5pt depth 0pt
\vrule height .5pt width 6pt depth 0pt}}\nolimits}

\def\theequation{\thesection.\arabic{equation}}
\def\theTh{\Roman{section}.\arabic{Th}}
\def\theProp{\Roman{section}.\arabic{Prop}}
\def\theCo{\Roman{section}.\arabic{Co}}

\def\theRm{\Roman{section}.\arabic{Rm}}

\newcommand{\reset}{\setcounter{equation}{0}\setcounter{Th}{0}\setcounter{Prop}{0}\setcounter{Co}{0}
\setcounter{Lm}{0}\setcounter{Rm}{0}}
\def\La{\Lambda}
\def\La{\Lambda}
\def\ti{\tilde}
\def\lf{\left}
\def\rg{\right}

\def\al{\alpha}
\def\la{\lambda}

\def\ep{\varepsilon}

\def\ds{\displaystyle}
\def\ov{\overline}
\def\Om{\Omega}
\def\om{\omega}
\def\p{\partial}

\def\res{\mathop{\hbox{\vrule height 7pt width .5pt 
depth 0pt\vrule height .5pt width 6pt depth 0pt}}\nolimits}

\begin{document}

\title{  Minmax Hierarchies and Minimal Surfaces in Manifolds}

\author{ Tristan Rivi\`ere\footnote{Department of Mathematics, ETH Zentrum,
CH-8093 Z\"urich, Switzerland.}}

\maketitle

{\bf Abstract :}{\it We introduce a general scheme that permits to generate successive min-max problems for producing critical points of higher and higher indices to Palais-Smale Functionals in Banach manifolds equipped with Finsler structures. We call the resulting tree of minmax problems a {\bf minmax hierarchy}. Using the  viscosity approach to the minmax theory of minimal surfaces introduced by the author
in a series of recent works, we explain how this scheme can be deformed for producing smooth minimal surfaces of strictly increasing area  in arbitrary codimension.
We implement this scheme to the case of the 3-dimensional sphere. In particular we are giving a min-max characterization of  the Clifford Torus and conjecture what are the next minimal surfaces  to come in the $S^3$ hierarchy.  Among other results we prove here the lower semi continuity of the Morse Index in the viscosity method below an area level.}

\medskip

\noindent{\bf Math. Class. 49Q05, 53A10,  58E12, 49Q10}

\section{Introduction}

A classical variational approach to the eigenvalue problem for the Laplacian on a closed oriented riemannian manifold $(N^n,h)$  is given by the {\it Rayleigh quotient method}.
It can be sketched as follows.

Introduce the {\it Hilbert Sphere}
\[
{\mathfrak S}:=\lf\{u\in W^{1,2}(N^n)\ ;\ \int_{N^n}|u|^2 \ dvol_h=1\rg\}
\]
We consider first
\[
\la_1:=\inf_{u\in {\mathfrak S}}E(u):=\int_{N^n}|du|^2 \ dvol_h\quad\mbox{ and }\quad{\mathcal C}_1:=\lf\{u\in {\mathfrak S}\ ;\ \Delta_h u=\la_1 u\rg\}
\]
where $\Delta_h$ is the positive {\it Beltrami Laplace Operator} on $N^n$. Iteratively we introduce $E_{k-1}:=\mbox{Span}\,{{\mathcal C}_{k-1}}\oplus E_{k-2}$ (with the convention $E_0=\{0\}$)
\[
\la_k:=\inf_{u\in {\mathfrak S}\cap E_{k-1}^\perp}E(u):=\int_{N^n}|du|^2 \ dvol_h\quad\mbox{ and }\quad{\mathcal C}_k:=\lf\{u\in {\mathfrak S}\ ;\ \Delta_h u=\la_k u\rg\}\simeq S^{n_k-1}
\]
where $E_{k-1}^\perp$ is the space orthogonal to $E_{k-1}$ for the $L^2$ scalar product and $n_{k}:=\mbox{dim} \,\mbox{Span}({\mathcal C}_{k})$.

Assume now that we don't want to make use neither of the linear nature of the problem nor on the existence of the scalar product. An alternative way to obtain the {\it Laplace Eigenspaces}
and {\it Laplace Eigenvalues} is given by what we call a {\it Minmax Hierarchy}. A {\it Minmax Hierarchy} in this  framework is the following iterative construction.
Starting from $\la_1$, ${\mathcal C}_1=\{-u^1,u^1\}$ and $1=n_1:=\mbox{dim} \,\mbox{Span}({\mathcal C}_1)$ which were obtained by a strict minimization of the {\it Dirichlet energy} in ${\mathfrak S}$, in order to produce $\la_2$ without using the scalar product one could proceed as follows. Introduce
\[
\mbox{Sweep}_1:=\lf\{ u\in \mbox{Lip}([-1,+1],{\mathfrak S}) \quad\mbox{ s. t. }\quad u_{-1}=-u^{1}\quad\mbox{ and }\quad u_{+1}:=u^1\rg\}
\]
One has
\[
\la_2=\inf_{u\in \mbox{Sweep}_1}\ \max_{y\in[-1,+1]} \, E(u_y)
\]
Using the classical {\it Palais deformation theory} one produces critical point of $E$ in ${\mathfrak S}$ realizing $\la_2$. Introduce ${\mathcal C}_2:=\lf\{u\in {\mathcal S}\ ;\ \Delta_h u=\la_2 u\rg\}$ and call $n_2:=\mbox{dim} \,\mbox{Span}({\mathcal C}_{2})$ and $N_2=n_1+n_2$. It is clear that by taking
a minmax based on the space of 1 dimensional paths connecting elements from ${\mathcal C}_2$ would lead to nothing since points in ${\mathcal C}_2$ can be connected within ${\mathcal C}_2\simeq S^{n_2-1}$ (i.e. $n_2>1$). We introduce instead $Y_2:=B^{n_2}\times B^{n_1}$
\[
\mbox{Sweep}_{N_2}:=\lf\{ 
\begin{array}{l}
\ds u\in \mbox{Lip}(Y_2,{\mathcal S}) \ ;\ \forall \, z\in \p B^{n_2}\quad u_{(z,\cdot)}\in \mbox{Sweep}_1 \\[3mm]
(u_{\cdot})^{-1}({\mathcal C}^2)\cap \p Y_2=\p B^{n_2}\times \{0\}\quad;\quad \mbox{deg}(u_{|\p B^{n_2}\times \{0\}})=+1
\end{array}
\rg\}
\]
One then proves (see subsection~\ref{rayleigh})
\[
\la_3=\inf_{u\in \mbox{Sweep}_{N_2}}\max_{y\in Y_2} E(u_y)
\]
We construct inductively 
\begin{itemize}
\item[i)] A sequence of integers $n_k\in {\N}^\ast$ 
\[
n_k:=\mbox{dim} \,({\mathcal C}_k)+1
\]
\item[ii)] A sequence $\mbox{Sweep}_{N_k}\subset \mbox{Lip}(Y_k,{\mathcal S})$ where $Y_k:= B^{n_k}\times B^{n_{k-1}}\cdots B^{n_1} $ characterized as follows
\[
\forall\, u\in\mbox{Sweep}_{N_k}\quad \forall z\in  \p B^{n_k}\quad\quad u_{(z,\cdot)}\in \mbox{Sweep}_{N_{k-1}}
\]
\item[iii)] we have $u^{-1}({\mathcal C}_{k-1})\cap \lf(\p B^{n_k}\times Y_{k-1}\rg)= \p B^{n_k}\times\{0\}$, 
\[
\mbox{deg}(u_{|_{\p B^{n_k}\times\{0\}}})=+1
\]
\item[iv)] for $\mbox{dim}\,Y_{k-1}>0$ (i.e. $k>2$)
\[
\max_{y\in B^{n_k}\times\ \p Y_{k-1}} \int_{N^n}|du_y|^2 \ dvol_h< \la_{k-1}
\]
\end{itemize}
With this definition for $\mbox{Sweep}_{N_k}$ one proves (see subsection~\ref{rayleigh})
\[
\la_k=\inf_{u\in \mbox{Sweep}_{N_k}}\max_{y\in Y_k} \int_{N^n}|du_y|^2 \ dvol_h
\]

\medskip

The goal of the present work is first to extend the notion of {\it Minmax hierarchies} to the general non-linear framework of a Banach Manifold for Lagrangians satisfying the {\it Palais-Smale condition}. In the second part of the work we establish that the above definition of a {\it Minmax Hierarchy} extends to the framework
of {\it minimal surfaces} in $M^m$ a closed sub-manifold of an euclidian space ${\R}^Q$ and to illustrate this definition by concrete construction of Hierarchies for {\it minimal surfaces}
in the sphere $S^3$. The {\it Hilbert sphere} ${\mathfrak S}$ introduced for the eigenvalue problem above is replaced by
the union of the Hilbert manifolds made of $W^{3,2}-$immersions of all the oriented closed surfaces. 
\[
\mbox{Imm}(M^m):=\bigcup_{g\in {\N}}\mbox{Imm}(\Sigma^g,M^m)
\] 
More precisely we are considering the quotients of these spaces by, $\mbox{Diff}^{\,\ast}_+(\Sigma^g)$ the marked positive diffeomorphisms of the $\Sigma_g$, the positive diffeomorphisms fixing a given point
\[
{\mathfrak M}(M^m):=\bigcup_{g\in {\N}}\mbox{Imm}(\Sigma^g,M^m)/\mbox{Diff}^{\,\ast}_+(\Sigma^g)\quad.
\]
The main idea is to take advantage of the topology of the space of critical points at the rank $k-1$ in order
to produce a minmax problem with a Width 
\[
W_k(M^m):=\inf_{\vec{\Phi}\in \mbox{Sweep}_{N_k}}\ \max_{y\in Y_k}\, \mbox{Area}\,(\vec{\Phi}(y))
\]
 strictly larger than the $k-1-$one :
 \be
 \label{00.strict}
 W_{k-1}(M^m)<W_k(M^m)
 \ee
 The sweep-out {\it admissible families} $\mbox{Sweep}_{N_k}$ are made of maps from an $N_k$ dimensional polyhedral chain into ${\mathfrak M}$ and the successive integer $N_k$ are called  {\it Minmax Indices} of the hierarchy. In theorem~\ref{th-I.1} we prove that the definition of the {\it Minmax Hierarchy} corresponding to the set of conditions  i)...iv) but for
 the space ${\mathfrak M}$ and the area gives the strict inequality (\ref{00.strict}). Then, we implement the {\it viscosity method}\footnote{The advantage of the viscosity method over
 previous minmax methods for minimal surfaces such as the GMT approach of Almgren-Pitts or the Dirichlet energy method of Colding Minicozzi is threefold :
 it gives regularity in arbitrary co-dimensions, it gives automatically a genus upper-bound, it gives automatically an index upper-bound in terms of the dimension
 of the parameter space.} introduced
 in \cite{Riv-minmax} in order to produce smooth, possibly branched, minimal immersions realizing the widths $W_{k}(M^m)$. In the 
 case of the eigenvalue problem above, the simplicity of the topology of the {\it Hilbert sphere} ${\mathfrak S}$ makes the successive existence of the
 spaces $\mbox{Sweep}_{N_k}$  rather straightforward. This is not the case anymore for the space ${\mathfrak M}$. The complexity of the topology of this
 {\it Hilbert manifold} is the main difficulty for producing an infinite hierarchy and this question should be studied more systematically in forthcoming works. 
 
 In order to illustrate the notion of {\it Minmax Hierarchy} in the framework of minimal surfaces we are taking the simplest possible case : that is $M^m=S^3$. We construct the 2 first elements of the hierarchy which should be branching staring at the 2nd level and we have 
  \[
 W_1(S^3)=4\pi<W_2(S^3)=2\pi^2< W_3(S^3)
 \]
 moreover the {\it minmax indices}  are 
 \[
 N_1=1\quad,\quad N_2=5\quad\mbox{ and }\quad N_3=8\quad.
  \]
 We are proving in section V that the first {\it width} $W_1(S^3)$ is exactly achieved by the {\it geodesic spheres} (index 1), the second by the {\it Clifford tori} (index 5) and we conjecture that the third and the forth are achieved by minimal surfaces of genus  2 both and index respectively  8 and 9. 
 
 \medskip
 
 While it is clear that $W_1(S^n)=4\pi$ for any $n\ge 3$ it would be interesting to study the value $W_2(S^4)$. We clearly have $W_2(S^4)\le 2\pi^2$. Observe that the viscosity method prevents the {\it Veronese Surface} of area $6\pi$ to be obtained in the {\it minmax hierarchy}, as we define it below, since one has excluded immersions of non-oriented surfaces (only an even cover of it, with area at least $12\pi$, can be obtained).

 \medskip
 
  It should be noted that  the passage from $\mbox{Sweep}_{N_{k-1}}(M^m)$ to $\mbox{Sweep}_{N_{k}}(M^m)$ requires the choice of a non zero element\footnote{Observe that in the definition of the hierarchy for the eigenvalue problem above we were considering the cohomology for the group $G=\Z$. Because of orientation problems explained below we
  had to reduce to the more coarse group $G={\Z}_2$ for the minimal surface case.}  of $H^{\ast}({\mathcal C}_{k-1},{\Z}_2)$. While this choice was more or less unique in the eigenvalue problem above since ${\mathcal C}_{k-1}\simeq S^{n_k-1}$, this is not anymore the case
  for the minimal surface problem. Hence there should be as many $\mbox{Sweep}_{N_{k}}(M^m)$ as choices of non zero elements in $H^{\ast}({\mathcal C}_{k-1},{\Z}_2)$. In that
  sense a {\it Minmax Hierarchy} should be seen as being a {\bf tree} more than a { sequence} - which was specific to the ``linear problem'' of finding the eigenvalues. There could moreover be non-isotopic components in $\mbox{Sweep}_{N_{k}}(M^m)$ each generating a new branch of the hierarchy.
 
Coming back to the $S^3$ case, starting at the rank $k\ge 3$ if the minimal surfaces realizing $W_k(S^3)$ have no special continuous symmetry coming from the action of $SO(4)$ (unlike the geodesic spheres and the Clifford tori for the rank respectively 1 and 2) hence we expect each connected components of ${\mathcal C}_k$ to be diffeomorphic to $SO(4)\simeq 
SO(3)\times S^3$. Consequently, the only possibilities\footnote{We are using that the cohomology groups of $SO(4)$ are given by the truncated polynomial algebra
\[
H^\ast(SO(4),{\Z}_2)={\Z}_2[b_1,b_3]/(b_1^4,b_3^2)
\]}
 for the jumps in the dimension of parameter space are given by
\[
n_k\in\{ 1, 2,3, 4,5, 6, 7\}  \quad.
\]
In generic situation, for a manifold $M^m$ equipped with a metric such that the minimal surfaces are isolated (see on that subject \cite{Wh1}) we expect
\[
n_1=m-2\quad\mbox{ and }\quad\forall \,k>1 \quad\quad n_k=1\quad.
\]
One would for instance take ``eversion type'' sweep-outs $\mbox{Sweep}_{N_{k}}(M^m)=\mbox{Sweep}_{m-3+k}(M^m)$ whose ``core'' is for instance a path connecting a minimal surface with a given orientation and the same surface with the opposite orientation.

\medskip

One has to bound the genus in each class $\mbox{Sweep}_{N_{k}}(M^m)$ in order to have the realization of $W_k(M^m)$. This can be achieved by imposing this bound in the definition
of $\mbox{Sweep}_{N_{k}}(M^m)$. The need to increase the topology of the considered surfaces in a {\it Minmax Hierarchy} is sometimes required for passing form the rank $k$ to the rank $k+1$. As a matter of example, we can take again the $S^3$ case. On can restrict to genus 0 surfaces in the definition of $\mbox{Sweep}_{1}(S^3)$ but in order for the next space $\mbox{Sweep}_{5}(S^3)$ not to be empty we have to allow the genus to be non zero (to be one in fact : we have $\mbox{Sweep}_{5}(S^2,S^3)=\emptyset$ while $\mbox{Sweep}_{5}(T^2,S^3)\ne \emptyset$). This explains why, by using the freedom of increasing the genus while passing from one level of the hierarchy to the next, one could hope to produce an infinite hierarchy for a given space $M^m$. This would give minimal surfaces realizing the infinitely many successive widths $W_k(M^m)$.
A major question would then be to understand asymptotically the behavior of $W_k(M^m)$ in a given branch of the hierarchy. In case $W_k(M^m)$ would have a {\bf sub-linear growth} with respect to $k$ one would immediately deduce the existence of infinitely many  closed minimal surfaces in $M^m$.

\medskip
 
 One could a-priori define {\it Minmax hierarchies} for {\it non oriented surfaces} or for {\it surfaces with boundary} (the last case being interesting in relation with the search of free boundary surfaces). One would then have to extend the viscosity method to these two frameworks. In the simplest framework $M^m=B^3$ one would consider for $\mbox{Sweep}_1(B^3)$ families of immersions of the disc $D^2$ intersecting $\p B^3$ orthogonally at their boundary and realizing a sweep out of $B^3$. The {\it free boundary oriented discs} realize obviously the first step of the hierarchy. Hence one has $$W_1(B^3)=\pi\quad,\quad n_1=1\quad,\quad {\mathcal C}_1\simeq S^2\quad\mbox{and }\quad n_2=3\quad.$$ One would probably need to ``increase the topology'' by considering immersions of the annulus $A$ instead of $D^2$ as described in section V and illustrated by the figures 6 and 7 . One would then  expect $W_2(B^3)$ to be achieved by a free boundary annulus of index at most $N_2=n_1+n_2=1+3=4$. The so called {\it critical catenoids} are then obvious candidates\footnote{See remark~\ref{crit-cat} which is making a more precise conjecture on this second level of the hierarchy for free boundary surfaces.} for realizing $W_2(B^3)$. Do we then have
 \[
 W_2(B^3)=2\,\pi\,\La^{-2}\quad ?\quad \mbox{ where }\quad\La\,\tanh\La=1\quad.
 \]
The space of oriented critical catenoids is diffeomorphic to two disjoint copies\footnote{Recall that
\[
H^\ast({\R}P^2,{\Z}_2)={\Z}_2[a]/a^3
\]
} of ${\R}P^2$ (the images by $SO(3)$ of a given critical catenoid with the two possible orientations). Exactly two values of $n_3$ can be then considered {\it a-priori}. Either $n_3=2$ or $n_3=3$, this would correspond
 to take the generator respectively of  $H^1({\R}P^2,{\Z}_2)$ and of $H^2({\R}P^2,{\Z}_2)$ or $n_3=1$, this last case would correspond to perform an eversion of the {\it critical catenoid} in $B^3$ preserving the boundary transverse to $\p B^3$. The two alternatives correspond to two branches of the {\it hierarchy tree}. The first alternatives $n_3=2,3$ should probably require to ``increase the topology'' as before by taking a genus zero surface with 3 boundary components while in the second case $n_3=1$
 one could restrict to  immersions of the annulus $A$ itself without having to increase the topology. This second case is of special interest because it could possibly produce a ``new'' free boundary annulus in $B^3$ with index $5$. For the remaining minmax of free boundary zero genus surfaces one should use non zero elements of
 \[
 H^\ast(SO(3),{\Z}_2)={\Z}_2[b]/b^4
 \]
hence we have $n_k\in\{1,2,3,4\}$.
\section{The general Definition of  Minmax Hierarchies.}
\reset
We shall denote ${\mathcal P}_N$ the category of $N-$dimensional compact manifolds orientable or non orientable, with or without boundary.


\subsection{The abstract scheme under the Palais Smale assumption.}

Let ${\mathfrak M}$ be an {\it Hilbert manifold} modeled on an Hilbert space ${\mathfrak H}$ or more generally a {\it Banach Manifold} modeled on a Banach space ${\mathfrak E}$  and assume that it is complete for the {\it Palais distance} $d_{\mathbf P}$ induced by the associated {\it Finsler Structure} $\|\cdot\|$ on $T\,{\mathfrak M}$. We shall restrict to the first case of
an Hilbert manifold \underbar{while considering exclusively} issues related to indices.

Let $E$ be a $C^2$ functional on ${\mathfrak M}$ which is assumed to satisfy the {\it Palais-Smale assumption} that is
\[
\forall\ \Phi_i\in {\mathfrak M}\quad\mbox{ s. t. }\quad\limsup_{i\rightarrow +\infty}E(\Phi_i)<+\infty\quad, \quad\lim_{i\rightarrow +\infty}\|DE(\Phi_i)\|_{\Phi_i}=0
\]
then
\[
\exists\, \Phi_{i'}\quad\mbox{ and }\quad \Phi_\infty\in {\mathfrak M}\quad\mbox{ s. t. }\quad d_{\mathbf P}(\Phi_{i'},\Phi_\infty)\longrightarrow 0\quad\mbox{ and }\quad DE(\Phi_\infty)=0
\]
A {\it Minmax Hierarchy} for the Functional $E$ requires the following objects

\begin{itemize}

\item[a)] a sequence (finite or infinite) of non zero integers $n_1,n_2,\cdots$
\item[b)] a sequence denoted 
\[
\mbox{Sweep}_{N_k}\subset\lf\{ (Y,\vec{\Phi})\ ;\ Y\in {\mathcal P}_{N_k}\quad\mbox{ and }\vec{\Phi}\in\mbox{Lip}_{\mathbf P}\lf(Y,{\mathfrak M}\rg) \rg\}\]
where $N_k:=n_1+\cdots+n_k$ and called {\it ${N_k}-$Sweepout space} of the hierarchy 
\item[c)] we have
\[
\forall\, (Y,\vec{\Phi})\in \mbox{Sweep}_{N_k}\quad\exists\, Z\in {\mathcal P}_{N_{k-1}}\quad\mbox{ s.t. }\ \p\,Y=\p( B^{n_k}\times Z)
\]
\item[d)] we have
\[
\forall\ x\in\p B^{n_k}\quad\quad \lf(Z,\vec{\Phi}(x,\cdot)\rg)\in \mbox{Sweep}_{N_{k-1}}
\]
\item[e)] a \underbar{strictly} increasing sequence of positive numbers (the {\it $k-$Widths} of the hierarchy)
\[
{W}_{k}:=\inf_{(Y,\vec{\Phi})\in \mbox{Sweep}_{N_k}}\quad\max_{y\in Y}\ E({\Phi}(y))
\]
\item[f)] for any homeomorphism $\Xi$ of $\mathfrak M$ satisfying
\[
\Xi({\Phi})={\Phi}\quad \mbox{ if }\quad E({\Phi})\le {W}_{k-1}-\eta_k
\]
for some $0<\eta_k<W_{k-1}-W_{k-2}$ we have
\[
\Xi\lf(\mbox{Sweep}_{N_k}\rg)\subset \mbox{Sweep}_{N_k}
\]
\end{itemize}

Using classical {\it Palais deformation theory} which applies to $E$  (see \cite{Riv-columb}) we obtain the existence of a sequence $\Phi_k$
such that
\[
E({\Phi}_k)={W}_k\quad.
\]
furthermore that the Banach Manifold ${\mathfrak M}$ is in fact {\it Hilbert}  and that $D^2E(\Phi)$ is Fredholm   we can ensure\footnote{We shall recall  the arguments leading to these assertions in the proof of 
theorem~\ref{th-I.1}} that $\mbox{Ind}({\Phi}_k)\le N_k$ where $\mbox{Ind}$ is the {\bf Morse Index} of $E$ (i.e. the maximal dimension of a vector space on which $D^2E$ is strictly negative).

\medskip

Of course the main issue in order to generate  such a hierarchy is to guarantee the series of strict inequalities between the successive $W_k$. We shall now explain a scheme that leads to such a series $W_k$.

 The first element in a hierarchy is an arbitrary {\it admissible familly} of the form
\[
\mbox{Sweep}_{n_1} \subset\lf\{
\begin{array}{l}
\ds(Y,\Phi)\in {\mathcal P}_{n_1}\times\mbox{Lip}(Y,{\mathfrak M})
\end{array}
\rg\}
\]
such that there exists $W_1>0$ satisfying
\[
\inf_{(Y,\Phi)\in\mbox{Sweep}_{n_1}}\max_{y\in Y}E(\Phi(y))= W_1
\]
and there exists $\eta>0$ such that for all homeomorphism $\Xi$ of ${\mathfrak M}$ equal to the identity for $E(\Phi)<W_1-\eta$ one has 
\[
\Xi(\mbox{Sweep}_{n_1})\subset \mbox{Sweep}_{n_1}\quad.
\]
Assuming now the hierarchy is constructed up to the order $k-1$, we introduce the notation for $l=1\cdots k-1$
\[
{\mathcal C}_l:=\lf\{ \Phi\in {\mathfrak M}\ ;\ \ ,\ E(\Phi)=W_{l}\quad\mbox{and }\quad DE(\Phi)=0\rg\}
\]
We are going to make the following assumption
\[
\mbox{(H1)}\quad\quad\quad {\mathcal C}_l \mbox{ \it is a smooth compact sub-manifold of }{\mathfrak M}
 \]
For any $\ep>0$ we denote
\[
{\mathcal O}_l(\ep):=\lf\{{\Phi}\in {\mathfrak M}\ ;\ d_{\mathbf P}({\Phi},{\mathcal C}_l)<\ep\rg\}
\]
Let $\ep_l>0$ be fixed such that $2\ep_l<\inf_{j<l}d_{\mathbf P}({\mathcal C}_l,{\mathcal C}_j)$ and
\[
\exists\ \pi_{l}\ \in\, \mbox{Lip}_{\mathbf P}\lf( {\mathcal O}_l(\ep_l) ,{\mathcal C}_l\rg)\quad\mbox{ s. t. }\quad \forall \ {\Phi}\in {\mathcal C}_l\quad\quad\pi_l({\Phi})={\Phi}
\]
as given by \cite{Lan}. The tubular neighborhood ${\mathcal O}_l(\ep_l) $ of ${\mathcal C}_l$  will simply be denoted ${\mathcal O}_l$. Because respectively of proposition~\ref{pr-A.1}  there exists $\delta_l>0$ such that
\be
\label{I.4}
\begin{array}{l}
\ds\forall\ (Y,\vec{\Phi})\in \mbox{Sweep}_{N_l}\quad\mbox{ , }\quad\max_{y\in Y}E({\Phi})\le {W}_l+\delta_l\\[3mm]
\ds\quad   \Longrightarrow\quad d_{\mathbf P}({\Phi}(Y),{\mathcal C}_l)<\ep_l\quad.
\end{array}
\ee
This being established we define $n_k$  as follows. Let $n_k\in{\N}^\ast$ such that\footnote{The reason why we are working
with ${\Z}_2$ cohomology comes from the fact that we are going to use Thom's resolution of Steenrod problem regarding the realization of ${\Z}_2-$homology classes
by continuous images of smooth manifolds (see \cite{Tho}).}
\[
H^{n_k-1}({\mathcal C}_{k-1},{\Z}_2)\ne 0
\]
and choose $\om_{k-1}$ being a non zero element of $H^{n_k-1}({\mathcal C}_{k-1},{\Z}_2)$. 

\medskip

Under the previous notations we define $\mbox{Sweep}_{N_k}$ to be the set of pairs $(Y,\vec{\Phi})$ such that

\begin{itemize}
\item[i)] 
\[
Y\in {\mathcal P}_{N_{k}}\quad,\quad {\Phi}\in\mbox{Lip}_{\mathbf P}(\ov{Y},{\mathfrak M})\quad.
\]
\item[ii)] There exists $Z\in {\mathcal P}_{N_{k-1}}$ s.t.
\[
\p Y=\p\lf(B^{n_k}\times Z\rg)
\]
\item[iii)] We have
\[
\forall\ x\in\p B^{n_k}\quad\quad \lf(Z,\vec{\Phi}(x,\cdot)\rg)\in \mbox{Sweep}_{N_{k-1}}
\]
and
\[
\max_{y\in \p Y}E({\Phi}(y))\le {W}_{k-1}+2^{-1}\,\delta_{k-1}
\]
\item[iv)] Let 
\[
\Om_{{\Phi}}:=\lf\{  y\in \p Y\ ;\ d_{\mathbf P}({\Phi}(y),{\mathcal C}_{k-1})<\ep_{k-1} \rg\}
\]
we have\footnote{Observe that $[\Om_{{\Phi}}\cap (\{x\}\times Z)]\in H_{N_{k-1}}(\Om_{{\Phi}},{\Z}_2)$ is independent of $x\in \p B^{n_k}$ for $n_k>1$.}

\[
[\Om_{{\Phi}}\cap (\{x\}\times Z)]\in H_{N_{k-1}}(\,\Om_{{\Phi}},\p\,\Om_{{\Phi}} ,{\Z}_2)\quad\mbox{ is {\bf Poincar\'e dual} to } (\pi_{k-1}\circ{\Phi})^\ast\om_{k-1}\in H^{n_k-1}(\Om_{{\Phi}},{\Z}_2)
\]
\item[v)]  We have
\[
\max_{y\in  B^{n_k}\times \p Z} E({\Phi}(y))\le\ {W}_{k-1}-\delta_{k-1}\quad.
\]
\end{itemize}

\medskip

\begin{figure}
\begin{center}
\includegraphics[width=17cm,height=23cm]{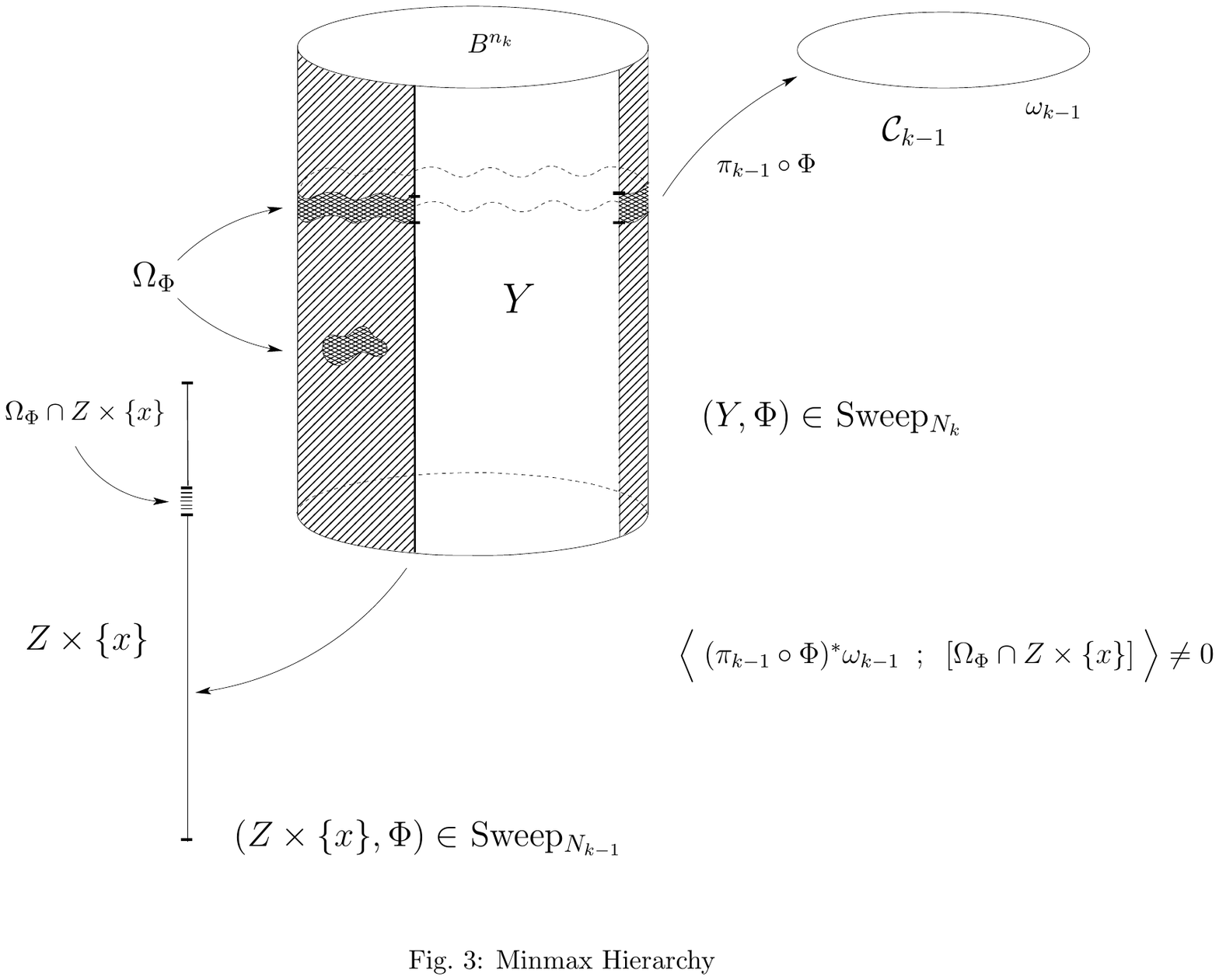}
\end{center}
\end{figure}

\begin{Th}
\label{th-I.1}
Under the hypothesis (H1) and assumptions i)...v) we have
\be
\label{strict-sigma}
{W}_{k-1}+\delta_{k-1}\le {W}_k
\ee
and $\{\mbox{Sweep}_{N_l}\}_{l\le k}$ defines a {\bf Minmax Hierarchy}. Consequently, for any $l\le k$ there exists ${\Phi}_l$ such that
\[
E({\Phi}_l)={W}_l\quad,\quad DE({\Phi}_l)=0\quad
\]
Assuming furthermore that the Banach Manifold ${\mathfrak M}$ is in fact {\it Hilbert}  and that $D^2E(\Phi)$ is Fredholm for any $\Phi$, we have
\be
\label{indice-palais}
\mbox{Ind}\,({\Phi}_l)\le N_l
\ee
where $\mbox{Ind}$ is the {\bf Morse Index} of $E$ (i.e. the maximal dimension of a vector space on which $D^2E$ is strictly negative).

\hfill $\Box$
\end{Th}
\noindent{\bf Proof of theorem~\ref{th-I.1}.} We shall first prove by induction  that ${W}_{k-1}+\delta_{k-1}\le {W}_k$. Assume this is not the case but assume $W_{l-1}+\delta_{l-1}\le W_l$ for $l\le k-1$. Choose $(Y,\vec{\Phi})\in \mbox{Sweep}_{N_k}$
such that
\be
\label{I.6}
\max_{y\in Y}\, E({\Phi}(y))\le {W}_{k-1}+\delta_{k-1}
\ee
Consider 
\[
\La_{{\Phi}}=\lf\{  y\in Y\ ;\ d_{\mathbf P}({\Phi}(y),{\mathcal C}_{k-1})<\ep_{k-1} \rg\}\quad.
\]
and denote $\Xi_{{\Phi}}:=\p \La_{{\Phi}}\cap\mbox{int}(Y)$. Since $y\rightarrow d_{\mathbf P}({\Phi}(y),{\mathcal C}_{k-1})$ is lipschitz we can assume, without loss of generality, that $\Xi_{{\Phi}}$ realizes an {\it homology manifold} which makes $(\Xi_{{\Phi}},\p \Xi_{{\Phi}})$ a {\it Poincar\'e duality pair}  .

 For $\ep_{k-1}$ chosen small enough we have
\[
F(\Phi)\le W_{k-1}-\delta_{k-1}\quad\Longrightarrow\quad d_{\mathbf P}(\Phi,{\mathcal C}_{k-1})>\ep_{k-1}
\]
thus v) implies
\[
\Xi_{{\Phi}}\cap (\p B^{n_k}\times Z)=\emptyset
\]
hence
\[
\p\, \Xi_{{\Phi}}=\p\,\Om_{{\Phi}}
\]
We shall now prove the following.

\medskip

\noindent{\bf Claim 1}
\[
[\p\, \Om_{{\Phi}}\cap (\{x\}\times Z)]\in \mbox{Im}\,\p_\ast \quad,
\]
where $\p_\ast$ is the boundary operator 
\[
\p_\ast\ :\ H_{N_{k-1}}(\Xi_{{\Phi}},\p\,\Xi_{{\Phi}},{\Z}_2)\quad\longrightarrow\quad H_{N_{k-1}-1}(\p\, \Xi_{{\Phi}},{\Z}_2)
\]
We recall the following {\it relative Poincar\'e duality} commutative diagram
\[
\begin{CD}
H^{p-1}(\p\, \Xi_{{\Phi}},{\Z}_2) @>\delta^\ast>> H^p(\Xi_{{\Phi}},\p\, \Xi_{{\Phi}},{\Z}_2)@>j^\ast>> H^p(\Xi_{{\Phi}},{\Z}_2)@>i^\ast>>H^p(\p\,\Xi_{{\Phi}},{\Z}_2)\\
@VVDV    @VVDV  @VVDV  @VVDV\\
H_{N_{k}-1-p}(\p\, \Xi_{{\Phi}},{\Z}_2)        @>i_\ast>> H_{N_k-p-1}(\Xi_{{\Phi}},{\Z}_2)  @>j_\ast >> H_{N_k-p-1}(\Xi_{{\Phi}},\p \,\Xi_{{\Phi}},{\Z}_2)@>\p_\ast >> H_{N_k-p-2}(\p\, \Xi_{{\Phi}},{\Z}_2)
\end{CD}
\]
where the vertical arrows are Poincar\'e isomorphisms simply denoted by $D$ and given, modulo a sign, by cap products respectively with $[\Xi_{{\Phi}}]$ or $[\p \,\Xi_{{\Phi}}]$. We apply this diagram
to the case $p=n_k-1$. The map $\pi_{k-1}\circ{\Phi}$ is well defined on $\Xi_{{\Phi}}$ and
\[
(\pi_{k-1}\circ{\Phi})^\ast\om_{k-1}\in H^{n_k-1}(\Xi_{{\Phi}},{\Z}_2)
\]
It is clear that the image of $(\pi_{k-1}\circ{\Phi})^\ast\om_{k-1}$ by the restriction map $i^\ast$ is $(\pi_{k-1}\circ{\Phi})^\ast\om_{k-1}$
itself where ${\Phi}$ is restricted to $Z\times \p B^{n_k}$. Let $U\in H_{N_{k-1}}(\Xi_{{\Phi}},\p\,\Xi_{{\Phi}},{\Z}_2)$ be the Poincar\'e dual to
$(\pi_{k-1}\circ{\Phi})^\ast\om_{k-1}\in H^{n_k-1}(\Xi_{{\Phi}},{\Z}_2)$. Because of the last part of the above diagram we have that
\[
\p_\ast U\quad\mbox{ is the Poincar\'e dual to }(\pi_{k-1}\circ{\Phi})^\ast\om_{k-1}\in H^{n_k-1}(\p\,\Xi_{{\Phi}},{\Z}_2)
\]
We shall now make use of the following lemma
\begin{Lm}
\label{rm-I.3} The assumption
\[
[\Om_{{\Phi}}\cap (\{x\}\times Z)]\in H_{N_{k-1}}(\Om_{{\Phi}},\p\Om_{{\Phi}},{\Z}_2)\quad\mbox{ is {\bf Poincar\'e dual} to } (\pi_{k-1}\circ{\Phi})^\ast\om_{k-1}\in H^{n_k-1}(\Om_{{\Phi}},{\Z}_2)
\]
implies 
\[
[\p\,\Om_{{\Phi}}\cap (\{x\}\times Z)]\in H_{N_{k-1}-1}(\p\,\Om_{{\Phi}},{\Z}_2)\quad\mbox{ is {\bf Poincar\'e dual} to } (\pi_{k-1}\circ{\Phi})^\ast\om_{k-1}\in H^{n_k-1}(\p\Om_{{\Phi}},{\Z}_2)
\]
\hfill$\Box$
\end{Lm}
\noindent{\bf Proof of lemma~\ref{rm-I.3} .} The lemma is a direct consequence of the following {\it relative Poincar\'e duality} commutative diagram 
\[
\begin{CD}
 H^{n_k-1}(\Om_{{\Phi}},{\Z}_2)@>i^\ast>>H^{n_k-1}(\p\Om_{{\Phi}},{\Z}_2)\\
  @VVDV  @VVDV\\
 H_{N_{k-1}}(\Om_{{\Phi}},\p \Om_{{\Phi}},{\Z}_2)@>\p_\ast >> H_{N_{k-1}-1}(\p \Om_{{\Phi}},{\Z}_2)
\end{CD}
\]
\hfill $\Box$

\noindent{\bf End of the proof of claim 1}. Hence, since by the previous lemma $[\p\, \Om_{{\Phi}}\cap (\{x\}\times Z)]$ is also {\it Poincar\'e dual} to $(\pi_{k-1}\circ{\Phi})^\ast\om_{k-1}\in H^{n_k-1}(\p\,\Xi_{{\Phi}},{\Z}_2)$, by uniqueness of the {\it Poincar\'e dual} we have
\[
[\p\,\Om_{{\Phi}}\cap (\{x\}\times Z)]=\p_\ast U
\]
 and  $[\p\,\Om_{{\Phi}}\cap (\{x\}\times Z)]$ is a boundary in $\Xi_{{\Phi}}$ and the claim is proved.  Using Thom's proof of {\it Steenrod Problem} on the realization of ${\Z}_2$ homology classes by continuous images of smooth un-oriented manifolds (see theorem III.2 in \cite{Tho}), we can assume that the concrete chain $U$ in $Y$ is the image  of an element in ${\mathcal P}_{N_{k-1}}$.  By an abuse of notation we identify $U$ with this element in ${\mathcal P}_{N_{k-1}}$. By ''pushing''  $U$  inside $Y\setminus \La_{{\Phi}}$, and summing this  with the homology manifolds $Z\times \{x\}\setminus \Om_{{\Phi}}\cap (\{x\}\times Z)$ we obtain and $(V,{\Phi})\in \mbox{Sweep}_{N_{k-1}}$. Because of (\ref{I.6}) we have 
\[
\max_{y\in V} E({\Phi}(y))\le {W}_{k-1}+\delta_{k-1}
\]
Using proposition~\ref{pr-A.1} we then have the existence of $y_\ep\in V$ such that $d_{\mathbf P}({\Phi}(y),{\mathcal C}_{k-1})<\ep_{k-1}$ which is a contradiction. Hence we have
\[
{W}_{k-1}+\delta_{k-1}\le {W}_k
\]
Consider a Pseudo-gradient for $E$ on ${\mathfrak M}^\ast:={\mathfrak M}\setminus \{\Phi\ ;\ DE(\Phi)=0\}$. We choose a cut-off for the action of the  Pseudo-gradient above the energy levels
${W}_{k-1}+\delta_{k-1}^\sigma/2$ in order for the flow to preserving the membership in $\mbox{Sweep}_{N_k}$. Following the classical {\it Palais deformation} arguments (see for instance \cite{Riv-columb}) we deduce the existence of $\Phi_k$ such that $E(\Phi_k)=W_k$ and $DE(\Phi_k)=0$.

Assuming now $\mathfrak M$ defines in fact an {\it Hilbert manifold} modeled on an Hilbert space ${\mathfrak H}$ and that $D^2E$ is everywhere {\it Fredholm}. Denote
\[
\p \,\mbox{Sweep}_{N_k}:=\lf\{
\begin{array}{l}
\ds (Z,\Psi)\in {\mathcal P}_{N_{k-1}}\times\mbox{Lip}(\p(B^{n_k}\times Z),{\mathfrak M})\quad;\quad\exists\, (Y,\Phi)\in \mbox{Sweep}_{N_k}\\[3mm]
\ds \p Y=\p(B^{n_k}\times Z)\quad;\quad\Phi=\Psi\quad\mbox{ on }\p Y
\end{array}
\rg\}
\]
By definition we have
\[
\mbox{Sweep}_{N_k}\subset\lf\{ C\in {\mathfrak C}_{N_k}\quad\exists (Z,\Psi)\in \p \,\mbox{Sweep}_{N_k}\quad \p C=\Psi_\ast[\p (B^{n_k}\times Z)]\rg\}
\]
where ${\mathfrak C}_{N_k}$ is the space of $N_k-$polyhedral chains in ${\mathfrak M}$. Observe that $[C]\in H_{N_k}({\mathfrak M}, \p (B^{n_k}\times Z))$ is non trivial. Indeed, assume there exists $D\subset B$ such that $\p D=\p C$ we would have $W_{k}<W_{k-1}+2^{-1} \delta_{k-1}$ which contradicts (\ref{strict-sigma}) . Hence $\mbox{Sweep}_{N_k}$ is by definition an {\it homological family} of dimension $N_k$ with boundary the cycles $\Psi_\ast[\p (B^{n_k}\times Z)]$ (see \cite{Gho}). Using corollary 10.5 of \cite{Gho} we obtain (\ref{indice-palais}).  This concludes
the proof of theorem~\ref{th-I.1}. \hfill $\Box$

\medskip

We observe that the topological condition regarding ${\Phi}$ on $\Om_{{\Phi}}$ is preserved by enlarging the set. Precisely the following lemma holds.

\begin{Lm}
\label{lm-I.4}
Let $V\subset Z\times\p B^{n_k}$ such that $\Om_{{\Phi}}\subset V$ and such that $\pi_{k-1}\circ{\Phi}$ extends continuously on $V$. Let $\om\in H^{n_{k-1}}(V)$ given by $\om:=j^\ast(\pi_{k-1}\circ{\Phi})^\ast\om_{k-1}$ where $j$ is the canonical inclusion map $j\, :\, \Om_{{\Phi}}\rightarrow V$. Assume
\[
[V\cap (\{x\}\times Z)]\in H_{N_{k-1}}(V,\p V,{\Z}_2)\quad\mbox{ is {\bf Poincar\'e dual} to } \om\in H^{n_k-1}(V,{\Z}_2)
\]
for some $x\in \p B^{n_k}$, then the condition iv) is satisfied.\hfill $\Box$
\end{Lm}
\noindent{\bf Proof of lemma~\ref{lm-I.4}.} The inclusion map $j$ induces a map on relative $N_{k-1}-$chains as follows
\[
j\,:\, C_{N_{k-1}}(V,\p V,{\Z}_2)\,\rightarrow\,C_{N_{k-1}}(\Om_{{\Phi}},\p\Om_{{\Phi}},{\Z}_2)\quad.
\] 
by restriction to $\Om_{\vec{\Phi}}$. We then have the restriction operator
\[
j_\ast\,:\, H_{N_{k-1}}(V,\p V,{\Z}_2)\,\rightarrow\,H_{N_{k-1}}(\Om_{{\Phi}},\p\Om_{{\Phi}},{\Z}_2)\quad.
\]
For any $\al\in H^{N_{k-1}}(\Om_{{\Phi}},\p\Om_{{\Phi}},{\Z})$ we have for the {\it cup product} $(\pi_{k-1}\circ{\Phi})^\ast\om_{k-1}\smile\al\in H^{N_{k-1}+n_k-1}(\Om_{{\Phi}},\p\Om_{{\Phi}},{\Z}_2)$
\[
\begin{array}{l}
\ds\lf<(\pi_{k-1}\circ{\Phi})^\ast\om_{k-1}\smile\al,\lf[\Om_\Phi\rg]\rg>=\lf<(\pi_{k-1}\circ{\Phi})^\ast\om_{k-1}\smile\al,j_\ast\lf[V\rg]\rg>\\[3mm]
\ds\quad\quad=\lf<j^\ast\lf((\pi_{k-1}\circ{\Phi})^\ast\om_{k-1}\smile\al\rg),\lf[V\rg]\rg>=\lf<\om\smile j^\ast\al,\lf[V\rg]\rg>=\lf<j^\ast\al,\om\frown\lf[V\rg]\rg>
\end{array}
\]
We are assuming that $\om$ is Dual to $[V\cap (\{x\}\times Z)]\in H_{N_{k-1}}(V,\p V,{\Z}_2)$ in other words
\[
[V\cap (\{x\}\times Z)]=\om\frown\lf[V\rg]
\]
Hence we have proved that $\forall\,\al\in H^{N_{k-1}}(\Om_{{\Phi}},\p\Om_{{\Phi}},{\Z})$
\[
\lf<(\pi_{k-1}\circ{\Phi})^\ast\om_{k-1}\smile\al,\lf[\Om_\Phi\rg]\rg>=\lf<j^\ast\al,[V\cap (\{x\}\times Z)]\rg>=\lf<\al,j_\ast[V\cap (\{x\}\times Z)]\rg>
=\lf<\al,[\Om_{{\Phi}}\cap (\{x\}\times Z)]\rg>
\]
Which implies
\[
[\Om_{{\Phi}}\cap (\{x\}\times Z)]=(\pi_{k-1}\circ{\Phi})^\ast\om_{k-1}\frown[\Om_{{\Phi}}]
\] 
Hence $[\Om_{{\Phi}}\cap (\{x\}\times Z)]$ is dual to $(\pi_{k-1}\circ{\Phi})^\ast\om_{k-1}$ and the lemma is proved.\hfill $\Box$

\begin{Rm}
\label{rm-I.3p}
For  $n_1=1$ one can afford to restrict to $(Y,\vec{\Phi})\in\mbox{Sweep}_1$ where $Y=(-1,+1)$ moreover one can replace $H^{n_2-1}({\mathcal C}_1,{\Z}_2)$ by $H^{n_2-1}({\mathcal C}_1,{\Z})$ (which is ``richer'') in the definition of $\mbox{Sweep}_{N^2}$. This is due to the fact that the chain $U$ in the proof of theorem~\ref{th-I.1} can be taken to be a segment homeomorphic to $(-1,+1)$ and that there is no orientation problem at this first level of the hierarchy in the case $n_1=1$. The whole proof
in this case for the passage from $k=1$ to $k=2$ is transposable word by word by replacing ${\Z}_2$ by ${\Z}$. 
\hfill$\Box$
\end{Rm}
\subsection{Minmax Hierarchies in the Linear case of Laplace Eigenvalues.}
\label{rayleigh}
The purpose of this section is to express the {\it Rayleigh quotient} for computing the eigenvalues of the laplacian on a closed surface as a particular case of the {\it Minmax hierarchies}.
We adopt the notations from the introduction and we define the {\it Hierarchy} by i), ii), iii) and iv). It is clear that it is a particular case of the hierarchies defined in the previous section
and hence, because of theorem~\ref{th-I.1} we have a strictly increasing sequence of eigenvalues $\mu_k$ such that
\[
\mu_k:=\inf_{u\in \mbox{Sweep}_{N_k}}\max_{y\in Y_k} E(u_y)
\]
The problem is to show that we indeed ``capture'' all the successive eigenvalues of the Laplacian. It suffices to show then that
\be
\label{eigenv}
\mu_k\le\la_k\quad.
\ee
We prove (\ref{eigenv}) by induction. Let $F_l$ be the vector eigenspace associated to $\la_l$. By induction assumption we have $n_l=\mbox{dim} F_l=\mbox{dim}\,{\mathcal C}_l+1$
We take for $\om_l$ the generator of $H^{n_l-1}({\mathcal C}_l,{\Z}_2)\simeq{\Z}_2$ (we obviously have ${\mathcal C}_l\simeq S^{n_l-1}$. Let $u_1\cdots u_k$ be an arbitrary
choice of eigenfunctions in ${\mathfrak S}$ for $\la_1<\cdots<\la_k$. We denote by $SO(F_l)$ the space of positive isometries of the euclidian spaces $F_l$.

\medskip

Let $u\ :\ (-1,+1)^k\times SO(F_2)\times\cdots SO(F_{k-1})\longrightarrow {\mathfrak S}$ given by
\[
\begin{array}{l}
u(t_1\cdots t_k,R_2\cdots R_{k-1}):=\\[3mm]
\ds\lf\{
\cos(\frac{\pi t_k}{2})\, u_k+\sin(\frac{\pi t_k}{2})\lf(\cos(\frac{\pi t_{k-1}}{2})\, R_{k-1}\,u_{k-1}+\sin(\frac{\pi t_{k-1}}{2})\,\lf(\cdots\lf(\cos(\frac{\pi t_2}{2})\, R_{2}\,u_{2}+\sin(\frac{\pi t_2}{2})\, u_1\rg)\rg)\rg)
\rg\}
\end{array}
\]
It is straightforward to check that, assuming $\mu_l=\la_l$ for $l<k$ we have $u\in \mbox{Sweep}_{N_k}$ and moreover
\[
\max_{(t_1\cdots t_k,R_2\cdots R_{k-1})\in (-1,+1)^k\times SO(F_2)\times\cdots SO(F_{k-1})}E\lf(u(t_1\cdots t_k,R_2\cdots R_k)\rg)\le \la_k
\]
This proves (\ref{eigenv}).

\section{Index Control in the Viscosity Method.}
\reset
Let $M^m$ be a closed sub-manifold of an euclidian space ${\R}^Q$.
We denote by $\Sigma^g$ a closed oriented 2 dimensional manifold of genus $g$.
\[
\mbox{Imm}(\Sigma^g,M^m):=\lf\{
\begin{array}{c} \vec{\Phi}\in W^{2,4}(\Sigma^g,M^m)\ ;\ \mbox{rank}\,(d\vec{\Phi}_x)=2\quad\forall\, x\in \Sigma^g
\end{array}
\rg\}
\]
We equip $\mbox{Imm}(\Sigma^g,M^m)$ with the $W^{3,2}$ topology  which makes it a {\it Hilbert manifold} (see \cite{Riv-columb}).
From \cite{Pink} we know that every pair $\vec{\Phi}_1$ and $\vec{\Phi}_2$ of embeddings of $\Sigma^g$  in $\mbox{Imm}(\Sigma^g,M^m)$ are regular homotopic : there exists a diffeomorphism $\Psi$ of $\Sigma^g$ such that $\vec{\Phi}_1$ and $\vec{\Phi}_2\circ\Psi$ are in the same path connected
component of $\mbox{Imm}(\Sigma^g,M^m)$.
We denote by $\mbox{Imm}_0(\Sigma^g,M^m)$ the subspace of $\mbox{Imm}(\Sigma^g,M^m)$ in the same path connected component of an embedding. 


\medskip

We shall denote 
\[
\mbox{Imm}(M^m):=\bigcup_{g\in {\N}}\mbox{Imm}_0(\Sigma^g,M^m)\quad\mbox{ and }\quad\mbox{Imm}^{g_0}(M^m):=\bigcup_{g\le g_0}\mbox{Imm}_0(\Sigma^g,M^m)
\]
For any element $\vec{\Phi}$ in $\mbox{Imm}_0(M^m)$ we denote by $\Sigma_{\vec{\Phi}}$ the abstract surface such that $\vec{\Phi}\in \mbox{Imm}(\Sigma_{\vec{\Phi}},M^m)$.

\medskip

\subsection{The Bundle of $W^{3,2}-$immersions of surfaces over the Hilbert Manifold of Immersed $W^{3,2}-$Surfaces.}

Let $\mbox{Diff}_+(\Sigma^g)$ be the topological group of positive  $W^{3,2}-$diffeomorphisms of $\Sigma^g$,  isotopic to the identity. This can be seen as an open subspace of $W^{3,2}(\Sigma^g,\Sigma^g)$ which itself defines a Hilbert Manifold (see \cite{Riv-columb}). For $g=0$ we are marking 3 distinct points, that we denote $a_1,a_2,a_3$, for $g=1$ we are marking one point  that we denote $a$ and for $g>1$ no point is marked. We denote by $\mbox{Diff}^\ast_+(\Sigma^g)$ the sub-group of $\mbox{Diff}_+(\Sigma^g)$ which are fixing the marked points. In particular for $g>1$ we have $\mbox{Diff}^\ast_+(\Sigma^g)=\mbox{Diff}_+(\Sigma^g)$. Recall nevertheless that for any diffeomorphism $\Psi$ isotopic to the identity the {\it Lefschetz Number} $L(\Psi)$ is given by definition by
\[
L(\Psi)=\mbox{Tr}(\Psi| H_0(\Sigma^g))-\mbox{Tr}(\Psi| H_1(\Sigma^g))+\mbox{Tr}(\Psi| H_2(\Sigma^g))=2-2g
\]
Hence for $g>1$ we have $L(\Psi)\ne 0$ and then $\Psi$ must have at least one fixed point. Due to lemma 1.2 in \cite{CMM} we deduce that for any $g\in {\N}$ the action of $\mbox{Diff}^{\, \ast}_+(\Sigma^g)$ on ${\mathcal M}^g:=\mbox{Imm}_0(\Sigma^g,M^m)$ is {\it free}. We aim to study the differential structure of this quotient. Since $\mbox{Diff}^{\, \ast}_+(\Sigma^g)$ misses
to be a {\it Banach Lie Group} but is only a {\it topological Group}\footnote{On the space of $W^{3,2}-$diffeomorphisms the right multiplication is smooth but the left multiplication is not differentiable, the inverse mapping is not $C^1$, there is no canonical chart in the neighborhood of the identity, the exponential map is continuous but not $C^1$, it
is not locally surjective in a neighborhood of the identity, the Bracket operation in the Tangent space to the identity is not continuous....etc see a description of all these ``pathological behavior'' for instance in \cite{EbMa} or \cite{Mil}} with an {\it Hilbert Manifold} structure the existence of a differentiable Hilbert structure on the quotient $\mbox{Imm}_0(\Sigma^g,M^m)/\mbox{Diff}^{\, \ast}_+(\Sigma^g)$ is not the result of classical
consideration and deserves to be studied with care (Progresses in this direction are given in \cite{BiFi} for $W^{3,2}-$embeddings but we are not going to follow this approach and the one
we choose is more specific to surfaces but more precise too) This is the goal of the next result : theorem~\ref{th-hilbert}.

\medskip

 For any $\vec{\Phi}\in {\mathcal M}^g:=\mbox{Imm}_0(\Sigma^g,M^m)$ we denote $P_{\vec{\Phi}}$ the $L^2-$ orthogonal projection from $(\wedge^{(0,1)}\Sigma)^{\otimes^2}$ onto the space  of {\it holomorphic quadratic form} $\mbox{Hol}_Q(\Sigma^g,g_{\vec{\Phi}}))$ on $(\Sigma^g,g_{\vec{\Phi}})$  and by $P_{\vec{\Phi}}^\perp:=\mbox{Id}-P_{\vec{\Phi}}$. Define the linear map
\[
\begin{array}{rcl}
\ds \ov{D}^\ast_{\vec{\Phi}}\ :\  T_{\vec{\Phi}}{\mathcal M}^g &\longrightarrow &\ds W^{2,2}\lf(   (\wedge^{(0,1)}\Sigma)^{\otimes^2} \rg)\\[5mm]
\ds \vec{w} &\longrightarrow & P_{\vec{\Phi}}^\perp\lf(\ov{\p}\vec{w}\ \dot{\otimes}\ \ov{\p}\vec{\Phi}\rg)
\end{array}
\]
where, in local complex coordinates for $\vec{\Phi}$ we denote
\[
 \ov{\p}\vec{w}\ \dot{\otimes}\ \ov{\p}\vec{\Phi}:={\p}_{\ov{z}}\vec{w}\, \cdot\, {\p}_{\ov{z}}\vec{\Phi}\ \ d\ov{z}\otimes d\ov{z}
\]
We are now going to prove the following theorem 

\begin{Th}
\label{th-hilbert}
Let $\vec{\Phi}\in\mbox{Imm}_0(\Sigma^g,M^m)$, then there exists an open neighborhood ${\mathcal O}_{\vec{\Phi}}$ of $\vec{\Phi}$ in $\mbox{Imm}_0(\Sigma^g,M^m)$ invariant under the action of $\mbox{Diff}^{\, \ast}_+(\Sigma^g)$ and
two smooth maps on ${\mathcal O}_{\vec{\Phi}}$, equivariant under the action of $\mbox{Diff}^{\, \ast}_+(\Sigma^g)$,
\[
{\vec{w}}_{\vec{\Phi}}\ :\ {\mathcal O}_{\vec{\Phi}}\longrightarrow \mbox{Ker}( \ov{D}^\ast_{\vec{\Phi}})\subset T_{\vec{\Phi}}{\mathcal M}^g 
\]
and
\[
{\Psi}_{\vec{\Phi}}\ :\ {\mathcal O}_{\vec{\Phi}}\longrightarrow  \mbox{Diff}^{\, \ast}_+(\Sigma^g)
\]
satisfying 
\[
\forall\ \vec{\Xi}\in {\mathcal O}_{\vec{\Phi}}\quad\quad\vec{\Xi}\circ\Psi_{\vec{\Phi}}(\vec{\Xi})=\pi_{M^m}\lf(\vec{\Phi}+\vec{w}_{\vec{\Phi}}(\vec{\Xi})\rg)
\]
where $\pi_{M^m}$ is the orthogonal projection onto $M^m$ and $\Upsilon_{\vec{\Phi}}:=({\vec{w}}_{\vec{\Phi}},{\Psi}_{\vec{\Phi}})$ realizes a diffeomorphism from ${\mathcal O}_{\vec{\Phi}}$ onto  $U_{\vec{\Phi}}\times \mbox{Diff}^{\, \ast}_+(\Sigma^g)$ where $U_{\vec{\Phi}}$ is a neighborhood of  $0$ in $\mbox{Ker}( \ov{D}^\ast_{\vec{\Phi}})$. Moreover the map $\Upsilon_{\vec{\Phi}}$ satisfies the following equivariance property : $\forall \, \vec{\Xi}\in {\mathcal O}_{\vec{\Phi}}$ and forall $\Psi_0\in \mbox{Diff}^{\, \ast}_+(\Sigma^g)$ one has
\[
\Psi_{\vec{\Phi}}(\vec{\Xi}\circ{\Psi}_0)=\Psi_0^{-1}\circ\Psi_{\vec{\Phi}}(\vec{\Xi})\quad\mbox{ and }\quad\vec{w}_{\vec{\Phi}}(\vec{\Xi}\circ\Psi_0)=\vec{w}_{\vec{\Phi}}(\vec{\Xi})\quad.
\]

The space $\mbox{Imm}_0(\Sigma^g,M^m)/\mbox{Diff}^{\, \ast}_+(\Sigma^g)$ is Hausdorff and defines a Hilbert Manifold such that the projection map
\[
\mbox{Imm}_0(\Sigma^g,M^m)\longrightarrow{\mathfrak M}_{g}(M^m):=\mbox{Imm}_0(\Sigma^g,M^m)/\mbox{Diff}^{\, \ast}_+(\Sigma^g)
\]
defines a $\mbox{Diff}^{\, \ast}_+(\Sigma^g)-$bundle for which $(\Upsilon_{\vec{\Phi}})_{\vec{\Phi}}$ represents a local trivialization. \hfill $\Box$\end{Th}

\noindent{\bf Proof of theorem~\ref{th-hilbert}.} We first construct $\Upsilon_{\vec{\Phi}}$ in a neighborhood of $\vec{\Phi}$. A basis of neighborhoods  of $\vec{\Phi}$ is given by
\[
{\mathcal V}^\ep_{\vec{\Phi}}:=\lf\{\vec{\Xi}=\pi_{M^m}\lf(\vec{\Phi}+\vec{v}\rg)\quad;\quad \vec{v}\in \Gamma^{3,2}(\vec{\Phi}^\ast TM^m)\cap W^{3,2}(\Sigma,{\R}^Q)\ \mbox{ and }\ \|\vec{v}\|_{W^{3,2}}<\ep\rg\}\quad.
\]
for $\ep>0$ small enough and where $\Gamma^{3,2}(\vec{\Phi}^\ast TM^m)$ denotes the $W^{3,2}-$sections of the pull-back bundle $\vec{\Phi}^\ast TM^m$, that is the sub-vector space
of $\vec{v}\in W^{3,2}(\Sigma^g,{\R}^Q)$ such that $\vec{v}(x)\in T_{\vec{\Phi}(x)}\Sigma^g$ for any $x\in \Sigma^g$.

For any  $\vec{v}\in \Gamma^{3,2}(\vec{\Phi}^\ast TM^m)$ we consider the tensor in $\Gamma^{2,2}((T^\ast\Sigma^g)^{(0,1)}\otimes(T\Sigma^g)^{(1,0)})$ given by
\[
  \ov{D}^\ast_{\vec{\Phi}}\vec{v}\ \res\, g^{-1}_{\vec{\Phi}} \quad\mbox{ where }\quad g_{\vec{\Phi}}^{-1}=e^{-2\la}\ [\p_z\otimes\p_{\ov{z}}+\p_{\ov{z}}\otimes \p_z]
\]
and where $\res$ is the contraction operator between covariant and contravariant tensors. We denote
\be
\label{mathcali}
{\mathcal I}:=\lf\{ \ov{D}^\ast_{\vec{\Phi}}\vec{v}\ \res\, g^{-1}_{\vec{\Phi}} \ ;\quad \vec{v}\in \Gamma^{3,2}(\vec{\Phi}^\ast TM^m)\rg\}\quad.
\ee
Recall the definition of the $\ov{\p}$ operator on $\wedge^{(1,0)}\Sigma^g$ given in local coordinates by
\[
\ov{\p}\lf(a\,\p_z\rg)=\p_{\ov{z}}a\ d\ov{z}\otimes\p_z\quad.
\]
Denote $\mbox{Hol}_1(\Sigma^g)$ the finite dimensional subspace of $\Gamma^{3,2}(\wedge^{(1,0)}\Sigma^g)$ made of holomorphic sections\footnote{Due to Riemann-Hurwitz Theorem, the  holomorphic tangent bundle $T^{(1,0)}\Sigma^g$ which is the inverse of the canonical bundle of the Riemann surface defined by $(\Sigma,g_{\vec{\Phi}})$ has a degree given by
\[
\mbox{deg}\lf(T^{(1,0)}\Sigma^g\rg)=\int_{\Sigma^g}c_1\lf(T^{(1,0)}\Sigma^g\rg)=2-2g
\]
 and therefore $\mbox{Hol}_1(\Sigma^g)\ne 0$ $\Rightarrow$ $g<2$.} of $T^{(1,0)}\Sigma^g$. 
We shall now prove the following lemma.

\begin{Lm}
\label{lm-delbar}
Under the previous notations we have that $\forall\ \vec{v}\in W^{3,2}(\Sigma^g,{\R}^Q)$
\be
\label{delbar}
\forall\ \vec{v}\in W^{3,2}(\Sigma,{\R}^Q)\quad \exists\, !\, f\in (\mbox{Hol}_1(\Sigma^g))^\perp\cap \Gamma^{3,2}(\wedge^{(1,0)}\Sigma^g)\quad\quad\mbox{s. t. }\quad\ov{\p}f=\ov{D}^\ast_{\vec{\Phi}}\vec{v}\ \res\, g^{-1}_{\vec{\Phi}}\quad.
\ee
Moreover
\be
\label{delbar-estim}
\|f\|_{W^{3,2}}\le C_{\vec{\Phi}}\ \|\vec{v}\|_{W^{3,2}}\quad.
\ee
\hfill $\Box$
\end{Lm}
\noindent{\bf Proof of lemma~\ref{lm-delbar}.}
We have for any $\al=a\, \p_z\in \Gamma^{3,2}(\wedge^{(1,0)}\Sigma^g)$ and $\beta=b\ d\ov{z}\otimes\p_z\in \Gamma^{2,2}((T^\ast\Sigma^g)^{(0,1)}\otimes(T\Sigma^g)^{(1,0)})$
\[
\begin{array}{l}
\ds\int_{\Sigma^g}\lf<\ov{\p}\lf(a\,\p_z\rg), b\ d\ov{z}\otimes \p_z\rg>_{g_{\vec{\Phi}}}\ d\mbox{vol}_{g_{\vec{\Phi}}}=\Re\lf[\frac{i}{2}\int_{\Sigma^g} {\p}\ov{a}\  b\  e^{2\la}\, dz\wedge d\ov{z}\rg]\\[5mm]
\ds\quad=\Re\lf[\frac{i}{2}\int_{\Sigma^g} d[\ov{a}\  b\  e^{2\la}]\wedge d\ov{z}\rg]-\Re\lf[\frac{i}{2}\int_{\Sigma^g} \ov{a}\  {\p}_z[b\  e^{2\la}]\ dz\wedge d\ov{z}\rg]\\[5mm]
\ds\quad=\int_{\Sigma^g} \lf<\al,\lf(\p\lf(\beta\res g_{\vec{\Phi}}\rg)\res_2\, g_{\vec{\Phi}}^{-1}\rg)\res g_{\vec{\Phi}}^{-1}\rg>_{g_{\vec{\Phi}}}\ d\mbox{vol}_{g_{\vec{\Phi}}}
\end{array}
\]
where
\[
g_{\vec{\Phi}}^{-1}=e^{-2\la}\ [\p_z\otimes\p_{\ov{z}}+\p_{\ov{z}}\otimes \p_z]\quad\mbox{ , }\quad (b \ dz\otimes d\ov{z}\otimes d\ov{z})\res_2 g_{\vec{\Phi}}^{-1}=e^{-2\la}\, b\ d\ov{z}\quad\mbox{and }\quad (e^{-2\la}\, b\ d\ov{z})\res g_{\vec{\Phi}}^{-1}=e^{-4\la}\, f\ \p_z
\]
Hence we have proved that the adjoint of $\ov{\p}$ on $\Gamma((T^\ast\Sigma^g)^{(0,1)}\otimes(T\Sigma^g)^{(1,0)})$ is given by
\[
\ov{\p}^\ast\ :\ \beta\in \Gamma((T^\ast\Sigma^g)^{(0,1)}\otimes(T\Sigma^g)^{(1,0)})\ \longrightarrow\ \ov{\p}^\ast\beta=\lf(\p\lf(\beta\res g_{\vec{\Phi}}\rg)\res_2\, g_{\vec{\Phi}}^{-1}\rg)\res g_{\vec{\Phi}}^{-1}\in \Gamma(\wedge^{(1,0)}\Sigma^g)
\]
We have $\mbox{Im}\,\ov{\p}=(\mbox{Ker}\,\ov{\p}^\ast)^\perp$. We have that 
\be
\label{kerdelbarstar}
\mbox{Ker}\,\ov{\p}^\ast=\lf\{\beta\in \Gamma((T^\ast\Sigma^g)^{(0,1)}\otimes(T\Sigma^g)^{(1,0)})\quad;\quad\beta\res g_{\vec{\Phi}}\in \mbox{Hol}_Q(\Sigma^g,g_{\vec{\Phi}}))\rg\}
\ee
Observe that
\[
\forall\,\gamma\in \Gamma\lf(\lf((T^\ast\Sigma^g)^{(0,1)}\rg)^{\otimes^2}\rg)\ ,\ \forall\,\beta\in \Gamma((T^\ast\Sigma^g)^{(0,1)}\otimes(T\Sigma^g)^{(1,0)})\quad
\lf<\gamma\res g^{-1}_{\vec{\Phi}},\beta\rg>_{g_{\vec{\Phi}}}=\lf<\gamma,\beta\res g_{\vec{\Phi}}\rg>_{g_{\vec{\Phi}}}
\]
This implies 
\be
\label{carac}
\gamma\res g^{-1}_{\vec{\Phi}}\in \mbox{Im}\,\ov{\p}\quad\Longleftrightarrow\quad \gamma\in \lf(\mbox{Hol}_Q(\Sigma^g,g_{\vec{\Phi}})\rg)^\perp\quad.
\ee
We deduce (\ref{delbar}) from (\ref{carac}) and (\ref{delbar-estim}) follows by classical estimates for elliptic complexes in Sobolev Spaces.\hfill $\Box$

\medskip

\noindent{\bf Continuation of the proof of theorem~\ref{th-hilbert}.} To $f\in(\mbox{Hol}_1(\Sigma^g))^\perp\cap \Gamma^{3,2}(\wedge^{(1,0)}\Sigma^g)$ solving (\ref{delbar}) we assign
\[
X:=2\,\Re(f)=2\,\Re((f_1+i\,f_2)\ \p_z)=(f_1\,\p_{x_1}+f_2\,\p_{x_2})=X_1\, \p_{x_1}+X_2\,\p_{x_2}\quad.
\]
Observe that, if we denote $\vec{X}:=d\vec{\Phi}\cdot X$, we have
\[
\ov{\p}\lf(\vec{X}\cdot \ov{\p}\vec{\Phi}\res g_{\vec{\Phi}}^{-1}\rg)=\ov{\p}\lf(e^{2\la} (X_1+i\,X_2) \ d\ov{z}\res g_{\vec{\Phi}}^{-1}\rg)=\ov{\p} f
\]
Observe also that, since $\vec{X}$ is tangent to the immersion $\vec{X}\cdot\ov{\p}\lf( \ov{\p}\vec{\Phi}\res g_{\vec{\Phi}}^{-1}\rg)=0$ hence
\[
\ov{\p}f=(\ov{\p}\vec{X}\cdot\ov{\p}\vec{\Phi})\res g_{\vec{\Phi}}^{-1}
\]
Using $\mbox{Im}\,\ov{\p}=(\mbox{Ker}\,\ov{\p}^\ast)^\perp$ and the characterization of $\mbox{Ker}\,\ov{\p}^\ast$ given by (\ref{kerdelbarstar}), we have $$\ov{\p}\vec{X}\cdot\ov{\p}\vec{\Phi}=P^\perp_{\vec{\Phi}}(\ov{\p}\vec{X}\cdot\ov{\p}\vec{\Phi})$$ and hence
\be
\label{rep}
\ov{\p}f=\ov{D}_{\vec{\Phi}}^\ast\vec{X} \res g_{\vec{\Phi}}^{-1}
\ee
We denote
\[
\lf\{
\begin{array}{l}
\ds{\mathcal X}^{3,2}(S^2):=\lf\{{X}\in \Gamma^{3,2}(TS^2)\ ;\ X(a_i)=0\quad i=1,2,3\rg\}\quad,\\[3mm]
\ds{\mathcal X}^{3,2}(T^2):=\lf\{{X}\in \Gamma^{3,2}(T T^2)\ ;\ X(a)=0\rg\}\quad,\\[3mm]
\ds{\mathcal X}^{3,2}(\Sigma^g)=\Gamma^{3,2}(\Sigma^g)\quad\mbox{ for }g>1
\end{array}
\rg.
\]
The space of Holomorphic Vector Field on $T^{(1,0)}S^2$ is a $3-$dimensional complex vector space given in ${\C}$, after the stereographic projection, by
\[
h(z)= (\al+\beta\, z+\gamma\, z^2)\ \p_z\quad\mbox{ where }(\al,\beta,\gamma)\in {\C}^3\quad.
\]
Whereas the space of Holomorphic Vector Field on $T^{(1,0)}T^2$ is a $1-$dimensional complex vector space given in ${\C}$ by
\[
h(z)=\al\ \p_z\quad\mbox{ where }\al\in {\C}\quad.
\]
while for $g>1$ we have $\mbox{Hol}_1(\Sigma^g)=\{0\}$. Hence, for any $g\in {\N}$ and any
\be
\label{hol-1}
f\in(\mbox{Hol}_1(\Sigma^g))^\perp\cap \Gamma^{3,2}(\wedge^{(1,0)}\Sigma^g)\quad\exists\ !\
h_f\in \mbox{Hol}_1(\Sigma^g)\quad\mbox{ s. t. }\quad\Re(f+h_f)\in {\mathcal X}^{3,2}(\Sigma^g)
\ee
moreover the map $f\rightarrow h_f$ from $(\mbox{Hol}_1(\Sigma^g))^\perp\cap \Gamma^{3,2}(\wedge^{(1,0)}\Sigma^g)$ into $\mbox{Hol}_1(\Sigma^g)$ is linear and smooth.
Hence we can summarize what we have proved so far in the following lemma.
\begin{Lm}
\label{lm-decomposition} Let $\vec{\Phi}$ be a $W^{3,2}-$immersion. Then the following holds
\[
\begin{array}{l}
\ds\forall\, \vec{v}\in \Gamma^{3,2}(\vec{\Phi}^\ast TM^m)\quad\exists \ !\ X\in{\mathcal X}^{3,2}(\Sigma^g) \quad\mbox{ s. t. }\quad\\[3mm]
\ds\ov{\p}\lf(X-\,i\,X^\perp\rg)=\ov{D}_{\vec{\Phi}}^\ast\vec{X}\res g_{\vec{\Phi}}^{-1}=\ov{D}_{\vec{\Phi}}^\ast\vec{v}\res g_{\vec{\Phi}}^{-1}
\end{array}
\]
where $\vec{X}=d\vec{\Phi}\cdot X$ and such that
\[
\|{X}\|_{W^{3,2}}\le C_{\vec{\Phi}}\ \|\vec{v}\|_{W^{3,2}}
\]
\hfill $\Box$
\end{Lm}
\noindent{\bf End of the proof of theorem~\ref{th-hilbert}.} Let $g_0$ be a smooth reference metric on $\Sigma^g$ and denote by $\exp^{g_0}$ the smooth
exponential map from $T\Sigma$ into $\Sigma$ associated to $g_0$. Let $\ep>0$ small and denote respectively
\[
{\mathcal X}^{3,2}_\ep(\Sigma^g):=\lf\{{X}\in {\mathcal X}^{3,2}(\Sigma^g)\ ;\ \|X\|_{W^{3,2}}<\ep\rg\}\quad,
\]
and
\[
{\mathcal D}^\perp_\ep:=\lf\{\Psi\in \mbox{Diff}_+^\ast(\Sigma)\ ;\quad \exists\, X\in {\mathcal X}^{3,2}_\ep(\Sigma^g)\quad \mbox{s.t}\quad \Psi(x)= \exp^{g_0}_x(X(x)) \rg\}\quad.
\]
We define
\[
\begin{array}{rcl}
\Lambda_{\vec{\Phi}}\ :\ {\mathcal V}_{\vec{\Phi}}^\ep\times {\mathcal D}^\perp_\ep &\longrightarrow& \Gamma^{2,2}((T^\ast\Sigma )^{(0,1)}\otimes (T^\ast\Sigma )^{(0,1)})\\[3mm]
(\vec{\Xi},\Psi) &\longrightarrow& \ov{D}_{\vec{\Phi}}^\ast\lf(\vec{\Xi}\circ\Psi\rg)\res g^{-1}_{\vec{\Phi}}
\end{array}
\]
The map is clearly $C^1$ and lemma~\ref{lm-decomposition}  gives that
\[
\lf.\p_{\Psi}\Lambda_{\vec{\Phi}}\rg|_{(\vec{\Phi},0)}\cdot X=\ov{D}_{\vec{\Phi}}^\ast\lf(d\vec{\Phi}\cdot X\rg)\res g^{-1}_{\vec{\Phi}}
\]
realizes an isomorphism between ${\mathcal X}^{3,2}$ and ${\mathcal I}$ (defined in (\ref{mathcali})). The implicit function theorem gives then the existence
of a $C^1$ map $\Psi_{\vec{\Phi}}(\vec{\Xi})$ defined in a neighborhood of $\vec{\Phi}$ such that
\[
\ov{D}_{\vec{\Phi}}^\ast\lf(\vec{\Xi}\circ\Psi_{\vec{\Phi}}(\vec{\Xi})\rg)\res g^{-1}_{\vec{\Phi}}=0
\]
and we denote $\vec{w}_{\vec{\Phi}}(\vec{\Xi}):=\vec{\Xi}\circ\Psi_{\vec{\Phi}}(\vec{\Xi})-\vec{\Phi}$. Because of the local uniqueness of $\Psi_{\vec{\Phi}}(\vec{\Xi})$ given by
the implicit function theorem, for any element $\Psi_0\in \mbox{Diff}_+^\ast(\Sigma)$ close to the identity and $\vec{\Xi}$ close enough to $\vec{\Phi}$ one has trivially
\[
\ov{D}_{\vec{\Phi}}^\ast\lf(\vec{\Xi}\circ\Psi_0\circ\Psi_0^{-1}\circ\Psi_{\vec{\Phi}}(\vec{\Xi})\rg)\res g^{-1}_{\vec{\Phi}}=0
\]
hence we deduce the equivariance property
\[
\Psi_{\vec{\Phi}}(\vec{\Xi}\circ{\Psi}_0)=\Psi_0^{-1}\circ\Psi_{\vec{\Phi}}(\vec{\Xi})\quad\mbox{ and }\quad\vec{w}_{\vec{\Phi}}(\vec{\Xi}\circ\Psi_0):=\vec{\Xi}\circ\Psi_{\vec{\Phi}}(\vec{\Xi})-\vec{\Phi}=\vec{w}_{\vec{\Phi}}(\vec{\Xi})
\]
This permits to extend $\Upsilon_{\vec{\Phi}}:=(\vec{w}_{\vec{\Phi}},\Psi_{\vec{\Phi}})$ on a neighborhood ${\mathcal O}_{\vec{\Phi}}$ of $\vec{\Phi}$ invariant under the action of $\mbox{Diff}_+^\ast(\Sigma)$.

\medskip

We are now proving the Hausdorff property for ${\mathfrak M}_{g}(\Sigma^g,M^m):=\mbox{Imm}_0(\Sigma^g,M^m)/\mbox{Diff}^{\, \ast}_+(\Sigma^g)$. Following classical considerations
(see the arguments in the proof of lemma 2.9.9 of \cite{Var}) it suffices to prove that 
\[
\Gamma:=\lf\{ (\vec{\Phi},\vec{\Phi}\circ\Psi)\quad;\quad\vec{\Phi}\in \mbox{Imm}_0(\Sigma^g,M^m)\ \mbox{ and }\ \Psi\in \mbox{Diff}^{\, \ast}_+(\Sigma^g)\rg\}
\]
is closed in $(\mbox{Imm}_0(\Sigma^g,M^m))^2$. This follows fro the first part of the proof of the theorem. Let $(\vec{\Phi}_k,\vec{\Xi}_k:=\vec{\Phi}_k\circ\Psi_k)\rightarrow (\vec{\Phi}_\infty,\vec{\Xi}_\infty)$ in $W^{3,2}$. For $k$ large enough both $\vec{\Phi}_k$ and $\vec{\Xi}_k$ are included in ${\mathcal O}_{\vec{\Phi}_\infty}$. Because of the continuity of the map
 $\vec{w}_{\vec{\Phi}_\infty}$ we have respectively
 \[
 \vec{w}_{\vec{\Phi}_\infty}(\vec{\Phi}_k)\rightarrow  \vec{w}_{\vec{\Phi}_\infty}(\vec{\Phi}_\infty)=0\quad\mbox{ and }\quad\vec{w}_{\vec{\Phi}_\infty}(\vec{\Xi}_k)\rightarrow  \vec{w}_{\vec{\Phi}_\infty}(\vec{\Xi}_\infty)
 \]
The equivariance of $\vec{w}_{\vec{\Phi}_\infty}$ gives $\vec{w}_{\vec{\Phi}_\infty}(\vec{\Xi}_k)=\vec{w}_{\vec{\Phi}_\infty}(\vec{\Phi}_k)$ hence  $\vec{w}_{\vec{\Phi}_\infty}(\vec{\Xi}_\infty)=0$.
Thus $\vec{\Xi}_\infty\circ\Psi_{\vec{\Phi}_\infty}=\vec{\Phi}_\infty$ and this shows that $\Gamma$ is closed and then ${\mathfrak M}(\Sigma^g,M^m)$ defines an Hausdorff {\it Hilbert Manifold}
and theorem~\ref{th-hilbert} is proved. \hfill $\Box$

\subsection{The Relaxed Area $A^\sigma$ and the Fredholm Property for $D^2A^\sigma(\vec{\Phi})$ on $T_{\vec{\Phi}}{\mathfrak M}$.}
We shall denote
\[
{\mathfrak M}(M^m):=\bigcup_{g\in {\N}}\mbox{Imm}_0(\Sigma^g,M^m)/\mbox{Diff}^{\, \ast}_+(\Sigma^g)\quad.
\]
and
\[
{\mathfrak M}^{g_0}(M^m):=\bigcup_{g\le g_0}\mbox{Imm}_0(\Sigma^g,M^m)/\mbox{Diff}^{\, \ast}_+(\Sigma^g)\quad.
\]
\medskip

For any immersion $\vec{\Phi}\in \mbox{Imm}(M^m)$ we denote
\[
F(\vec{\Phi}):=\int_{\Sigma_{\vec{\Phi}}}\lf[1+|{\vec{\mathbb I}}_{\vec{\Phi}}|^2\rg]^2\ dvol_{g_{\vec{\Phi}}}
\]
where ${\vec{\mathbb I}}_{\vec{\Phi}}$ is the second fundamental form of the immersion $\vec{\Phi}$ in $M^m$. 

\medskip

Observe that 
\be
\label{I.001}
\forall\ g_0\in{\N} \quad\quad\exists\ C_{g_0}>0\quad\quad F(\vec{\Phi})<C_{g_0}\quad\Longrightarrow\quad\mbox{genus}\,(\Sigma_{\vec{\Phi}})\le g_0
\ee
This is a direct consequence of {\it Gauss Bonnet theorem} and {\it Cauchy Schwartz inequality}.

It is clear that $F(\vec{\Phi})$ only depends on the equivalence class $[\vec{\Phi}]$
of $\vec{\Phi}$ in ${\mathfrak M}(\Sigma^g,M^m)$. Since $F$ is a smooth functional on $\mbox{Imm}_0(\Sigma^g,M^m)$ (see \cite{Riv-minmax}) it descends to a smooth functional on 
${\mathfrak M}(\Sigma^g,M^m)$. We shall now prove the following theorem.
\begin{Th}
\label{th-fredholm}
Let $[\vec{\Phi}]$ be a critical point of $F$ in ${\mathfrak M}(M^m)$. Then the second derivative of $F$ at $[\vec{\Phi}]$ defines a Fredholm operator .\hfill $\Box$
\end{Th}
\noindent{\bf Proof of theorem~\ref{th-fredholm}.} From \cite{KLL} we know that $\vec{\Phi}$ is smooth in a conformal coordinates. We shall be working in the chart in the neighborhood
of $[\vec{\Phi}]$ in ${\mathfrak M}(\Sigma^g,M^m)$ given
by $\vec{w}_{\vec{\Phi}}$ from theorem~\ref{th-hilbert}. In other words  we identify 

\be
\label{fred-1}
T_{[\vec{\Phi}]}{\mathfrak M}\simeq\lf\{\vec{w}\in \Gamma^{3,2}(\vec{\Phi}^\ast T M^m)\quad;\quad P^\perp_{\vec{\Phi}}\lf(\ov{\p}\vec{w}\,\dot{\otimes}\,\ov{\p}\vec{\Phi}\rg)=0\rg\}
\ee
For such a $\vec{w}$ we denote by $q_{\vec{w}}$ the holomorphic quadratic form given by
\[
\ov{\p}\vec{w}\,\dot{\otimes}\,\ov{\p}\vec{\Phi}=q_{\vec{w}}\quad.
\]
After contracting with the tensor $g_{\vec{\Phi}}^{-1}$, this equation becomes
\[
\ov{\p}\lf(\vec{w}\cdot\ov{\p}\vec{\Phi}\res g_{\vec{\Phi}}^{-1}\rg)=-\pi_{\vec{n}}(\vec{w})\cdot \vec{h}^0+q_{\vec{w}}\res g_{\vec{\Phi}}^{-1}
\]
where $\vec{h}^0$ is the trace free part of the second fundamental form, which is orthogonal to the tangent plane of the immersion and given in local coordinates by 
\[
\vec{h}^0_{\vec{\Phi}}=\p_{\ov{z}}\lf(e^{-2\la}\p_{\ov{z}}\vec{\Phi}\rg)\ d\ov{z}\otimes \p_z\quad.
\]
Using the characterization of $\mbox{Im}\,\ov{\p}=(\mbox{Ker}\,\ov{\p}^\ast)^\perp$ given by  (\ref{carac}) we deduce
\[
\ov{\p}\lf(\vec{w}\cdot\ov{\p}\vec{\Phi}\res g_{\vec{\Phi}}^{-1}\rg)=-{\mathfrak P}_{\vec{\Phi}}\lf(\pi_{\vec{n}}(\vec{w})\cdot \vec{h}^0_{\vec{\Phi}}\rg)
\]
where ${\mathfrak P}_{\vec{\Phi}}$ is the orthogonal projection onto $(\mbox{Hol}_Q(\Sigma^g,g_{\vec{\Phi}})\res  g_{\vec{\Phi}}^{-1})^\perp$. Denote $\vec{X}_{\vec{w}}$ the projection
of $\vec{w}$ onto the tangent plane (i.e. $\vec{X}_{\vec{w}}=\vec{w}-\pi_{\vec{n}}(\vec{w})$) and let $X_{\vec{w}}$ be the vector field on $\Sigma$ such that $d\vec{\Phi}\cdot X_{\vec{w}}=\vec{X}_{\vec{w}}$. Following the computations from the previous subsection we deduce
\be
\label{fred-2}
\ov{\p}\lf(X_{\vec{w}}-\,i\,X_{\vec{w}}^\perp\rg)=-{\mathfrak P}_{\vec{\Phi}}\lf(\pi_{\vec{n}}(\vec{w})\cdot \vec{h}^0_{\vec{\Phi}}\rg)
\ee
Denote
\[
\begin{array}{rcl}
\ds\pi_T\ :\ \Gamma^{3,2}(\vec{\Phi}^\ast TM^m) &\longrightarrow & \Gamma^{3,2}((T\Sigma)^{(1,0)})\\[3mm]
\ds\vec{w} &\longrightarrow & X_{\vec{w}}-\,i\,X_{\vec{w}}^\perp
\end{array}
\]
In view of the expression of the second derivative $D^2F$ given by (\ref{deuxi-18}) we have that, modulo compact operators (remembering that $\vec{\Phi}$ is smooth), We are reduced\footnote{The sum of a {\it Fredholm operator} with a {\it compact operator} is {\it Fredholm}.} to study the Fredholm nature of the operator generated by the following quadratic form
\[
\begin{array}{l}
\ds Q_{\vec{\Phi}}(\vec{w})= \,\, \int_{\Sigma}(1+|\vec{\mathbb I}_{\vec{\Phi}}|^2_{g_{\vec{\Phi}}}) \lf|\pi_{\vec{n}}\lf(D^{g_{\vec{\Phi}}}d\vec{w}\rg)\rg|_{g_{\vec{\Phi}}}^2\ dvol_{g_{\vec{\Phi}}}
+2\, \int_{\Sigma}\, \lf|\lf<\vec{\mathbb I}_{\vec{\Phi}},D^{g_{\vec{\Phi}}}d\vec{w}\rg>_{g_{\vec{\Phi}}}\rg|^2\ dvol_{g_{\vec{\Phi}}}
\end{array}
\]
combined with (\ref{fred-2}). Hence the Symbols of the generated operator, in local conformal coordinates, is given by
\[
\lf\{
\begin{array}{l}
2\,e^{-2\la} \pi_{\vec{n}}\circ\lf[(1+|\vec{\mathbb I}_{\vec{\Phi}}|^2_{g_{\vec{\Phi}}}) |\xi|^4+2\,e^{-4\la}\sum_{i,j,k,l} \vec{\mathbb I}_{kl}\otimes\vec{\mathbb I}_{kl} \ \xi_i\,\xi_j\,\xi_k\,\xi_l\rg]\circ\pi_{\vec{n}}\\[5mm]
(\xi_1+i\,\xi_2)\circ\pi_T
\end{array}
\rg.
\]
This is clearly the symbol defining an {\it elliptic operator} on $\Gamma^{3,2}(\vec{\Phi}^\ast TM^m)$ and $D^2F$ is {\it Fredholm} on $T_{[\vec{\Phi}]}{\mathfrak M}$. This concludes the proof
of theorem~\ref{th-fredholm}.\hfill $\Box$

\subsection{Lower Semi-Continuity of the Morse Index in the Viscosity Method.}

In the applications below, we shall mostly consider
the {\it area lagrangian}
\[
\mbox{Area}(\vec{\Phi})=\int_{\Sigma_{\vec{\Phi}}}\, dvol_{g_{\vec{\Phi}}}
\]
or its relaxations of the form
\[
A^\sigma(\vec{\Phi}):=\mbox{Area}(\vec{\Phi})+\sigma^2\, F(\vec{\Phi})=\mbox{Area}(\vec{\Phi})+\sigma^2\, \int_{\Sigma_{\vec{\Phi}}}\lf[1+|{\vec{\mathbb I}}_{\vec{\Phi}}|^2\rg]^2\ dvol_{g_{\vec{\Phi}}}\quad.
\]
where $\sigma>0$. The work \cite{Riv-minmax}  has been devoted to the asymptotic analysis of sequences of critical points of $A^{\sigma_k}$, with uniformly bounded 
$A^{\sigma_k}$ energy and satisfying the so called {\it entropy condition}
\be
\label{entropy}
\sigma_k^2\, F(\vec{\Phi}_{k})= o\lf(\frac{1}{\log\sigma_k^{-1}}\rg)\quad\mbox{ as $\sigma_k$ goes to zero.}
\ee
It is proved in this two work that, modulo extraction of a subsequence, the immersions $\vec{\Phi}_{k}$ {\it varifold converges} towards a 2-dimensional {\it integer rectifiable stationary varifold} $\mathbf v_\infty$ of $M^m$ which is parametrized\footnote{The regularity of such \underbar{parametrized}  integer rectifiable stationary varifold, which is a very peculiar subclass of integer rectifiable stationary varifold, is studied in \cite{PiR}. } precisely there exists a riemann surface $(S_\infty,h)$ a map $\vec{\Psi}_\infty\in W^{1,2}(S_\infty,M^m)$ which is almost everywhere conformal and an integer multiplicity $\theta\in L^\infty(S_\infty,{\N}^\ast)$ such that
\be
\label{integer-target}
\forall \ F\in C^\infty(M^m)\quad\quad\int_{S_\infty}\theta\ \lf[\lf<d(F\circ\vec{\Psi}_\infty)\cdot d \vec{\Psi}_\infty\rg>_h-\ F(\vec{\Psi}_\infty)\ \vec{\mathbb I}_{M^m}(d\vec{\Psi}_\infty,d\vec{\Psi}_\infty)_h\rg]\ dvol_h=0\quad.
\ee
where $\vec{\mathbb I}_{M^m}$ is the second fundamental form of $M^m\hookrightarrow {\R}^Q$. When 
\be
\label{multip-bound}
\limsup_{k\rightarrow +\infty}\mbox{Area}(\vec{\Phi}_k)<\frac{8\pi}{\|\vec{\mathbb I}_{M^m}\|_\infty}\quad,
\ee
then $\theta\equiv 1$ and according to \cite{Riv-reg} the map $\vec{\Psi}_\infty$ defines a smooth  minimal embedding. Indeed  we have by (\ref{integer-target}) that
${\mathbf v}_\infty$ has an $L^\infty$ generalized mean curvature in ${\R}^Q$ equal $|d\vec{\Psi}_\infty|^2_h\ dvol_h$ a.e. on $S_\infty$ to
\[
\vec{H}_{{\R}^Q}:= -\,|d\vec{\Psi}_\infty|^{-2}_h\ \vec{\mathbb I}(d\vec{\Psi}_\infty,d\vec{\Psi}_\infty)_h\ 
\]
 Under the assumptions (\ref{multip-bound}), using Li-Yau upper-bound of the density in ${\R}^Q$ (see \cite{LiY}), we have
\be
\label{upp}
\Theta_\infty(\vec{p})\le \frac{1}{4\pi}\int_{S_\infty}|\vec{H}_{{\R}^Q}|^2\ N\ 2^{-1}\ |d\vec{\Psi}_\infty|^2_h\ dvol_h\le\ \frac{\|\vec{\mathbb I}_{M^m}\|_\infty}{4\pi}\  \limsup_{k\rightarrow +\infty}\mbox{Area}(\vec{\Phi}_k)<2
\ee
where $\Theta_\infty$ is the density of the varifold ${\mathbf v}_\infty$ at $\vec{p}$. This gives that $\theta<2$ everywhere and hence $\theta=1$  a.e. on $S_\infty$ since ${\mathbf v}_\infty$ is integer. 

\medskip

 The question to compare the {\it Morse index} of the limiting surface for the area with the {\it Morse index} of the sequence $\vec{\Phi}_k$ for the relaxed functionals $A^{\sigma_k}$ was left open in these works. We are now giving an answer to that question assuming the upper bound (\ref{multip-bound}).

\medskip
\begin{Th}
\label{th-A.3}
Let $\vec{\Phi}_k$ be a sequence of immersions of a closed surface $\Sigma^g$, critical points of $A^{\sigma_k}$ and such that
\[
\limsup_{k\rightarrow+\infty}A^{\sigma_k}(\vec{\Phi}_k)<+\infty\quad\mbox{and }\quad\sigma_k^2\int_{\Sigma^g}(1+|\vec{\mathbb I}_{\vec{\Phi}_k}|^2)^2\ dvol_{g_{\vec{\Phi}_k}}=o\lf( \frac{1}{\log\sigma_k^{-1}} \rg)
\]
Then there exists a subsequence $\vec{\Phi}_{k'}$ such that the corresponding immersed surface converges in varifolds towards a parametrized integer rectifiable stationary varifold ${\mathbf v}_\infty:=(S_\infty,\vec{\Psi}_\infty,N)$ moreover we have
\[
\mbox{genus}\,(S_\infty)\le g\quad.
\]
Finally if (\ref{multip-bound}) holds or more generally if we know that $\vec{\Phi}_k$ bubble tree $W^{1,2}$ converges\footnote{We expect this to hold for any quasi-minimizing sequence and in particular for solutions to Minmax problems in general.} then  $\theta\equiv 1$ and $\vec{\Psi}_\infty$ is a minimal conformal  immersion satisfying
\be
\label{A.I.14}
\mbox{Ind } (\vec{\Psi}_\infty)\le \liminf_{k\rightarrow \infty} \mbox{Ind}^{\,\sigma_{k'}}(\vec{\Phi}_{k'})
\ee
where $\mbox{Ind} (\vec{\Psi}_\infty)$  is the maximal dimension of a subspace of $T_{\vec{\Psi}_\infty(S_\infty)}{\mathfrak M}$ on which $D^2\mbox{Area}(\vec{\Psi}_\infty)$ is strictly negative and $\mbox{Ind}^{\,\sigma_k}(\vec{\Phi}_k)$ is the maximal dimension of a subspace of $T_{\vec{\Phi}_k(\Sigma^g)}{\mathfrak M}$ on which $D^2 A^{\sigma_k}(\vec{\Phi}_k)$ is strictly negative.
\hfill $\Box$
\end{Th}
\noindent{\bf Proof of theorem~\ref{th-A.3}.} We shall assume that $\Sigma^g$ is connected. We shall present the computations  for $M^m=S^m$. The general constraint generates lower order terms whose abundance could mask the true reason why the theorem is true whereas the same terms in the $M^m=S^m$ case are easier to present. The first part of the theorem is the  main results of \cite{Riv-minmax} . It remains to prove the inequality
(\ref{A.I.14}) under the assumption (\ref{multip-bound}). The first derivative of the area of an immersion (possibly branched) of a closed surface $\Sigma$ into ${\R}^Q$ is given by (see \cite{Riv-minmax})
\[
D\mbox{Area}(\vec{\Phi})\cdot\vec{w}=\int_\Sigma \lf<d\vec{\Phi}\,;\,d\vec{w}\rg>_{g_{\vec{\Phi}}}\ d\mbox{vol}_{g_{\vec{\Phi}}}
\]
and the second derivative\footnote{ A reader familiar to the rich literature in geometry and geometric measure theory on minimal surface theory in 3 dimension might not recognize the most commonly used expression of the second derivative of the area by the mean of the {\it Jacobi field}. This classical presentation of $D^2\mbox{Area}$ has the advantage to ``reduce'' this operator to an operator on function by introducing the decomposition $\vec{w}=w\, \vec{n}$. This decomposition however is not ``analytically'' favorable since $\vec{n}$ has a-priori one degree of regularity less than $\vec{w}$. This observation is at the base of the analysis of the {\it Willmore functional} as it has been developed by the author in the recent years.   }
\[
D^2 \mbox{Area}(\vec{\Phi})\cdot(\vec{w},\vec{w})=\int_\Sigma \lf[\lf<d\vec{w}\,;\,d\vec{w}\rg>_{g_{\vec{\Phi}}}+\lf|\lf<d\vec{\Phi}\,;\,d\vec{w}\rg>_{g_{\vec{\Phi}}}\rg|^2 - 2^{-1}\lf|d\vec{\Phi}\dot{\otimes}\, d\vec{w}+d\vec{w}\dot{\otimes}\, d\vec{\Phi}\rg|^2\rg]\ d\mbox{vol}_{g_{\vec{\Phi}}}
\]
where in coordinates $d\vec{\Phi}\dot{\otimes} d\vec{w}:=\sum_{i,j}\p_{x_i}\vec{\Phi}\cdot\p_{x_j}\vec{w}\ dx_i\otimes dx_j$. 

\medskip

Since $\theta=1$  a.e. on $S_\infty$ and, following the proof of the main theorem of \cite{Riv-minmax}, 
we can extract a subsequence that we keep denting $\vec{\Phi}_{k}$ such that we have a {\it bubble tree } strong $W^{1,2}$ convergence
of $\vec{\Phi}_{k}$ towards a minimal (possibly branched) immersion $\vec{\Psi}_\infty$ of a surface $S_\infty$. More precisely, if one denotes $\{S^j_\infty\}_{j\in J}$ to be the connected components of $S_\infty$, for every $j\in J$ there exists $N^j$ points $\{a^{j,l}\}_{l=1\cdots N^j}$ (containing in particular the possible branched points of $\vec{\Psi}_\infty$ and a converging sequence of constant scalar curvature metrics  $h^j_{k}$ of volume one  and for any $\delta>0$ a sequence of conformal embeddings $\phi^j_{k}$ from $(S^j_\infty\setminus \cup_{l=1}^{N^j} B_\delta(a^{j,l}), h^j_{k})$ into $(\Sigma^g,g_{\vec{\Phi}_{k}})$ such that
\be
\label{A.I.15}
\vec{\Psi}^j_{k}:=\vec{\Phi}_{k}\circ \phi^j_{k}\longrightarrow {\vec{\Psi}_\infty}\quad\quad\mbox{ strongly in } W^{1,2}_{loc}(S^j_\infty\setminus \cup_{l=1}^{N^j} B_\delta(a^{j,l}))
\ee
We shall need the following intermediate lemma
\begin{Lm}
\label{czero}
Under the previous assumptions and notations we have
\be
\label{A.I.15-b}
\vec{\Psi}^j_{k}:=\vec{\Phi}_{k}\circ \phi^j_{k}\longrightarrow {\vec{\Psi}_\infty}\quad\quad\mbox{ strongly in } C^0_{loc}(S^j_\infty\setminus \cup_{l=1}^{N^j} B_\delta(a^{j,l}))
\ee
\hfill $\Box$
\end{Lm}
\noindent{\bf Proof of lemma~\ref{czero}.} Let $x_0\in S^j_\infty\setminus \cup_{l=1}^{N^j} B_\delta(a^{j,l})$. For $r$ small enough we have in a conformal chart
with respect to the limiting metric $h^j_\infty$
\be
\label{A.I.15-d}
\lf\|
\lf(\begin{array}{c}\p_{x_1}\vec{\Psi}_\infty\\[3mm]
\p_{x_2}\vec{\Psi}_\infty
\end{array}\rg)
- e^{\la(x_0)}\ \lf(
\begin{array}{c}
\vec{e}_1(x_0)\\[3mm]
\vec{e}_2(x_0)
\end{array}
\rg)
\rg\|_{L^\infty(B_r(x_0))}\le C\, r
\ee
where $(\vec{e}_1(x_0),\vec{e}_2(x_0))$ is an  orthonormal family in $T_{\vec{\Psi}_\infty(x_0)}M^m$. Let $\ep>0$. Because of the strong $W^{1,2}$ convergence, for $k$ large enough we have (omitting to write the upper index $j$)
\be
\label{A.I.15-e}
\int_{B_{2\ep^2}(x_0)}\lf|
\lf(\begin{array}{c}\p_{x_1}\vec{\Psi}_k\\[3mm]
\p_{x_2}\vec{\Psi}_k
\end{array}\rg)
- e^{\la(x_0)}\ \lf(
\begin{array}{c}
\vec{e}_1(x_0)\\[3mm]
\vec{e}_2(x_0)
\end{array}
\rg)
\rg|^2\ dx^2\le  C\ \ep^4
\ee
and
\[
\int_{B_{2\ep^2}(x_0)}\lf|\vec{\Psi}_k-\vec{\Psi}_\infty\rg|^2\ dx^2\le \ep^8
\]
Hence, using the {\it mean value theorem}, we have the existence of $r_k\in(\ep^2,2\ep^2)$ such that
\[
\int_{\p B_{r_k}(x_0)}\lf|
\lf(\begin{array}{c}\p_{x_1}\vec{\Psi}_k\\[3mm]
\p_{x_2}\vec{\Psi}_k
\end{array}\rg)
- e^{\la(x_0)}\ \lf(
\begin{array}{c}
\vec{e}_1(x_0)\\[3mm]
\vec{e}_2(x_0)
\end{array}
\rg)
\rg|\ dl\le  C\ \ep^2
\]
and
\be
\label{A.I.15-c}
\int_{\p B_{r_k}(x_0)}\lf| \vec{\Psi}_k-\vec{\Psi}_\infty\rg|\ dl\le C\ \ep^4
\ee
This implies in particular
\[
\int_0^{2\pi}\lf|\p_\theta\lf[\vec{\Psi}_k(r_k,\theta)- r_k\ e^{\la(x_0)}\ \cos\theta\ \vec{e}_1(x_0)-r_k\ e^{\la(x_0)}\ \sin\theta\ \vec{e}_2(x_0)\rg]\rg|\ d\theta\le C\ \ep^2
\]
which gives
\[
\|\vec{\Psi}_k(r_k,\theta)-\vec{\Psi}_k(r_k,0)- r_k\ e^{\la(x_0)}\ (\cos\theta-1)\ \vec{e}_1(x_0)-r_k\ e^{\la(x_0)}\ \sin\theta\ \vec{e}_2(x_0)\|_{L^\infty([0,2\pi])}\le C\ep^2
\]
From (\ref{A.I.15-c}) we deduce that there exists $\theta_k\in [0,2\pi]$ such that $|\vec{\Psi}_k(r_k,\theta_k)-\vec{\Psi}_\infty(r_k,\theta_k)|<\, C\, \ep^2$ hence we deduce, since
$\|\vec{\Psi}_\infty-\vec{\Psi}_\infty(x_0)\|_{L^\infty(B_{r_k}(x_0))}\le C\, \ep^2$ that
\[
\|\vec{\Psi}_k(r_k,\theta)-\vec{\Psi}_\infty(x_0)\|_{L^\infty([0,2\pi])}\le C\ep^2
\]
Assume there would exists a point $x_1\in B_{r_k}(x_0)$ such that
\[
|\vec{\Psi}_k(x_1)-\vec{\Psi}_\infty(x_0)|>\ep
\]
Then, $\epsilon$ being fixed, using the almost monotonicity formula lemma III.1 of \cite{Riv-minmax}, for $k$ large enough we would obtain 
\[
\mbox{Area}(\vec{\Psi}_k(B_{r_k}(x_0)))\ge c_0\ \ep^2
\]
Since $h_k^j$ converges strongly in any $C^l$ norm towards $h^j_\infty$ and since $\vec{\Psi}_k$ is conformal with respect to $h^j_k$, we would then deduce
\[
\int_{B_{r_k}(x_0)}|\nabla\vec{\Psi}_k|^2\ dx^2\ge c_1\ \ep^2
\]
But from (\ref{A.I.15-d}) and (\ref{A.I.15-e}) we have $\int_{B_{r_k}(x_0)}|\nabla\vec{\Psi}_k|^2\ dx^2\le C\ \ep^4$. This is a contradiction. Thus we have proved that
for $k$ large enough we have
\[
\|\vec{\Psi}_k-\vec{\Psi}_\infty\|_{L^\infty(B_{\ep^2}(x_0))}\le \ep
\]
which implies the lemma.\hfill $\Box$

\medskip

\noindent{\bf End of the proof of theorem~\ref{th-A.3}.}

For $\delta$ small enough and $k'$ large enough the  subdomains $\Omega_{k}^j(\delta):=\phi^j_{k}(S^j_\infty\setminus \cup_{l=1}^{N^j} B_\delta(a^{j,l}))$ are disjoint and 
\[
\lim_{\delta\rightarrow 0}\lim_{k\rightarrow +\infty}\mbox{Area}\lf(\vec{\Phi}_{k}\lf(\Sigma^g\setminus \bigcup_{j\in J}\Omega_{k}^j(\delta)\rg)\rg)=0
\]
Let $\vec{w}_1\cdots\vec{w}_N$ a family of $N$ independent smooth vectors in $W^{2,4}(\vec{\Psi}_\infty^\ast T{M}^m)$ on the Span of which
$D^2\mbox{Area}$ is strictly negative. We can assume without loss of generality that the $\vec{w}_i$ are $C^\infty$One modifies each of these vectors in the following way. For each $i\in\{1\cdots Q\}$ for each $j\in J$ and each $l\in\{1\cdots N^j\}$ we introduce (after identifying for each  $j$ and $l$ the tangent planes to $M^m$ around $\vec{\Phi}_\infty(a^{j,l})$ with the one at exactly $\vec{\Phi}_\infty(a^{j,l})$)
\[
\vec{w}^{\,\delta}_i(x)=\lf\{\begin{array}{l}
\vec{w}_i(x)\quad\quad\quad\quad \mbox{for }|a^{j,l}-x|\ge \sqrt{\delta}\\[3mm]
\vec{w}_i(x)\ \chi^{\,\delta}(|x-a^{j,l}|)\quad \quad\mbox{ for }\delta\le |a^{j,l}-x|\le \sqrt{\delta}\\[3mm]
0\quad\quad\quad\quad \mbox{ for } |a^{j,l}-x|\le \delta
\end{array}
\rg.
\] 
 where we take $\chi^{\,\delta}(s)$ to be  a slight smoothing of
$
\log({s}/{\delta})/\log({1}/{\sqrt{\delta}})
$
One verifies that $\vec{w}_i^{\,\delta}\in C^\infty(S_\infty)$ strongly converges towards $\vec{w}_i$ in $W^{1,2}(S_\infty,{\R}^Q)$ and therefore, in view of the explicit expression of $D^2 \mbox{Area}(\vec{\Psi}_\infty)\cdot(\vec{w},\vec{w})$, there exists $\delta$ small enough such that $\vec{w}^{\,\delta}_1\cdots\vec{w}^{\,\delta}_N$ a family of $N$ independent smooth vectors in $W^{2,4}(\vec{\Psi}_\infty^\ast TM^m)$ on the Span of which
$D^2\mbox{Area}$ is strictly negative. We fix such a $\delta$.

Let $\rho >0$ small enough such that for any $z\in M^m$ the map $\vec{\Psi}_\infty$ is injective
on each components of $\vec{\Psi}^{-1}_\infty(\ov{B^Q_\rho(z)})\subset S^j_\infty$. Let $\{\chi_s(z)\}_{s\in \{1\cdots N\}}$ be a finite smooth partition of unity of $M^m\subset {\R}^Q$ such that the support of every $\chi_s$ is included in an $m-$ball of radius $\rho$. We denote the connected components of $\vec{\Psi}_\infty^{-1}(\mbox{Supp}(\chi_s))$ in $S_\infty$ by $\Om_{s}^t$ for $t=1\cdots n_{s}$ and $\om_s^t$ are the corresponding characteristic functions. We have that $\chi_s(\vec{\Psi}_\infty(x))\ \om_s^t(x)$ is smooth for any $s\in \{1\cdots N\}$ and any
 $t\in\{1\cdots n_s\}$ and moreover
 \[
 d(\chi_s(\vec{\Psi}_\infty(x))\ \om_s^t(x))=d(\chi_s(\vec{\Psi}_\infty(x)))\ \om_s^t(x)
\]
 We can then write each $\vec{w}_i^{\,\delta}$ in the form
\[
\vec{w}_i^{\,\delta}(x)=\sum_{s=1}^N\chi_s(\vec{\Psi}_\infty(x))\sum_{t=1}^{n_{s}}\vec{v}_{t,s}(\vec{\Psi}_\infty(x))\ \om_s^t
\]
where $\vec{v}_{t,s}$ are smooth functions. For any $s=\in \{1\cdots N\}$ since the components $\ov{\Om_s^t}$ are disjoint to each other for $t\in\{1\cdots n_s\}$ We can include them in strictly larger disjoint open sets $\ov{\Om_s^t}\subset\ti{\Om_s^t}$ and we denote $\ti{\om}_s^t$ the corresponding characteristic functions. We still have of course
\[
\vec{w}_i^{\,\delta}(x)=\sum_{s=1}^N\chi_s(\vec{\Psi}_\infty(x))\sum_{t=1}^{n_{s}}\vec{v}_{t,s}(\vec{\Psi}_\infty(x))\ \ti{\om}_s^t
\]
moreover, since $\vec{\Psi}_k$ uniformly converges of  towards $\vec{\Psi}_\infty$ on $S^j_\infty\setminus \cup_{l=1}^{N^j} B_\delta(a^{j,l})$, for $k$ large enough, we have for every
$s$ and $t$
\[
 d(\chi_s(\vec{\Psi}_k(x))\ \tilde{\om}_s^t(x))=d(\chi_s(\vec{\Psi}_k(x)))\ \tilde{\om}_s^t(x)
\]
It is clear that 
\be
\label{A.I.16-a}
\vec{w}_{i,k}^{\,\delta}(x):=\sum_{s=1}^N\chi_s(\vec{\Psi}_k(x))\sum_{t=1}^{n_{s}}\vec{v}_{t,s}(\vec{\Psi}_k(x))\ \ti{\om}_s^t\ \longrightarrow\ \vec{w}_i^{\,\delta}(x)\mbox{ strongly in }W^{1,2}_{loc}(S^j_\infty\setminus \cup_{l=1}^{N^j} B_\delta(a^{j,l}))
\ee
Using the compositions with the  maps $(\phi^{j,k'})^{-1}$ we extend the $\vec{w}^{\,\delta}_{i,k}$, that we still denote $\vec{w}_{i,k}^{\,\delta}$ to the whole of $\Sigma^g$ by taking
$\vec{w}_{i,k}^{\,\delta}=0$ on $\Sigma^g\setminus \bigcup_{j\in J}\Omega_{k'}^j(\delta)$.  We see $\vec{w}^{\,\delta}_{i,k}$ as vectors in ${\R}^Q$ and we denote by $\pi^j_{k'}$ the map from $S^j_\infty\setminus \cup_{l=1}^{N^j} B_\delta(a^{j,l})$ into the space of projection matrices which to $x\in $ assigns the orthogonal projection from $T_{\vec{\Psi}^j_{k'}(x)}{\R}^Q$ into $T_{\vec{\Psi}^j_{k'}(x)}M^m$. In other words, let $P_z$ to be the $C^1$ map from $M^m$ into the space of $Q\times Q$ matrices which assigns the orthogonal projection onto $T_zM^m$, we have $\pi^j_{k'}(x):=P_{\vec{\Psi}^j_{k'}(x)}$ and we have
\be
\label{A.I.16}
\pi^j_{k'}\longrightarrow P_{\vec{\Psi}_\infty}\quad\quad\mbox{ strongly in }W^{1,2}_{loc}(S^j_\infty\setminus \cup_{l=1}^{N^j} B_\delta(a^{j,l}))
\ee
  On $S^j_\infty\setminus \cup_{l=1}^{N^j} B_\delta(a^{j,l})$  we denote  $\vec{u}_{i,k'}^{\,\delta}(x):=\pi^j_{k}(x)(\vec{w}_{i,k}^{\,\delta})$. Because of (\ref{A.I.16}) we have
 \be
\label{A.I.17}
  \vec{u}_{i,k}^{\,\delta}\longrightarrow \vec{w}_i^{\,\delta}\quad\quad\mbox{ strongly in }W^{1,2}(S^j_\infty)\quad.
\ee
Consider now the symmetric matrix
  \[
\begin{array}{l}
\ds D^2 \mbox{Area}(\vec{\Phi}_{k})(\vec{u}_{i,k}^{\,\delta},\vec{u}_{i',k}^{\,\delta})=\\[3mm]
\ds\quad\quad\sum_{j=1}^{\mbox{card}(J)}\int_{S^j_\infty} \lf[\lf<d\vec{u}_{i,k}^{\,\delta}\,;\,d\vec{u}_{i',k}^{\,\delta}\rg>_{g_{\vec{\Psi}^j_{k}}}+\lf<d\vec{\Psi}^j_{k}\,;\,d\vec{u}^{\,\delta}_{i,k}\rg>_{g_{\vec{\Psi}_{k}^j}}\lf<d\vec{\Psi}_{k}^j\,;\,d\vec{u}_{i',k}^{\,\delta}\rg>_{g_{\vec{\Psi}_{k}^j}} \rg]\ d\mbox{vol}_{g_{\vec{\Psi}_{k}^j}}\\[3mm]
\ds- \ 2^{-1} \sum_{j=1}^{\mbox{card}(J)}\int_{S^j_\infty} \lf<d\vec{\Psi}_{k}^j\dot{\otimes}\, d\vec{u}_{i,k}^{\,\delta}+d\vec{u}_{i,k}^{\,\delta}\dot{\otimes}\, d\vec{\Psi}_{k}^j, d\vec{\Psi}_{k}^j\dot{\otimes}\, d\vec{u}_{i',k}^{\,\delta}+d\vec{u}^{\,\delta}_{i',k}\dot{\otimes}\, d\vec{\Psi}_{k}^j\rg>\ d\mbox{vol}_{g_{\vec{\Psi}_{k}^j}}
\end{array}
\]
Let $f$ and $g$ be two smooth functions supported on $\vec{\Psi}_\infty(S_\infty^j\setminus\cup_{l=1}^{N^j} B_\delta(a^{j,l}) )$ then one has 
\[
\int_{S^j_\infty}<d(f(\vec{\Psi}_k)), d(g(\vec{\Psi}_k))>_{g_{\vec{\Psi}_{k}}}\ d\mbox{vol}_{g_{\vec{\Psi}_{k}}}=\int_{S^j_\infty}<d(f(\vec{\Psi}_k)), d(g(\vec{\Psi}_k))>_{h^j_k}\ d\mbox{vol}_{h^j_k}
\]
And since $h_{k}^j$ converges in any norms towards $h^j_\infty$, because of the strong $W^{1,2}$ convergence of $\vec{\Psi}_k$ on $S_\infty^j\setminus\cup_{l=1}^{N^j} B_\delta(a^{j,l})$ one has
\be
\label{A.I.18}
\int_{S^j_\infty}<d(f(\vec{\Psi}_k)), d(g(\vec{\Psi}_k))>_{g_{\vec{\Psi}_{k}}}\ d\mbox{vol}_{g_{\vec{\Psi}_{k}}}\longrightarrow \int_{S^j_\infty}<d(f(\vec{\Psi}_\infty)), d(g(\vec{\Psi}_\infty))>_{g_{\vec{\Psi}_{\infty}}}\ d\mbox{vol}_{g_{\vec{\Psi}_{\infty}}}
\ee
In a conformal chart  for $h_{k}^j$ we denote $e^{\la^j_{k'}}:=|\p_{x_1}\vec{\Psi}_{k'}^j|=|\p_{x_2}\vec{\Psi}_{k'}^j|$. Because of the strong $W^{1,2}$ convergence (\ref{A.I.15}) we have
\[
e^{\la^j_{k'}}\longrightarrow e^{\la^j_\infty}=|\p_{x_1}\vec{\Psi}_\infty|=|\p_{x_2}\vec{\Psi}_\infty|\quad\mbox{ a. e. in }\quad S^j_\infty\quad.
\]
Since $e^{\la^j_\infty}>0$ almost everywhere on $S^j_\infty$ we have $e^{-\la^j_{k'}}\longrightarrow e^{-\la^j_\infty}$ almost everywhere and then for $i=1,2$
\[
\p_{x_i}\vec{\Psi}^j_{k}/e^{\la^j_{k}}\longrightarrow \p_{x_i}\vec{\Psi}^j_{\infty}/e^{\la^j_{\infty}}\quad\mbox{ almost everywhere}
\]
Let $f$ ,$g$, $\phi$ and $\psi$ be 4 smooth functions on $M^m$ where $f$ and $g$ are supported on  $\vec{\Psi}_\infty(S_\infty^j\setminus\cup_{l=1}^{N^j} B_\delta(a^{j,l}) )$  one has in local conformal coordinates
\[
\begin{array}{l}
\ds<d(f(\vec{\Psi}^j_k))\otimes d(\phi(\vec{\Psi}^j_k)), d(g(\vec{\Psi}^j_k))\otimes d(\psi(\vec{\Psi}^j_k))>_{g_{\vec{\Psi}^j_{k}}}\ d\mbox{vol}_{g_{\vec{\Psi}^j_{k}}}=\\[3mm]
\ds\sum_{\mu,\nu=1,2} e^{-2\la^j_k} \p_{x_\mu}f(\vec{\Psi}^j_k)\ \p_{x_\nu}\phi(\vec{\Psi}^j_k)\ \p_{x_\mu}g(\vec{\Psi}^j_k)\ \p_{x_\nu}\psi(\vec{\Psi}^j_k)\ dx_1\wedge dx_2
\end{array}
\]
Because of the above
\[
e^{-2\la^j_k} \p_{x_\mu}f(\vec{\Psi}^j_k)\ \p_{x_\nu}\phi(\vec{\Psi}^j_k)\ \p_{x_\mu}g(\vec{\Psi}^j_k)\ \p_{x_\nu}\psi(\vec{\Psi}^j_k)\longrightarrow 
e^{-2\la^j_\infty} \p_{x_\mu}f(\vec{\Psi}^j_\infty)\ \p_{x_\nu}\phi(\vec{\Psi}^j_\infty)\ \p_{x_\mu}g(\vec{\Psi}^j_\infty)\ \p_{x_\nu}\psi(\vec{\Psi}^j_\infty)
\]
almost everywhere and we have moreover
\[
|e^{-2\la^j_\infty} \p_{x_\mu}f(\vec{\Psi}^j_\infty)\ \p_{x_\nu}\phi(\vec{\Psi}^j_\infty)\ \p_{x_\mu}g(\vec{\Psi}^j_\infty)\ \p_{x_\nu}\psi(\vec{\Psi}^j_\infty)|\le C\ |\nabla\vec{\Psi}^j_k|^2\rightarrow |\nabla\vec{\Psi}^j_\infty|^2\quad\mbox{strongly in }L^1
\]
Hence the generalized dominated convergence theorem implies
\[
\begin{array}{l}
\ds\int_{S^j_\infty}<d(f(\vec{\Psi}^j_k))\otimes d(\phi(\vec{\Psi}^j_k)), d(g(\vec{\Psi}^j_k))\otimes d(\psi(\vec{\Psi}^j_k))>_{g_{\vec{\Psi}^j_{k}}}\ d\mbox{vol}_{g_{\vec{\Psi}^j_{k}}}\\[3mm]
\ds\longrightarrow\quad\int_{S^j_\infty}<d(f(\vec{\Psi}^j_\infty))\otimes d(\phi(\vec{\Psi}^j_\infty)), d(g(\vec{\Psi}^j_\infty))\otimes d(\psi(\vec{\Psi}^j_\infty))>_{g_{\vec{\Psi}^j_{\infty}}}\ d\mbox{vol}_{g_{\vec{\Psi}^j_{\infty}}}
\end{array}
\]
Similarly we also have
\[
\begin{array}{l}
\ds\int_{S^j_\infty}\lf<d(f(\vec{\Psi}^j_k)), d(g(\vec{\Psi}^j_k))\rg>_{g_{\vec{\Psi}_k}}\ \lf<d(\phi(\vec{\Psi}^j_k)), d(\psi(\vec{\Psi}^j_k))\rg>_{g_{\vec{\Psi}_k}}  \ d\mbox{vol}_{g_{\vec{\Psi}^j_{k}}}\\[3mm]
\ds\longrightarrow\quad\int_{S^j_\infty}\lf<d(f(\vec{\Psi}^j_\infty)), d(g(\vec{\Psi}^j_\infty))\rg>_{g_{\vec{\Psi}_\infty}}\ \lf<d(\phi(\vec{\Psi}^j_\infty)), d(\psi(\vec{\Psi}^j_\infty))\rg>_{g_{\vec{\Psi}_\infty}}  \ d\mbox{vol}_{g_{\vec{\Psi}^j_{\infty}}}
\end{array}
\]
Combining all the above gives
\be
\label{A.I.19}
 D^2 \mbox{Area}(\vec{\Phi}_{k})(\vec{u}_{i,k}^{\,\delta},\vec{u}_{i',k}^{\,\delta})\quad\longrightarrow\quad D^2 \mbox{Area}(\vec{\Phi}_{\infty})(\vec{w}_{i}^{\,\delta},\vec{w}_{i'}^{\,\delta})
\ee
Hence, for $k$ large enough $(D^2 A(\vec{\Phi}_{k})(\vec{u}_{i,k}^{\,\delta},\vec{u}_{i',k}^{\,\delta}))_{i,i'=1\cdots N}$ defines a strictly negative quadratic form.

\medskip

Using now lemma~\ref{deuxi} below we deduce that for any $i,i'\in \{1\cdots N\}$
\be
\label{A.I.20}
\sigma_k^2 \lf|D^2 F(\vec{\Phi}_{k})(\vec{u}_{i,k}^{\,\delta},\vec{u}_{i',k}^{\,\delta})\rg|\le C\sigma_k^2 \lf[ F(\vec{\Phi}_{k})+\mbox{Area}(\vec{\Phi}_{k})^{1/4}\ F(\vec{\Phi}_{k})^{3/4}\rg]=o(1)
\ee
Combining (\ref{A.I.19}) and (\ref{A.I.20}) we obtain that for $k$ large enough $(D^2 A^{\sigma_k}(\vec{\Phi}_{k})(\vec{u}_{i,k}^{\,\delta},\vec{u}_{i',k}^{\,\delta}))_{i,i'=1\cdots N}$ defines a strictly negative quadratic form. This implies inequality (\ref{A.I.14}) and theorem~\ref{th-A.3} is proved.\hfill $\Box$

\section{Minmax Hierarchies for the area}
\reset
\subsection{Definition of the ${\mathbf F}-$distance on ${\mathfrak M}(M^m)$.}

For any immersion $\vec{\Phi}\in \mbox{Imm}(\Sigma^g,M^m)$ we can define the corresponding {\it oriented varifold} in $M^m$ as follows
\[
\forall \varphi\in C^0(\ti{G}_2(M^m))\quad V_{\vec{\Phi}}(\varphi):=\int_{\Sigma^g}\varphi\lf(\vec{\Phi}(x),\vec{\Phi}_\ast T_x\Sigma^g\rg)\ dvol_{g_{\vec{\Phi}}}
\]
where $\ti{G}_2(M^m)$ is the Grassmann bundle of oriented 2 planes in $TM^m$ over $M^m$. When $\varphi$ is just a function in $M^m$ we keep denoting
\[
V_{\vec{\Phi}}(\varphi):=\int_{\Sigma^g}\varphi\lf(\vec{\Phi}(x)\rg)\ dvol_{g_{\vec{\Phi}}}
\]
We call the {\it ${\mathbf F}-$distance} between 2 immersions $\vec{\Phi}$ and $\vec{\Psi}$ of respectively 2 oriented closed surfaces $\Sigma^g$ and $\Sigma^h$
\[
{\mathbf F}(\vec{\Phi},\vec{\Psi}):=\sup_{\|\varphi\|_{Lip}\le 1}V_{\vec{\Phi}}(\varphi)-V_{\vec{\Psi}}(\varphi)+{\mathcal F}\lf(\vec{\Phi}_\ast[\Sigma^g]-\vec{\Psi}_\ast[\Sigma^h]\rg)
\]
where ${\mathcal F}$ is the usual Flat norms between 2-cycles and $\vec{\Phi}_\ast[\Sigma^g]$ and $\vec{\Psi}_\ast[\Sigma^h]$ denote respectively
the push forwards  by $\vec{\Phi}$ and by $\vec{\Psi}$ of the currents of integration along respectively $\Sigma^g$ and  $\Sigma^h$. Observe that $V$ and ${\mathbf F}$ are independent of the oriented parametrization and hence these two functions ``descend'' to ${\mathfrak M}(M^m)$. We shall
also use the following function which defines a distance on ${\mathfrak M}(M^m)$.
\[
\ti{\mathbf F}(\vec{\Phi},\vec{\Psi}):={\mathbf F}(\vec{\Phi},\vec{\Psi})+|F(\vec{\Phi})-F(\vec{\Psi})|\quad.
\]
Given a metric space $(X,d)$ we say that a map $\vec{\Phi}$ is ${\mathbf F}-$Lipschitz from $X$ into ${\mathfrak M}(M^m)$ if it satisfies
\[
\exists \ K>0\quad\mbox{s.t.}\quad \forall\ x,y\in X\quad{\mathbf F}(\vec{\Phi}(x),\vec{\Phi}(y))\le K\ d(x,y)
\]
The space of ${\mathbf F}-$lipschitz maps from $X$ into ${\mathfrak M}(M^m)$ is denoted
\[
\mbox{Lip}_{\mathbf F}\lf(X,{\mathfrak M}(M^m)\rg)\quad.
\]
\subsection{The viscosity Limits of Minmax Hierarchies for the Area.}

We shall first introduce a modified hierarchy ${W}_k^\sigma$ for $A^\sigma$ which limit will be the hierarchy for the area. In this sub-section we are assuming a genus bound $g\le g_0$
and we are working in ${\mathfrak M}^{g_0}(M^m)$.

The first element in a hierarchy can be given by $n_1=N_1=m-2$ and is the same as in the previous subsection.
\[
\mbox{Sweep}_{N_1}(M^m):=\lf\{   
\begin{array}{c}
\ds (Y,\vec{\Phi})\ ;\ Y\in{\mathcal P}_{m-2}\quad \vec{\Phi}\in \mbox{Lip}(Y,{\mathfrak M}^{g_0}(M^m))\\[3mm]
\ds\vec{\Phi}_\ast\lf[Y\times \Sigma_{\vec{\Phi}}\rg]\quad\mbox{ generates }H_m(M^m,{\Z})\\[3mm]
\forall\ y\in \p Y\quad\quad A^\sigma(\vec{\Phi}(y))<W_1(M^m)
\end{array}\rg\}
\]
where $W_1(M^m)$ is the usual {\it Width } of a closed manifold that is the minimal maximal area needed to {\it sweep-out} $M^m$. Denote for $\sigma>0$
\[
W^\sigma_1(M^m):=\inf_{(Y,\vec{\Phi})\in \mbox{Sweep}^0_{N_1}(S_0,M^m)}\quad\quad\max_{y\in Y_0}\ A^\sigma(\vec{\Phi}(t))
\]
where
\[
Y_0:=\lf\{y\in Y\ ; \mbox{Area}(\vec{\Phi}(y))>2^{-1}\,W^0_0\rg\}
\]
and we will simply write $W_l(M^m)$ for $W^0_l(M^m)$. Assuming now the hierarchy is constructed up to the order $k-1$, we introduce the notation for $l=1\cdots k-1$
\[
{\mathcal C}_l:=\lf\{ T\in {\mathcal Z}_2(M^m)\ ;\ \ \mbox{Area}(T)={W}_l(M^m)\quad\mbox{and }\quad T\mbox{ is a current of integration on a minimal surface }\rg\}
\]
We are going to make the following assumption
\[
\mbox{(H1)}\quad\quad\quad {\mathcal C}_l \mbox{ \it is a smooth compact sub-manifold of }{\mathcal Z}_2(M^m)\mbox{\it equipped with the flat norm \footnote{This should be satisfied for generic $M^m$ in ${\R}^Q$ for $m=3$ - see \cite{Wh1}.}. }
 \] 
 For any $\ep>0$ we denote
\[
{\mathcal O}_l(\ep):=\lf\{\vec{\Phi}\in {\mathfrak M}^{g_0}\ ;\ d_{\mathbf F}(\vec{\Phi},{\mathcal C}_l)<\ep\rg\}
\]
Let $\ep_l>0$ be fixed such that 
\be
\label{I.6-b}
\exists\ \pi_{l}\ \in\, \mbox{Lip}_{\mathcal F}\lf( {\mathcal O}_l(\ep_l) ,{\mathcal C}_l\rg)\quad\mbox{ s. t. }\quad \forall \ \vec{\Phi}\in {\mathcal C}_l\quad\quad\pi_l(\vec{\Phi})=\vec{\Phi}
\ee
as given by \cite{Lan}. The tubular neighborhood ${\mathcal O}_l(\ep_l) $ of ${\mathcal C}_l$ for the ${\mathbf F}-$distance will be denoted ${\mathcal O}_l$. 
We shall denote
\[
\mbox{Sweep}_{N_l}(M^m)_0:=\lf\{\vec{\Phi}\in\mbox{Sweep}_{N_l}(M^m)\quad;\quad A^{\tau_l}(\vec{\Phi})\le W^{\tau_l}_l(M^m)+\delta_l\rg\}
\]
where $\tau_l$ and $\delta_l$ are given by theorem~\ref{th-A.2} for $\ep:=\ep_l/2$ and are characterized by the following assertion
\be
\label{I.6-c}
\forall\, (Y,\vec{\Phi})\in\mbox{Sweep}_{N_l}(M^m)\quad\quad A^{\tau_l}(\vec{\Phi})\le W^{\tau_l}_l(M^m)+\delta_l\quad\Longrightarrow\quad d_{{\mathcal F}}(\vec{\Phi}(Y),{\mathcal C}_l)<2^{-1}\,\ep_l
\ee

\medskip

We define $n_k$  as follows. Let $n_k\in{\N}^\ast$ such that
\[
H^{n_k-1}({\mathcal C}_{k-1},{\Z}_2)\ne 0
\]
and choose $\om_{k-1}$ being a non zero element of $H^{n_k-1}({\mathcal C}_{k-1},{\Z}_2)$. 

\medskip

Under the previous notations we define $\mbox{Sweep}_{N_k}(M^m)$ to be the set of pairs $(Y,\vec{\Phi})$ such that

\begin{itemize}
\item[i)] 
\[
Y\in {\mathcal P}_{N_{k}}\quad,\quad \vec{\Phi}\in\mbox{Lip}(\ov{Y},{\mathfrak M}^{g_0})
\]
\item[ii)] There exists $Z\in {\mathcal P}_{N_{k-1}}$ s.t.
\[
\p Y=\p\lf(B^{n_k}\times Z\rg)
\]
\item[iii)] there exists
\[
\vec{\Psi}\in \mbox{Lip}_{\mathbf F}\lf(\p B^{n_k}, \mbox{Sweep}_{N_{k-1}}(M^m)_{0}\rg)
\]
and
\[
\max_{x\in \p B^{n_k}} d_{\mathbf F}(\vec{\Phi}(x,\cdot),\vec{\Psi}(x))<2^{-1}\,\ep_{k-1}
\]
\item[iv)] Let 
\[
\Om_{\vec{\Phi}}:=\lf\{  y\in \p Y\ ;\ d_{\mathbf F}(\vec{\Phi}(y),{\mathcal C}_{k-1})<\ep_{k-1} \rg\}
\]
we have

\[
[\Om_{\vec{\Phi}}\cap (\{x\}\times Z)]\in H_{N_{k-1}}(\,\Om_{\vec{\Phi}},\p\,\Om_{\vec{\Phi}} ,{\Z}_2)\quad\mbox{ is {\bf Poincar\'e dual} to } (\pi_{k-1}\circ\vec{\Phi})^\ast\om_{k-1}\in H^{n_k-1}(\Om_{\vec{\Phi}},{\Z}_2)
\]
\item[v)]  We have $\forall \, x\in B^{n_k}$ 
\[
\limsup_{z\rightarrow \p Z}\mbox{Area}(\vec{\Phi}(x,z))<W_{k-1}-\delta_{k-1}
\]
\end{itemize}
We shall consider the following important property
\[
(\mbox{P}_k)\quad\quad\forall\, \ep>0\quad\exists\,\delta>0\quad\quad\forall A\in \mbox{Sweep}_{N_{k}}(M^m)\quad\max_{\vec{\Phi}\in A}\mbox{Area}(\vec{\Phi})<W_k(M^m)+\delta\quad\Longrightarrow d_{\mathbf F}(A,{\mathcal C}_k)<\ep
\]
The following proposition holds
\begin{Th}
\label{pr-II.2}
Assume $(\mbox{P}_l)$ holds for any $l\le k$ then
\be
\label{I.8}
W_{k}(M^m)<W_{k+1}(M^m)
\ee
and for $\ep_{k}>0$, $\tau_k>0$ and $\delta_k>0$ satisfying (\ref{I.6-b}) and (\ref{I.6-c}) chosen small enough in the definition of $\mbox{Sweep}_{N_k}(M^m)_0$ then $\mbox{Sweep}_{N_{k+1}}(M^m)$ is an admissible family that is to say : There exists $\eta>0$ such that
 for any homeomorphism $\Xi$ of ${\mathfrak M}^{g_0}$ satisfying
\[
\forall\,\vec{\Phi}\in {\mathfrak M}^{g_0}\quad\mbox{ s. t. }\quad\mbox{Area}(\vec{\Phi})\le W_{k}(M^m)-\eta\quad \Longrightarrow\quad\Xi(\vec{\Phi})=\vec{\Phi}
\]
then
\[
\Xi\lf(\mbox{Sweep}_{N_{k+1}}(M^m)\rg)\subset \mbox{Sweep}_{N_{k+1}}(M^m)\quad.
\]
Thus there exists a parametrized interger rectifiable stationary varifold ${\mathbf v}_{k+1}:=(S_{k+1},\vec{\Psi}_{k+1},\theta_{k+1})$  such that 
\[
\mbox{genus}\,(\Sigma_{\vec{\Psi}_{k+1}})\le g_0\quad,\quad\mbox{Area}\,(\vec{\Psi})=W_{k+1}(M^m)
\]
Assuming
\[
W_{k+1}(M^m)<\frac{8\pi}{\|\vec{\mathbb I}_{M^m}\|_\infty}
\]
then $\theta_{k+1}\equiv 1$, $\vec{\Psi}_{k+1}$ is a smooth conformal minimal embedding and
\[
\quad\mbox{Ind}\,(\vec{\Psi}_{k+1})\le N_{k+1}\quad.
\]
\hfill $\Box$
\end{Th}
Before proving the previous theorem we observe that the main result of \cite{Pitts} gives
\begin{Th}
\label{th-Pitts}
For $m=3$ the hypothesis $(\mbox{P}_k)$ holds for all $k$.
\hfill $\Box$
\end{Th}
\begin{Rm}
\label{rm-I.1}
In concrete situations (the next section will be illustrating that fact) the most delicate part in going from $k-1$ to $k$ in the construction of the hierarchy
is not really to fix $\om_k$ and $n_k$ (there are sometimes several choices) or to construct $\vec{\Phi}_x$ all along the boundary of $Y$ in such a way that iv) holds. The difficulty is more to find $Y$,  and some extension $\vec{\Phi}$ inside $Y$. The freedom to change the topology of the surface and to pass from a given surface at the level $k-1$ to another surface, usually of higher genus, at the level $k$ is given in order to ease the construction of $\vec{\Phi}$ inside $Y$ as it is illustrated again in the next section.
\end{Rm}
\begin{Rm}
\label{rm-I.2}
In generic situations one could expect ${\mathcal C}_l$ to be made of isolated points (see the main result of \cite{Wh1}). It contains in particular one minimal surface and the same one with opposite orientation. Then one would take $n_l=1$ and $\p B^1=\{-1\}\cup\{+1\}$ and we would be considering the boundary data to be on $Z\times\{-1\}$ an arbitrary element of $\mbox{Sweep}_{N_l}(M^m)$ satisfying (\ref{I.6-c}) and at $Z\times\{+1\}$ we would take exactly the same element of  $\mbox{Sweep}_{N_l}(M^m)$ but with opposite orientation. The elements in $\mbox{Sweep}_{N_{l+1}}(M^m)$ would then consist of {\bf eversions} of elements in $\mbox{Sweep}_{N_l}(M^m)$.
\end{Rm}

\medskip

\noindent{\bf Proof of theorem~\ref{pr-II.2}.} The fact that property $(\mbox{P}_{k-1})$ implies (\ref{I.8}) is established by following word by word the arguments in the
first part of the proof of theorem~\ref{th-I.1}. 

Let
\[
W^\sigma_k(M^m):=\inf_{(Y,\vec{\Phi})\in \mbox{Sweep}_{N_{k}}(M^m)}\max_{y\in Y}A^\sigma(\vec{\Phi})\quad.
\]
and we shall simply use the notation $W_k(M^m)$ for $W_k^0(M^m)$. It is clear that $W_k^\sigma(M^m)$ is an increasing function of $\sigma$ and $k$. Moreover
\be
\label{I.7}
\lim_{\sigma\rightarrow 0}W_k^\sigma(M^m)=W_k(M^m)\quad.
\ee
We assume that $\delta_{k-1}$ and $\tau_{k-1}$ have been chosen small enough in the definition of $\mbox{Sweep}_{N_{k-1}}(M^m)_0$ in such a way that
\be
 \label{I.7.0}
 \forall x\in \p B^{n_k}\quad\max_{z\in Z} \mbox{Area}(\vec{\Psi}(x,z))\le W_{k-1}(M^m)+4^{-1}\,(W_k(M^m)-W_{k-1}(M^m))
 \ee
This implies for a fixed $\vec{\Phi}\in \mbox{Sweep}_{N_{k}}(M^m)$ and $\sigma$ small enough
\be
\label{I.7.1}
 \forall x\in \p B^{n_k}\quad\max_{z\in Z} A^\sigma(\vec{\Phi}(x,z))\le W_{k-1}(M^m)+2^{-1}\,(W_k(M^m)-W_{k-1}(M^m))
\ee
 We consider the subpart of $\mbox{Sweep}_{N_{k}}(M^m)$ such that the inequality (\ref{I.7.1}) holds (this can only increase the value of $W^\sigma_k(M^m)$ but not the limiting value of $W_k(M^m)$. Using  \cite{Riv-columb} (working with gradients multiplied by cut-offs above the $A^\sigma$ energy level  $W_{k-1}(M^m)+2^{-1}\,(W_k(M^m)-W_{k-1}(M^m))<W^\sigma_k(M^m)$) combined with corollary 10.5 of \cite{Gho} we obtain the existence of $\sigma_i\rightarrow 0$ and $\vec{\Phi}^{\sigma_i}\in {\mathfrak M}^{g_0}$ such that
\[
\lf\{
\begin{array}{l}
\ds W_k^{\sigma_i}(M^m)=A^{\sigma_i}(\vec{\Phi}^{\sigma_i})\quad,\quad  DA^{\sigma_i}(\vec{\Phi}^{\sigma_i})=0\quad\mbox{ and }\quad\mbox{Ind}^{\sigma_i}\,(\vec{\Phi}^{\sigma_i})\le N_k\\[3mm]
\ds (\sigma^i)^2\int_{\Sigma_{\vec{\Phi}^{\sigma_i}}}(1+|{\mathbb I}_{\vec{\Phi}^{\sigma_i}}|^2)^2\ dvol_{\vec{\Phi}^{\sigma_i}}=o\lf(\frac{1}{\log (\sigma^i)^{-1}}\rg)\quad,
\end{array}
\rg.
\]
where $\mbox{Ind}^{\sigma_i}\,(\vec{\Phi}^{\sigma_i})$ is the largest dimension of a sub-vector space of $W^{3,2}({(\vec{\Phi}^{\sigma_i})}^\ast TM^m)$ on which $D^2 A^{\sigma_i}(\vec{\Phi}^{\sigma^i})$ is strictly negative. 
Using now theorem~\ref{th-A.3} we obtain theorem~\ref{pr-II.2}.\hfill $\Box$
\section{An example of Hierarchies of Minmax Problems for the Area in $S^3$.}
\begin{Lm}
\label{lm-II.1}
For any $\vec{\Phi}\in \mbox{Imm}^0(\Sigma^g,S^3)$ there exists $\vec{\Phi}_{(t)}\in \mbox{Lip}_{\mathbf F}\lf((-1,1),\mbox{Imm}^0(\Sigma^g,S^3)\rg)$ such that
\be
\label{II.1}
\vec{\Phi}_\ast[(-1,1)\times \Sigma_g]=\om_{S^3}
\ee
 and 
 \be
 \label{II.2}
 \exists\ t_0\in (-1,1)\quad \vec{\Phi}_{t_0}=\vec{\Phi}
 \ee
 and 
 \be
\label{II.2-a}
\lim_{t\rightarrow \pm 1}|V_{\vec{\Phi}(t)}|(\ti{G}_2(S^3))=0
\ee
\hfill $\Box$
\end{Lm}
\noindent{\bf Proof of lemma~\ref{lm-II.1}.}
\begin{Dfi}
\label{df-II.2}
We call a sweep-out of $S^3$ by $\Sigma^g$ the space of $\vec{\Phi}_{(t)}\in \mbox{Lip}_{\mathbf F}\lf((-1,1),\mbox{Imm}^0(\Sigma^g,S^3)\rg)$ satisfying (\ref{II.1}) and (\ref{II.2-a}).
We denote this space by $\mbox{Sweep}_1(\Sigma^g,S^3)$. For any $\vec{\Phi}\in \mbox{Imm}^0(\Sigma^g,S^3)$, the family of  $\vec{\Phi}_{(t)}\in \mbox{Lip}_{\mathbf F}\lf((-1,1),\mbox{Imm}^0(\Sigma^g,S^3)\rg)$ satisfying (\ref{II.1}), (\ref{II.2}) and (\ref{II.2-a}) is denoted $\mbox{Sweep}_1(\vec{\Phi})$ and is equipped with the metric issued by ${\mathbf F}$
\[
{\mathbf F}(\vec{\Phi}(\cdot),\vec{\Psi}(\cdot)):=\sup_{t\in(-1,1)}{\mathbf F}(\vec{\Phi}(t),\vec{\Psi}(t)).
\] 
Observe that because of (\ref{II.2-a}) we have
\[
\lim_{t,t'\rightarrow \pm 1}{\mathbf F}(\vec{\Phi}(t),\vec{\Phi}(t'))=0
\]
\hfill $\Box$
\end{Dfi}
We denote for any $\vec{\Phi}\in \mbox{Imm}^0(\Sigma^g,S^3)$
\[
W_1(\vec{\Phi}):=\inf_{\vec{\Phi}_{(t)}\in \mbox{Sweep}_1(\vec{\Phi})}\ \max_{t\in(0,1)}\mbox{Area}(\vec{\Phi}_{(t)})
\]
From \cite{Lev} and \cite{Gro} (see also \cite{Ros}) combined with \cite{BDF} we deduce the following result.
\begin{Th}
\label{th-II.3}
We have
\be
\label{II.3}
W_1(S^3)=\inf_{\vec{\Phi}\in \mbox{Imm}^0(\Sigma^g,S^3)}W_1(\vec{\Phi})=4\pi
\ee
moreover $W_1(\vec{\Phi})=4\pi$ if and only if there exists $\vec{\Phi}_{(t)}\in \mbox{Sweep}_1(\vec{\Phi})$ and $t_0\in (-1,+1)$ such that
$\vec{\Phi}_{t_0}(\Sigma^g)$ is a geodesic 2-sphere. Denote ${\mathcal S}$ the space of geodesic oriented 2-spheres in $S^3$ identified
with their current of integration. We have moreover\footnote{We are using the fact that $|S^3|/2=\pi^2$.} 
\[
\begin{array}{l}
\ds\forall\, \ep>0\quad\exists\, \delta>0\quad\forall\, \vec{\Phi}\in\mbox{Imm}^0(\Sigma^g,S^3)\quad\forall\, t\in (-1,1)\quad\mbox{s.t. }\quad {\mathbf M}\lf(\vec{\Phi}_\ast[(-1,t)\times \Sigma_g]\rg)\in (-\delta+\pi^2,\pi^2+\delta)\\[3mm]
\ds\quad\mbox{then}\quad\quad\mbox{Area}(\vec{\Phi}_{(t)})\le W_1(S^3)+\delta\quad\Longrightarrow\quad {\mathbf F}\lf((\vec{\Phi}_{(t)})_\ast[\Sigma^g],S\rg)\le\ep\quad.
\end{array}
\]
for some $S\in {\mathcal S}$. \hfill $\Box$
\end{Th}
Observe that ${\mathcal S}$ canonically identifies to the Grassman manifold of {\it oriented } 3-planes in ${\R}^4$ :
\[
\tilde{G}_{3}({\R}^4)\simeq  S^3\quad.
\]
 For every $\vec{\Phi}_{(t)}\in\mbox{Sweep}_1(\vec{\Phi})$ and every $t\in(-1,1)$ we denote by
\[
\mbox{Vol}_{\vec{\Phi}}(t):={\mathbf M}\lf(\vec{\Phi}_\ast[(-1,t)\times \Sigma_g]\rg)
\]
and let $\om_{S^3}$ to be the class in $H^3(S^3,{\Z}_2)$ such that 
\[
<\om_{S^3},[S^3]>_{H^3,H_3}=+1\quad.
\]
With the above notations we have the following proposition.
\begin{Prop}
\label{pr-II.4}
Let $\ep>0$ such that there exists a lipschitz projection onto the space of geodesic sphere $ {\mathcal S}$ in the ${\mathcal F}-$neighborhood  of $ {\mathcal S}$. Denote $\delta_0>0$ the constant given by theorem~\ref{th-II.3}. For any $g>0$ there exists $\vec{\Phi}\in \mbox{Lip}_{\mathbf F}(B^4,\mbox{Sweep}_1(\Sigma^g,S^3))$ such that
\[
\sup_{y\in B^4\times \{\pm 1\}}W_1(\vec{\Phi}_{(y)})\le 4 \pi-\delta_0
\]
moreover for any $z\in\p B^4$
\[
\pi_{\mathcal S}^\ast\,\om_{S^3}\frown\lf[\mbox{Vol}_{\vec{\Phi}}^{\ \delta}\rg]=\lf[\mbox{Vol}^{\ \delta}_{\vec{\Phi}}\cap \{z\}\times(-1,+1)\rg]
\]
 for some $0<\delta<\delta_0$ where $\mbox{Vol}_{\vec{\Phi}}^{\ \delta}:=\mbox{Vol}^{-1}_{\vec{\Phi}}([\pi^2-\delta,\pi^2+\delta])$. An element $\vec{\Phi}\in \mbox{Lip}_{\mathbf F}(B^4,\mbox{Sweep}_1(\Sigma^g,S^3))$ satisfying all the above is called a $g-$sweep-out of order 5 in $S^3$. The space of of $g-$sweep-out of order 5 in $S^3$ is denoted by $\mbox{Sweep}_5(\Sigma^g,S^3)$. We have that
\[
\mbox{Sweep}_{\,5}(T^2,S^3)\ne \emptyset\quad.
\]
Moreover
\be
\label{in-2picarre}
\inf_{\vec{\Phi}\in \mbox{Sweep}_5(T^2,S^3)}\quad\max_{y\in (-1,1)\times B^4}\mbox{Area}_{S^3}\,(\vec{\Phi}(y))< 8\pi\quad.
\ee
\hfill $\Box$
\end{Prop}
Before proving this proposition in the following subsection we first shall make the link with the general definition of {\bf Minmax Hierarchies}
as defined in the previous section.
For any $Y\in {\mathcal P}_1$ and $\vec{\Phi}\in  \mbox{Lip}_{\mathbf F}(Y,Imm^0(\Sigma^g,S^3))$ one has 
\[
4\pi\le \max_{t\in Y}\mbox{Area}(\vec{\Phi}(t))
\]
therefore it suffices to restrict to domains $Y=(-1,+1)$ and the definition of $\mbox{Sweep}_1(S^3)$ above does not include the mention
of the domain $Y$. We have then
\[
n_1=1\quad\mbox{ and }\quad W^0_1(S^3)=4\pi
\]
 Then we have $n_2=4$ and ${\mathcal C}_2^0={\mathcal S}\simeq S^3$ and $\pi_2=\pi_{\mathcal S}$. We choose $\om_2=\om_{S^3}$ where $\om_{S^3}$ is the class in $H^3(S^3,{\Z}_2)$ such that 
\[
<\om_{S^3},[S^3]>_{H^3,H_3}\ne 0\quad.
\]
For every $\ep>0$ there exists $\delta>0$ such that
\[
\Xi_{\vec{\Phi}}:=\lf\{ (z,t)\in \p B^4\times (-1,+1)\quad;\quad d_{\mathbf F}(\vec{\Phi}(z,t),{\mathcal C}_1^0)<\ep\rg\}\subset \mbox{Vol}^{-1}_{\vec{\Phi}}([\pi^2-\delta,\pi^2+\delta])\quad.
\]
Hence we can use both remark~\ref{rm-I.3} and lemma~\ref{lm-I.4} to deduce that condition iv) is satisfied. Now, we have imposed $Y=B^4\times(-1,1)$ and $ Z=(-1,1)$. A-priori this is too restrictive in order to apply theorem~\ref{th-I.1}. However the argument applies because
$N_{k-1}=n_0=1$. Then the chain $U$ is one dimensional and can be modified to coincide with a segment and belong to the class
$\mbox{Sweep}_1(S^3)$. Thus we have
\be
\label{II.4}
4\pi=W^0_1(S^3)<W^0_2(S^3)=\inf_{\vec{\Phi}\in \mbox{Sweep}_{\,5}(T^2,S^3)}\ \max_{y\in B^4\times(-1,1)}\mbox{Area}\,(\vec{\Phi}(y))<8\pi
\ee
and $\mbox{Sweep}_1(S^2,S^3)$ and $\mbox{Sweep}_{\,5}(T^2,S^3)$ realizes a {\bf Minmax hierarchy}. We shall prove bellow that
\[
W^0_2(S^3)=2\pi^2\quad\mbox{ and }\quad{\mathcal C}^0_1=\{\mbox{Clifford Tori}\}\simeq {\mathbb C}{\mathbb P}^1\times {\mathbb C}{\mathbb P}^1
\]
We will then comment on the possibility how to extend this hierarchy to an uncountable one but we will first spend the next subsection to prove proposition~\ref{pr-II.4}.

\subsection{Explicit constructions of Sweep-outs of various orders of  $B^3$ and $S^3$ by surfaces of various topologies.}

Before proving the proposition above we are going to prove a series of intermediate lemma. Let ${\mathcal C}$ be the set of oriented large circles (i.e. closed geodesics) in $S^2$ :
\[
{\mathcal C}\simeq S^2
\]
and, in a small enough neighborhood of ${\mathcal C}$ for the Flat norm (i.e. ${\mathcal F}^{-1}(\cdot,{\mathcal C})<\ep_{{\mathcal C}}$) within the space of oriented 1 dimensional cycles ${\mathcal Z}_1(S^2)$, we 
consider a ${\mathcal F}-$lipschitz projection $\pi_{\mathcal C}$ . We denote by $SP_2(\mbox{Imm}^0(S^1,S^2))$ the 2-fold symmetric product of the space $\mbox{Imm}^0(S^1,S^2)$ of immersions of circles in $S^2$ isotopic\footnote{There are exactly two connected components of the immersions of $S^1$ in $S^2$, the one of the equator and the one of the eight.} to any large circle.

\begin{Lm}
\label{lm-II.5}
There exists a map in $\vec{\Xi}\in \mbox{Lip}_{\mathbf F}\lf(\ov{B^3},SP_2(\mbox{Imm}^0(S^1,S^2))\rg)$ such that 
\be
\label{be1}
\forall\, \sigma\in \p B^3\quad\quad {\mathcal F}\lf(\vec{\Xi}(\sigma)_\ast \lf([S^1]+[S^1]\rg),{\mathcal C}\rg)<\ep_{\mathcal C}
\ee
and the map
\be
\label{be2}
\sigma\in \p B^3\ \ \longrightarrow\ \pi_{\mathcal C}\lf(\vec{\Xi}(\sigma)_\ast \lf([S^1]+[S^1]\rg)\rg)\in {\mathcal C}\simeq S^2
\ee
has a degree $\pm 1$.
\hfill $\Box$
\end{Lm}
The existence of such a map is illustrated by figure 1 which is showing a ``slice'' of such a map restricted to the half disc $B^3\cup\{x_2=0, x_1>0\}$
and which is taken to be axially symmetric. It is important to observe that the two discs cannot be ``indexed'' continuously on the whole ball $B^3$ and that it is
necessary to take the symmetric product $SP_2(\mbox{Imm}^0(S^1,S^2))$ .
\begin{figure}
\begin{center}
\includegraphics[width=15cm,height=21cm]{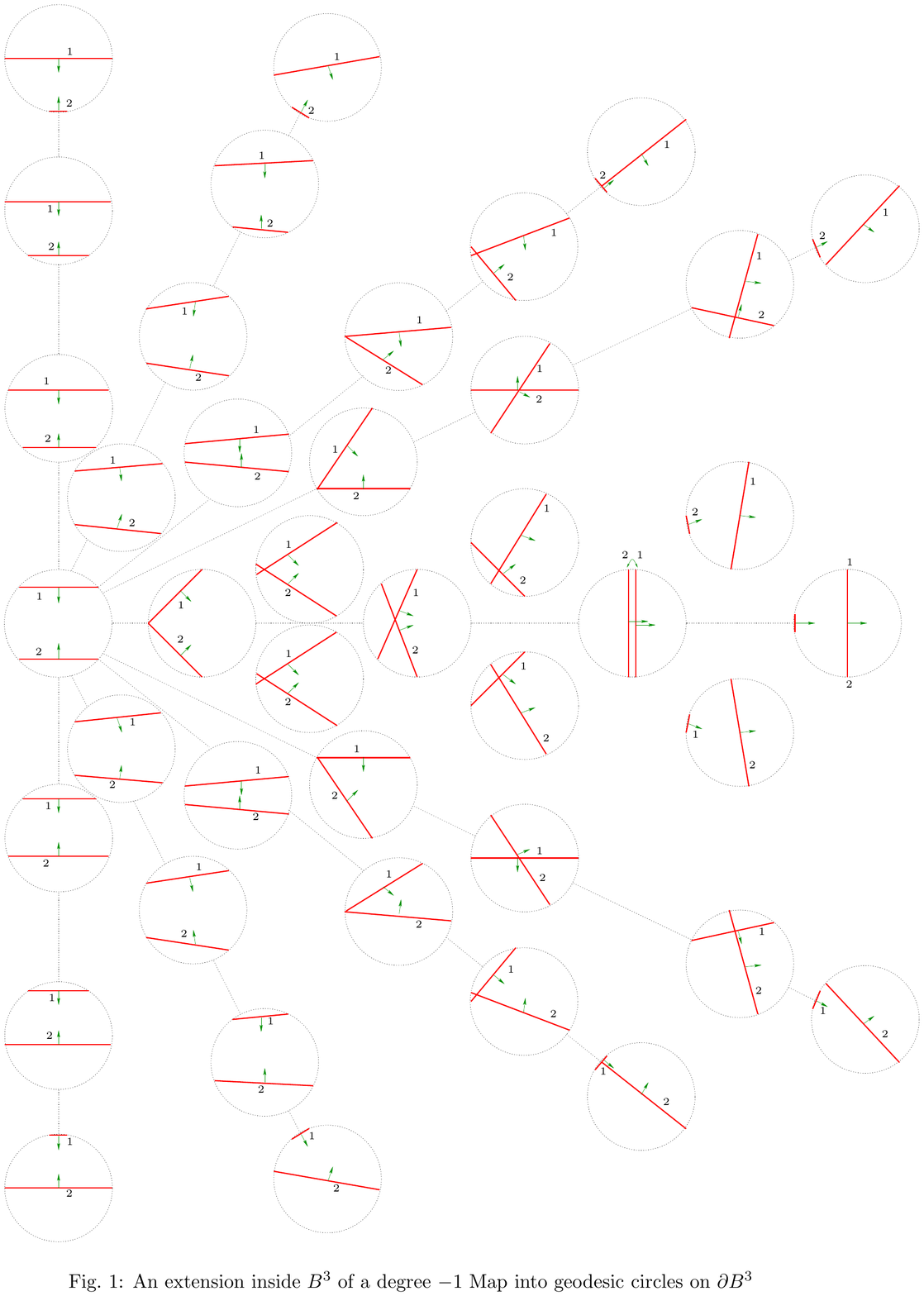}
\end{center}
\end{figure}
\medskip

\noindent{\bf Proof of lemma~\ref{lm-II.5}.} We shall restrict to circles in $S^2$ of constant curvatures (i.e. boundaries of geodesic balls in $S^2$). This space identifies canonically to ${\mathfrak C}:=S^2\times (-1,+1)$. Recall that the Symmetric product operation is an homotopy functor (see \cite{Hat} page 481) and hence, since ${\mathfrak C}$ is obviously homotopy equivalent to ${\mathcal C}\simeq S^2$ we have that
\[
SP_2({\mathfrak C})\simeq SP_2(S^2)\simeq {\C}{\mathbb P}^2
\]
The last homotopy equivalence (due to \cite{Lia}) and can be seen as follows. To a pair of points in $S^2\simeq {\C}{\mathbb P}^1\simeq {\mathbb C}\cup\{\infty\}$ we assign the family of order 2 polynomials having these points as roots. The space of order two polynomials identifies to $C{\mathbb P}^2$
\[
\forall\, \{a_+,a_-\}\in SP_2({\mathbb C})\quad\quad\varphi(\lf\{a_+,a_-\rg\})=[1, a_++a_-,a_+a_-]\in {\C}{\mathbb P}^2
\]
and $\varphi(a_+,\infty)=[0,{1}/{a_+},1]$ and $\varphi(\infty,\infty)=[0,0,1]$. On $\p B^3\simeq {\mathbb C}\cup\{\infty\}$ we take
\[
\vec{\Xi}(\sigma):=\{(\ov{\sigma},0),(\sigma,1)\}
\]
Observe that the map
\[
\varphi\circ\vec{\Xi}(\sigma)=[1,2\,\Re(\sigma), |\sigma|^2]
\]
is pulling back the K\"ahler form (that we express in the chart $[z_0,z_1,z_2]\rightarrow (z_1/z_0,z_2/z_0)=(w_1,w_2)$
\[
\om_{C{\mathbb P}^2}=\frac{i}{2\pi |w|^4}\lf(|w|^2 dw_1\wedge d\ov{w_1}+|w|^2 dw_2\wedge d\ov{w_2}-\ov{w_1}\, w_2\ dw_1\wedge d\ov{w_2}-\ov{w_2}\, w_1\ dw_2\wedge d\ov{w_1}   \rg)
\]
to zero. Hence $\varphi\circ\vec{\Phi}$ is homotopically trivial and can be extended continuously throughout $B^3$ and this proves the lemma.

\medskip

An alternative proof of lemma~\ref{lm-II.5} is also given by the following explicit example : We assume that $\vec{\Xi}$ is axially symmetric in the following sense. Writing 
\[
\vec{\Xi}(x,y,z)=\lf\{  (\sigma_1,t_1),(\sigma_2,t_2) \rg\}
\]
we choose first
\be
\label{xi-1}
\vec{\Xi}(R_\phi^z(x,y,z)):=\lf\{(R^z_\phi(\sigma_1),t_1), (R^z_\phi(\sigma_2),t_2)\rg\}
\ee
where $R^z_\phi$ is the rotation of axis $Oz$ and angle $\phi$. Then we choose 
\begin{itemize}
\item[i)] For $\theta\in [-\pi/2,\pi/2]$ and $0<r<1/2$ we set
\be
\label{xi-2}
\vec{\Xi}(r\cos\theta,0,r\sin\theta):=\lf\{((\cos\theta_r,0,\sin\theta_r), -1/2), (\cos\theta_r,0,-\sin\theta_r), -1/2)\rg\}
\ee
where $\theta_r:=(1-2r)\pi/2+2\, r\,|\theta|$
\item[ii)] For $\theta\in [-\pi/2,\pi/2]$ and $1/2\le r<1$ we set
\be
\label{xi-3}
\vec{\Xi}(r\cos\theta,0,r\sin\theta):=\lf\{((\cos\theta,0,\sin\theta), -r), (\cos\theta,0,-\sin\theta), -1+r)\rg\}
\ee
\end{itemize}
\hfill $\Box$

\medskip

We introduce now the space
\[
SP_2^\ast(\mbox{Imm}^0(S^1,S^2)):=\lf\{ \{\gamma,\gamma'\}\in SP_2(\mbox{Imm}^0(S^1,S^2))\quad;\quad \gamma(S^1)\cap\gamma'(S^1)=\emptyset\rg\}\quad.
\]

Let $A$ be the following annulus  $A:=B_1\setminus B_{1/2}$ and we denote by $\mbox{Imm}_T^0(A, B^3)$ the space of immersions of the annulus $A$ which are sending the boundary
of $A$ transversally to the boundary of $B^3$ and regular homotopic in this class to the intersection of $B^3$ with the catenoid $(\cosh s\,\cos\theta,\,\cosh s\, \sin\theta,  s)$. We shall now deform the element $\vec{\Xi}:=\{\gamma,\gamma'\}\in SP^\ast_2(\mbox{Imm}^0(S^1,S^2))$ constructed in the proof of lemma~\ref{lm-II.5}  into a new map $\vec{\Xi}_\ast$
taking values into $SP_2^\ast(\mbox{Imm}^0(S^1,S^2))$ this time while preserving the main properties (\ref{be1}) and (\ref{be2}). Moreover we are going to assign to this map a map $\vec{\Psi}\in \mbox{Lip}_{\mathbf F}(B^3, \mbox{Imm}_T^0(A, B^3))$ whose boundary is $\vec{\Xi}_\ast$. The set of oriented
flat discs bounding a large circle in $S^2$ will be denoted ${\mathcal D}$ and identifies to the space ${\mathcal C}\simeq S^2$ above.
We call a 1 {\it sweep-out} of $B^3$ by $A$ a map $\vec{\Phi}_{(t)}\in \mbox{Lip}_{\mathbf F}((-1,+1),\mbox{Imm}_T^0(A, B^3))$ such that
\[
\vec{\Phi}_\ast\lf([(-1,+1)\times A]\rg)\quad\mbox{generates }H_3(B^3,\p B^3,{\Z})
\]
This space is denoted $\mbox{Sweep}_1(A,B^3)$.

\begin{Dfi}
\label{df-II.5e}
The space of maps $\vec{\Phi}\in \mbox{Lip}_{\mathbf F}(B^3, \mbox{Sweep}_1(A,B^3))\subset \mbox{Lip}_{\mathbf F}(B^3\times(-1,+1), \mbox{Imm}_T^0(A,B^3)) $ such that 
there exists $0<\ep<\ep_{\mathcal D}$  and all $z\in \p B^3$ satisfying
 \[
 (\pi_{\mathcal D}\circ\vec{\Phi})^\ast\om_{\mathcal D}\frown [\Om_{\vec{\Phi}}]=[\Om_{\vec{\Phi}}\cap\{z\}\times (-1,+1)]\quad\quad\mbox{in }H_1(\Om_{\vec{\Phi}},\p \Om_{\vec{\Phi}},{\Z}_2)
 \]
 where
 \[
 \Om_{\vec{\Phi}}:=\lf\{(z,t)\in \p B^3\times (-1,+1)\quad;\quad d_{\mathbf F}(\vec{\Phi}(z,t),{\mathcal D})<\ep\rg\}
 \]
 Such a space is called $4-$sweep-outs of $B^3$ by $A$. The set of  $4-$sweep-outs of $B^3$ by $A$ is denoted $\mbox{Sweep}_4(A,B^3)$.\hfill $\Box$
\end{Dfi}

 We have the following lemma.
\begin{Lm}
\label{lm-II.5e}
\[
\mbox{Sweep}_4(A,B^3)\ne\emptyset\quad.
\]
Moreover,  there exists $\vec{\Phi}\in \mbox{Sweep}_4(A,B^3)$ such that
\be
\label{II.4}
\max_{(z,t)\in B^3\times(-1,+1)}\mbox{Area}_{S^3}(\vec{\Phi}(z,t))< 4\,\pi
\ee
where $\mbox{Area}_{S^3}$ is the area taken with respect to the metric obtained by pulling back the standard metric on $S^3$ using the stereographic projection.
\hfill $\Box$
\end{Lm}

\begin{figure}
\begin{center}
\includegraphics[width=18cm,height=24cm]{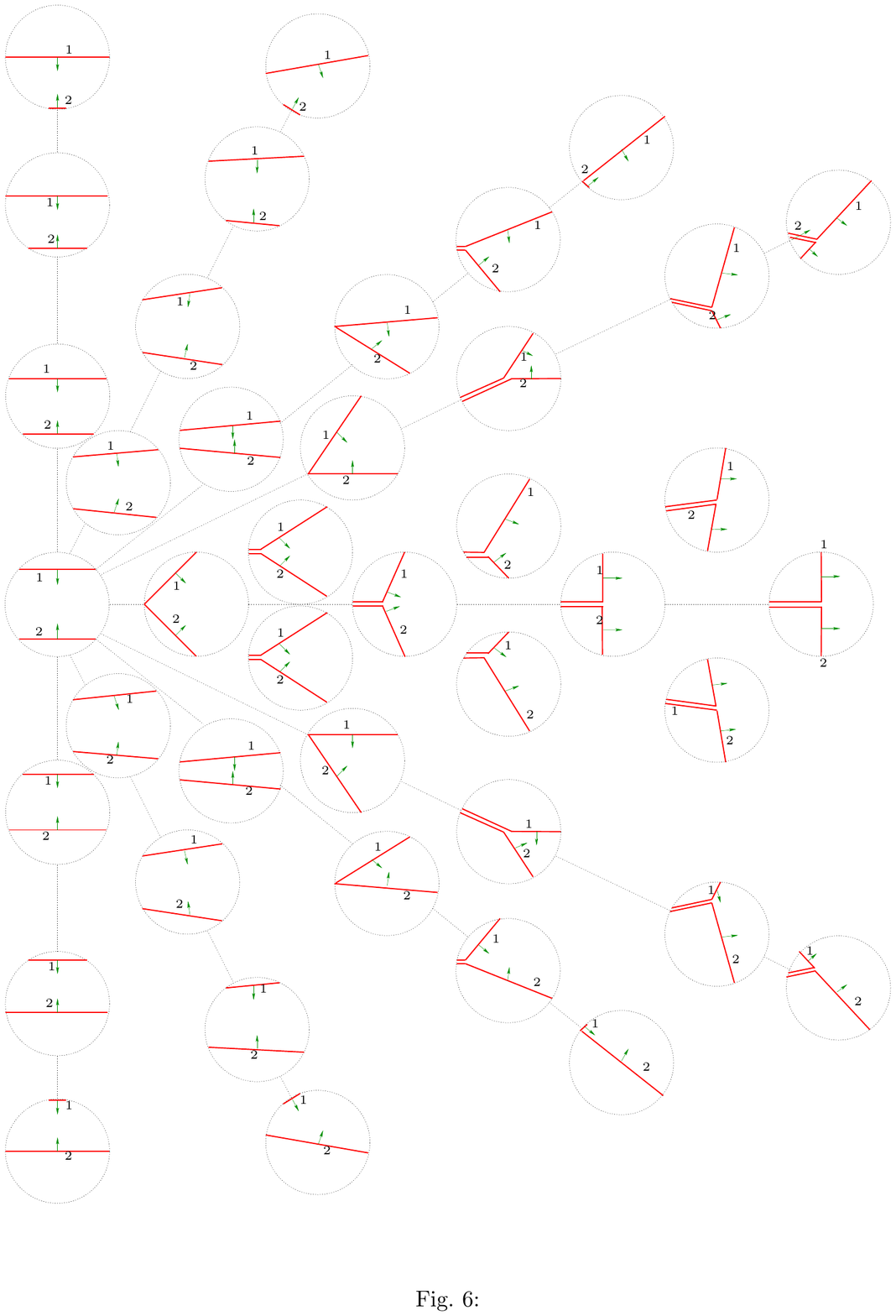}
\end{center}
\end{figure}
\noindent{\bf Proof of lemma~\ref{lm-II.5e}.} Let $\vec{\Xi}$ be the map given by (\ref{xi-1}), (\ref{xi-2}) and (\ref{xi-3}). We shall first construct a modification $\vec{\Xi}_\ast$ of $\vec{\Xi}$
which takes values in $SP_2^\ast(\mbox{Imm}^0(S^1,S^2))$. We do it as follows. We first construct a map $\vec{\Psi}$ from the half disc $D^2_+$ given by 
\[
D^2_+:=\lf\{(x,0,z)\ ; \ x^2+z^2<1\quad x>0\rg\}
\]
into $SP_2^\ast(\mbox{Imm}_T([0,1],D^2))$ where $\mbox{Imm}_T([0,1],D^2)$ is the space of smooth immersions of segments in the discs without boundary inside the disc and cutting
$\p D^2$ transversally. We shall denote by $\vec{\Psi}$ this map.

\medskip

\noindent{\bf Step 1.} {\it Construction of $\vec{\Psi}\in \mbox{Lip}_{\mathbf F}(D^2_+\times (-1,1),SP_2^\ast(\mbox{Imm}_T([0,1],D^2))$.}

\medskip

Let $\Om(r,\theta)$ be the following subdomains of $D^2$

\begin{itemize}
\item[i)] For $\theta\in [-\pi/2,\pi/2]$ and $0<r<1/2$ we set
\be
\label{xi-2a}
\begin{array}{l}
\ds\Om(r,\theta):=\lf\{(x_1,x_2)\in D^2\ ;\ x_1\,\cos\theta_r+x_2\sin\theta_r>-1/2\rg\}\\[3mm]
\ds\quad\quad\quad\ \cap\lf\{  (x_1,x_2)\in D^2\ ;\ x_1\,\cos\theta_r-x_2\sin\theta_r>-1/2    \rg\}
\end{array}
\ee
where $\theta_r:=(1-2r)\pi/2+2\, r\,|\theta|$
\item[ii)] For $\theta\in [-\pi/2,\pi/2]$ and $1/2\le r<1$ we set
\be
\label{xi-3a}
\begin{array}{l}
\ds\Om(r,\theta):=\lf\{(x_1,x_2)\in D^2\ ;\ x_1\,\cos\theta+x_2\sin\theta>-r\rg\}\\[3mm]
\ds\quad\quad\quad\ \cap\lf\{  (x_1,x_2)\in D^2\ ;\ x_1\,\cos\theta-x_2\sin\theta>-1+r    \rg\}
\end{array}
\ee
\end{itemize}
Since $\Om(r,\theta)$ is obtained by intersecting two half planes and one disc, it is convex. The intersection of the boundary of $\Om(r,\theta)$ with the interior of $D^2$ is made either of two disjoint segments, one single segments or the union of  two segments intersecting at their ends inside $D^2$. It is not difficult to verify that the domain of points $(x,z)=(r\,\cos\theta,r\,\sin\theta)$
where $D^2\setminus\Om(r,\theta)$ is connected is diffeomorphic to a disc. We denote this domain by $U$. On $\p U\cap \{(x,z)\ ;\ x^2+z^2<1\ , \ x>0\}$ the two segments given by $\vec{\Xi}(x,0,z)$
intersect at exactly  one point which is on the boundary of $D^2$. We denote  by $J(x,z)$ this point. Observe that $\p U\cap \{(x,z)\ ;\ x^2+z^2<1\ , \ x>0\}$ is a connected strict subset of $\p U\simeq S^1$. For that reason\footnote{We have indeed a trivial fibration over $U\simeq D^2$ and the boundary value one subsegment of $\p U$ can be extended to the whole of $U$ continuously by parallel transport} one can produce a continuous map $S=(P_1,P_2)\in C^0(U,(D^2)^2)$ from $U$ into the space of oriented segments in $D^2$ such that for any $(x,z)\in U$
\[
P_1(x,z)\in D^2\quad\, P_2(x,z)\in \p\Om(r,\theta)\cap D^2\quad\mbox{ s. t. } (P_1(x,z),P_2(x,z))\subset D^2\setminus \Om(r,\theta)
\]
and
\[
P_1(x,z)=P_2(x,z)=J(x,z)\quad\quad\mbox{ on }\p U\cap  \{(x,z)\ ;\ x^2+z^2<1\ , \ x>0\}\quad.
\]
We split the segments $[P_1,P_2]$ into two parallel segments in such a way to construct a ``canal'' $C_\ep(r,\theta)$ joining $\Om(r,\theta)$ to $\p D^2$ inside $D^2\setminus \Om(r,\theta)$ splitting this last set into two components (see figure 6). The boundary of the new two components of $D^2\setminus (\Om(r,\theta)\cup C_\ep(r,\theta))$ realizes the map $\vec{\Psi}$ restricted to the half disc $\{(x,z)\ ;\ x^2+z^2<1\ , \ x>0\}$. The extension on $B^3$ is obtained by the imposed axial symmetry (\ref{xi-1})

\medskip

For $t\in (0,1)$ one modifies continuously $\vec{\Psi}$ by taking $\Om_t(r,\theta)$ to be the following subdomains of $D^2$
\begin{itemize}
\item[i)] For $\theta\in [-\pi/2,\pi/2]$ and $0<r<1/2$ we set
\be
\label{xi-2ab}
\begin{array}{l}
\ds\Om_t(r,\theta):=\lf\{(x_1,x_2)\in D^2\ ;\ x_1\,\cos\theta_r+x_2\sin\theta_r>-1/2-t/2\rg\}\\[3mm]
\ds\quad\quad\quad\ \cap\lf\{  (x_1,x_2)\in D^2\ ;\ x_1\,\cos\theta_r-x_2\sin\theta_r>-1/2-t/2    \rg\}
\end{array}
\ee
where $\theta_r:=(1-2r)\pi/2+2\, r\,|\theta|$
\item[ii)] For $\theta\in [-\pi/2,\pi/2]$ and $1/2\le r<1$ we set
\be
\label{xi-3ab}
\begin{array}{l}
\ds\Om_t(r,\theta):=\lf\{(x_1,x_2)\in D^2\ ;\ x_1\,\cos\theta+x_2\sin\theta>-(1-t)\,r-t\rg\}\\[3mm]
\ds\quad\quad\quad\ \cap\lf\{  (x_1,x_2)\in D^2\ ;\ x_1\,\cos\theta-x_2\sin\theta>-1+r (1-t)   \rg\}
\end{array}
\ee
\end{itemize}
the rest of the construction of $\vec{\Psi}$ for $t\in(0,1)$ is following in a continuous way all the above steps for constructing $\vec{\Psi}$ at $t=0$ replacing continuously $\Om_t(r,\theta)$ by $\Om(r,\theta)$ as well as taking the size of the canal $\ep_t$ tending to zero as $t\rightarrow +1$ (see a slice for $\theta=0$ in the upper half part of figure 7).

\medskip

The construction of $\vec{\Psi}$ for $t<0$ is following a slightly different procedure. We consider the same sets $\Om_t(r,\theta)$ as defined by the formulas (\ref{xi-2ab}) and (\ref{xi-3ab})
and where $t\in(-1,0)$. Starting from the sets $\Om_t(r,\theta)\cup C_{\ep_t}(r,\theta)$ we continuously translate this set such that it is not changed on the vertical axis but
in such a way that 
\[
\lim_{t\rightarrow -1}\max_{r,\theta}\mbox{Area}\lf(\Om_t(r,\theta)\cup C_{\ep_t}(r,\theta)\rg)=0
\]
(see a slice for $\theta=0$ in the lower half part of figure 7).

\begin{figure}
\begin{center}
\includegraphics[width=18cm,height=24cm]{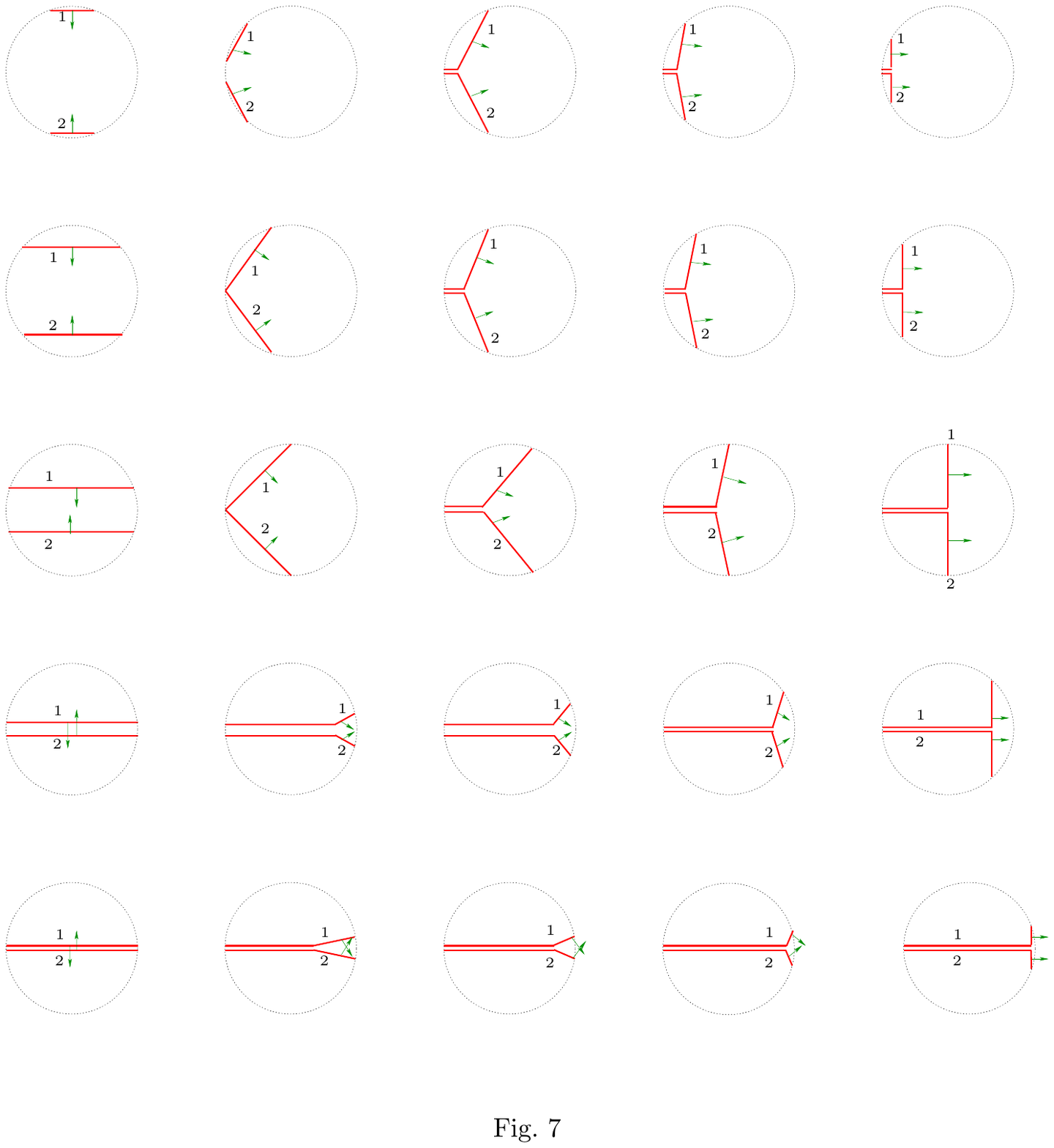}
\end{center}
\end{figure}

\medskip 

\noindent{\bf Step 2.} {\it Construction of $\vec{\Xi}_\ast\in Lip_{\mathbf F}(B^3\times (-1,+1),SP_2^\ast(\mbox{Imm}_T(S^1,S^2))$.}

\medskip

We first construct the map $\vec{\Xi}_\ast$ on the half disc $D^2_+:=\lf\{(x,0,z)\ ; \ x^2+z^2<1\quad x>0\rg\}$. We define $\vec{\Xi}_\ast(x,0,z)$ to be given by the union
of portion of vertical  circles in $S^2=\{(x_1,x_2,x_3)\ ;\ x_1^2+x^2+x^2_3=1\}$ given by the intersection between $S^2$ and union of vertical  portions of planes whose intersection with
$D^2:=\{(x_1,x_2)\ ;\ x_1^2+x_2^2<1\}$ is exactly given by the map $\vec{\Psi}$ previously constructed. In other words, if Let $\pi_3$ denotes the canonical projection from $B^3$ into $D^2$ given by 
\[
\pi_3\ ;\ (x_1,x_2,x_3)\in B^3\quad\longrightarrow \quad (x_1,x_2)\in D^2
\]
We have
\be
\label{fund-ident}
\forall (x,0,z)\in D^2_+\quad\forall\ t\in (-1,+1)\quad \pi_3\circ \vec{\Xi}_\ast(x,0,z,t)=\vec{\Phi}(x,0,z,t)\quad.
\ee

 The complete map $\vec{\Xi}_\ast$ is finally  constructed from it's restriction on the half disc $D^2_+$  by imposing the rotation invariance (\ref{xi-1}).
 
 \medskip
 
 \noindent{\bf Step 3.} {\it Construction of $\vec{\Phi}\in \mbox{Sweep}_4(A,B^3)$ satisfying (\ref{II.4}).} 

\medskip

We first construct the map from the half disc $D^2_+:=\lf\{(x,0,z)\ ; \ x^2+z^2<1\quad x>0\rg\}$. We first bound the non-parallel parts of  $\vec{\Xi}_\ast(x,0,z,t)$ by their
vertical parts. This constructs partly $\vec{\Phi}(x,0,z,t)$. The parallel parts of $\vec{\Xi}_\ast(0,0,z,t)$ are also bounded by the corresponding vertical discs. Now coming to the ``canal parts'' $\pi_3^{-1}(\p C_{\ep_t}(r,\theta))\cap S^2$ are bounded by a tiny catenoid of area $O(\ep_t)$ included in $B^3$ and also bounding the vertical lines given by the boundary of the non parallel parts of  $\vec{\Phi}_\ast(x,0,z,t)$. Finally, we extends continuously the small catenoid junctions to the points $(x,0,z)$ where $C_{\ep_t}(r,\theta)=\emptyset$ by small tubes of area $O(\ep_t)$ connecting the two vertical discs bounded by $\vec{\Xi}_\ast(x,0,z,t)$ and included in $\pi^{-1}_3(\Om_t(r,\theta))$ in such a way that the tubes are axially symmetric for $(0,0,z)$ on the $z$ axis. This defines $\vec{\Phi}(x,0,z)$ for every $(x,0,z)\in D^2_+$. The complete map $\vec{\Phi}$ is finally  constructed from it's restriction on the half disc $D^2_+$  by imposing the rotation invariance (\ref{xi-1}). One easily verify that $\vec{\Phi}$ as it has been constructed is an element of $\mbox{Sweep}_4(A,B^3)$.

Modulo an error bounded by $O(\ep_t)$ the whole area of the resulting immersed annulus $\vec{\Phi}(x,y,z)$ is given for every point $(x,y,z)$ by a pair of subsets of vertical discs included in $B^3$ which are at a distance bounded form below by a positive constant. The maximal area for a flat disc in $B^3$ for the metric induced by the pull-back of the $S^3$ metric by the inverse of the stereographic projection is exactly achieved by discs passing through the center and is equal to $2\pi$. Hence we have established (\ref{II.4}). This concludes the proof of
 lemma~\ref{lm-II.5e}. \hfill $\Box$

\medskip

\noindent{\bf Proof of proposition~\ref{pr-II.4}.} Denote $\ep_1\cdots \ep_4$ the canonical basis of ${\R}^4$. Let $\phi\in [0,2\pi]$ and $\psi\in [0,\pi)$ we denote
\[
\lf\{
\begin{array}{l}
\ds\ep_1^{\phi,\psi}:=\cos\psi(\cos\phi\,\ep_1+\sin\phi\,\ep_2)+\sin\psi\, \ep_4\ \\[3mm]
\ds\ep_2^{\phi,\psi}=-\sin\psi\,\ep_1+\cos\psi\,\ep_2\\[3mm]
\ds \ep_3^{\phi,\psi}=\ep_3\ \\[3mm]
\ds\ep_4^{\phi,\psi}=-\sin\psi(\cos\phi\,\ep_1+\sin\phi\,\ep_2)+\cos\psi\, \ep_4\ \\[3mm]
\end{array}
\rg.
\]
in such a way that $\ep_i^{0,0}=\ep_i$.
Let $\Pi^\psi_\phi$ be the stereographic projection with respect to $\ep_4^{\phi,\psi}$ from $S^3$ onto $\mbox{Span}\lf\{\ep_1^{\phi,\psi},\ep_2^{\phi,\psi},\ep_3^{\phi,\psi}\rg\}$.
We interpret this projection as taking values into ${\R}^3$ by using the isomorphism 
\[
x_1\, \ep_1^{\phi,\psi}+x_2\,\ep_2^{\phi,\psi}+x_3\,\ep_3^{\phi,\psi}\longrightarrow\ (x_1 ,x_2,x_3)
\]
 Consider the plane in ${\R}^3$ passing through the origin and orthogonal
to the vector of coordinates $(\cos\theta,0,\sin\theta)$. It's preimage by $\Pi^\psi_\phi$ is a geodesic $S^2$ in $S^3$ (i.e. an element of ${\mathcal S}$) given by the intersection between $S^3$ and the hyperplane
of ${\R}^4$ passing through the origin and orthogonal to
\be
\label{degre}
\begin{array}{l}
\ds\lf.d\lf(\Pi^\psi_\phi\rg)^{-1}\rg|_{(0,0,0)}\cdot(\cos\theta,0,\sin\theta)=\cos\theta\ \ep_1^{\phi,\psi}+\sin\theta\ \ep_3^{\phi,\psi}\\[3mm]
\ds\quad=\lf(\cos\theta\, \cos\psi\,\cos\phi,\cos\theta\, \cos\psi\,\sin\phi , \sin\theta , \cos\theta\,\sin\psi\rg)
\end{array}
\ee
Hence the map given in spherical coordinates by
\be
\label{geod-sphere}
(\phi,\psi,\theta)\in [0,2\pi)\times [0,\pi)\times [-\pi/2,\pi/2]\ \longrightarrow \ \lf(\Pi^\psi_\phi\rg)^{-1}\lf(\lf\{(\cos\theta,0,\sin\theta),0\rg\}\rg)\in {\mathcal S}\quad\mbox{is a deg. 1 map.}
\ee 

\medskip

Let $\vec{\Phi}(x,t)$ be the map constructed in the proof of
lemma~\ref{lm-II.5e}. For all $(\phi,\psi)$ in $[0,2\pi]$ $r\in [0,1]$ we denote
\[
\vec{\La}(r\,\cos\theta\ \ep_1^{\phi,\psi}+r\,\sin\theta\ \ep_3^{\phi,\psi},t):=(\Pi^\psi_\phi)^{-1}\lf\{\vec{\Phi}(r\, \cos\theta, 0,r\, \sin\theta, t)\rg\}
\]
We have that for every $(r,\theta,\phi,\psi,t)$ the image $\vec{\La}(r\,\cos\theta\ \ep_1^\phi+r\,\sin\theta\ \ep_3^\psi,t)$ is an annulus intersecting
transversally the geodesic $2-$sphere given by $v\in S^3$ s.t. $\ep_4^{\phi\psi}\cdot v=0$. We complete such an immersion by using a reflexion 
with respect to the hyperplane $x\in {\R}^4$ such that $\ep_4^\psi\cdot x=0$ in order to obtain an oriented torus (after possibly  smoothing the resulting
surface along the edge realized by the connecting parts in case $\vec{\La}(r\,\cos\theta\ \ep_1^\phi+r\,\sin\theta\ \ep_3^\psi,t)$ is not intersecting $\ep_4^\psi\cdot x=0$ orthogonally).
We observe that the resulting map, due to (\ref{geod-sphere}) belongs to $\mbox{Sweep}_5(T^2,S^3)$ and because of (\ref{II.4}) it satisfies 
\[
\max_{y\in (-1,1)\times B^4}\mbox{Area}_{S^3}\,(\vec{\La}(y))< 8\pi\quad.
\]
This concludes the proof of  proposition~\ref{pr-II.4}.
\hfill $\Box$

\subsection{Area Minmax in $\mbox{Sweep}_5(T^2,S^3)$.}

\subsubsection{The second level of the $S^3$ hierarchy and a minmax characterization of Clifford Tori.}
\begin{Th}
\label{th-II.6}
We have 
\[
W_2(S^3):=\inf_{\vec{\Phi}\in \mbox{Sweep}_5(T^2,S^3)}\quad\max_{y\in (-1,1)\times B^4}\mbox{Area}_{S^3}\,(\vec{\Phi}(y))=2\pi^2
\]
and is achieved exactly and exclusively by  the isometric image of the Clifford Torus. 

\hfill $\Box$
\end{Th}
\noindent{\bf Proof of theorem~\ref{th-II.6}.}  Let $\vec{\Phi}$ be the element of $\mbox{Sweep}_5(T^2,S^3)$ constructed in the proof of proposition~\ref{pr-II.4} and satisfying (\ref{in-2picarre}). Introduce the family $\mbox{Sweep}_5(S^3)$ to be the subspace of pairs $(Y,\Psi)$ such that $Y\in {\mathcal P}_5$ and $\mbox{Lip}(Y,\mbox{Imm}^1(S^3))$ and satisfying i)... v) in the definition of {\it Minmax hierarchy}. We clearly have $\vec{\Phi}\in \mbox{Sweep}_5(S^3)$ and hence $W_2(S^3)< 8\pi$. Using theorem~\ref{pr-II.2} and theorem~\ref{th-Pitts}
we deduce that 
\[
W_1(S^3)<W_2(S^3)<8\pi
\]
and $W_2(S^3)$ is achieved by a smooth minimal surface of genus less than one and index less than 5. It cannot be neither a geodesic sphere which has area $4\pi$ or 
several geodesic sphere (or a multiple covering of a geodesic sphere) that would have more that $8\pi$ area ($8\pi>2\pi^2$). Hence, using Urbano's result \cite{Urb}
we deduce that the surface is a {\it Clifford torus} and $W_2(S^3)=2\pi^2$. This concludes the proof of theorem~\ref{th-II.6}.\hfill $\Box$

\medskip

\begin{Con}
\label{crit-cat}
In a similar manner we expect the critical catenoids to be the unique solutions of the following minmax problem
\[
W_2(B^3)=\inf_{\vec{\Phi}\in \mbox{Sweep}_4(A,S^3)}\quad\max_{y\in (-1,1)\times B^3}\mbox{Area}_{\,{\R}^3}\,(\vec{\Phi}(y))
\]
where $\mbox{Sweep}_4(A,S^3)$ is the {\it admissible family } defined in definition~\ref{df-II.5e} which we proved to be non empty in lemma~\ref{lm-II.5e}.
Hence we should have $W_2(B^3)=2\,\pi\,\La^{-2}$ where $\La\,\tanh\La=1$.\hfill $\Box$
\end{Con}
\begin{Rm}
\label{rm-crit-cat}
In  recent works (see \cite{Dev} and \cite{Tra}) the critical catenoids have been characterized as being the unique free boundary surfaces, different from the flat discs, of index less or equal than 4. This result
should play the role of Urbano's result, in the proof of theorem~\ref{th-II.6}, in the way to prove the above conjecture. What is missing at this stage is the free boundary version of
our result theorem~\ref{th-A.3}.\hfill $\Box$
\end{Rm}
\subsubsection{An alternative  proof of theorem~\ref{th-II.6} using Bryant's Computations of the Conformal Volumes of the CMC Clifford tori.}  We are giving now an argument that leads to a new proof of theorem~\ref{th-II.6} modulo some weakening
of the notion of $\mbox{Sweep}_5(T^2,S^3)$ and allow some singularities of zero measure.

\medskip

Using in particular remark~\ref{rm-I.3p} we can define $\mbox{Sweep}_5(T^2,S^3)$ to be the elements $$\vec{\Phi}\in \mbox{Lip}_{\mathbf F}(B^4, \mbox{Sweep}_1(T^2,S^3))$$ (where the
definition of $\mbox{Sweep}_1(T^2,S^3)$ is the one we gave in the previous sub-sections) satisfying

\begin{itemize}
\item{i)} There exists $K$ compact in $\p B^4$ such that
\[
\forall\ z\in \p B^4\setminus K\quad\quad \lim_{\rho\rightarrow 1}\vec{\Phi}(\rho z)\rightarrow \vec{\Psi}(z)\quad
\]
where $\forall\ z\in \p B^4\setminus K$
\[
\vec{\Psi}\in \mbox{Sweep}_1(S^2,S^3)\quad\mbox{ and }\quad\sup_{t\in (-1,+1)}\mbox{Area}(\vec{\Psi}(z,t))\le 4\pi+\delta_{\mathcal S}/2
\]
\item{ii)} Let $\La_{\vec{\Phi}}:=\{(z,t)\in B^4\times (-1,+1)\quad;\quad d_{\mathbf F}(\vec{\Phi}(z,t),{\mathcal S})<\ep_{{\mathcal S}}\}$ then
\[
\lim_{t\rightarrow 1}{\mathcal H}^3\lf(\pi_{\mathcal S}\lf(\La_{\vec{\Phi}}\cap \lf(B^4\times (t,+1)\cup B^4\times (-1,-t)\rg) \rg)\rg)=0
\]
\end{itemize}

\medskip

 It is not difficult to see that our main result theorem~\ref{pr-II.2} extends to this weaker notion of hierarchy including singularities. We are now making use of the computation of the conformal volume of the CMC Clifford tori given in \cite{Bry}. For any $b\in {\R}^+_\ast$ we consider the family of Clifford Tori $\mbox{Cl}_b$ in $S^3$ given by
\be
\label{II.4-a}
\mbox{Cl}_b\quad :\ (\theta,\phi)\in {\R}/2\pi{\Z}\times {\R}/2\pi{\Z}\longrightarrow\ \frac{1}{\sqrt{1+b^2}}(e^{i\theta},b\, e^{i\phi})
\ee
(for $b=1$ this projection is nothing but the Willmore torus).  For $t\in (-1,1)$ let $b_t:=(1+t)/(1-t)$. Observe that
\[
t\in(-1,+1)\ \longrightarrow \ \mbox{Cl}_{b_t} \ \in \mbox{Sweep}_1(T^2,S^3)
\]
Now for $\vec{a}\in B^4$ we consider the conformal transformation from $S^3$ into $S^3$ given by
\be
\label{moebius}
\varphi_{\vec{a}}(\vec{z}):=(1-|\vec{a}|^2)\frac{\vec{z}-\vec{a}}{|\vec{z}-\vec{a}|^2}-\vec{a}
\ee
and we introduce the map
\[
\vec{\Phi}\ :\ (\vec{a},t)\in B^4\times(-1,+1)\ \longrightarrow\ \varphi_{\vec{a}}\circ \mbox{Cl}_{b_t} \quad.
\]
Observe that
\[
\vec{\Phi}\in \mbox{Lip}_{\mathbf F}(B^4,\mbox{Sweep}_1(T^2,S^3))\quad\mbox{ and }\forall\, t\in (-1,+1)\quad\limsup_{|\vec{a}|\rightarrow +1}\mbox{Area}(\Phi(\vec{a},t))\le 4\pi=W_1(S^3)
\]
Denote for $s<1$
\[
\Om_{\vec{\Phi}}(s,t):=\lf\{\vec{a}\in \p B^4_s(0)\quad; \quad d_{\mathbf F}(\vec{\Phi}(\vec{a},t),\mathcal{S})<\ep_{{\mathcal S}}\rg\}
\]
It is not difficult to prove that
\[
\Om_{\vec{\Phi}}(t):=\lim_{s\rightarrow 1}\Om_{\vec{\Phi}}(s,t)=\lf\{(\frac{1}{\sqrt{1+b_t^2}}(i\,e^{i\theta},b_t\,i\, e^{i\phi})\ ;\ (\theta,\phi)\in {\R}/2\pi{\Z}\times {\R}/2\pi{\Z}\rg\}\quad.
\]
Introducing now $\Om_{\vec{\Phi}}:=\Om_{\vec{\Phi}}((-1,+1))$ since the map
\[
(t,\theta,\phi)\in (-1,+1)\times{\R}/2\pi{\Z}\times {\R}/2\pi{\Z}\ \longrightarrow \frac{1}{\sqrt{1+b_t^2}}(i\,e^{i\theta},b_t\,i\, e^{i\phi})\in S^3\setminus \Gamma_\pm
\]
where $\Gamma_{+}:={\C}\oplus 0\cap S^3$ and $\Gamma_{-}:=0\oplus {\C}\cap S^3$ is an homeomorphism and one verifies that for  every $z\in \p B^4\setminus K$
where $K:=\Gamma_+\cup\Gamma_-$ 

$$[\Om_{\vec{\Phi}}\cap \{z\}\times (-1,+1)]\in H_1(\Om_{\vec{\Phi}},\p\,\Om_{\vec{\Phi}},{\Z})\quad\mbox{ is Poincar\'e dual to } \quad\vec{\Phi}^\ast\,\pi_{\mathcal S}^\ast\,\om_{\mathcal S}\quad.$$
From \cite{Bry} we have
\be
\label{II.5}
\lf\{
\begin{array}{l}
\ds\forall \ t\ge\frac{\sqrt{2}-1}{\sqrt{2}+1}\quad\quad\max_{\vec{a}\in B^4}\mbox{ Area }\lf(\vec{\Psi}(\vec{a},t)\rg)=\frac{8}{3}\ \sqrt{\frac{2}{3}}\ \pi^2\,\frac{\sqrt{1+t^2}}{1+t}\\[5mm]
\ds\forall \ t\le\frac{-\sqrt{2}+1}{\sqrt{2}+1}\quad\quad\max_{\vec{a}\in B^4}\mbox{ Area }\lf(\vec{\Psi}(\vec{a},t)\rg)=\frac{8}{3}\ \sqrt{\frac{2}{3}}\ \pi^2\,\frac{\sqrt{1+t^2}}{1-t}\\[5mm]
\ds\forall\ t\in\lf[\frac{-\sqrt{2}+1}{\sqrt{2}+1},\frac{\sqrt{2}-1}{\sqrt{2}+1}\rg]\quad\quad\max_{\vec{a}\in B^4}\mbox{ Area }\lf(\vec{\Psi}(\vec{a},t)\rg)=2\ \pi^2\ \frac{1-t^2}{1+t^2}
\end{array}
\rg.
\ee
Hence from (\ref{II.5}) we deduce
\[
\max_{y\in B^4\times (-1,+1)}\ \mbox{Area}(\vec{\Phi}(y))\le \frac{8\,\pi^2}{3\,\sqrt{2}}<2\,(4\pi)=8\pi\quad.
\]
The rest of the proof follows as before after having verified that theorem~\ref{pr-II.2} holds for this weaker notion of $\mbox{Sweep}_5(T^2,S^3)$ and observed that the conditions i) and ii) are fulfilled for this 5 dimensional family of conformal transformations of {\it CMC Clifford tori}.\hfill $\Box$

\medskip

\begin{Rm}
\label{rm-willmore} One could wonder if the min-max characterization of the Clifford Tori given by theorem~\ref{th-II.6} is related to the one given by Marques and Neves in their proof
of the Willmore conjecture. For any minimal immersion of non zero genus $\vec{\Phi}\, :\, \Sigma^g\, \rightarrow \, S^3$, the associated canonical family as it has been de-singularized
in \cite{MN} along $\vec{\Phi}(\Sigma^g)$ can be seen as an element of $\mbox{Sweep}_5(S^3)$ taking values into the 2-rectifiable cycles ${\mathcal Z}_2(S^3)$ instead of $\mbox{Imm}(\Sigma^g,S^3)$.
Moreover it is continuous with respect to the Flat distance ${\mathcal F}$ and not necessary with respect to the ${\mathbf F}$ distance (when the focal set of  $\vec{\Phi}(\Sigma^g)$
has non zero 2 dimension). In order to obtain the Willmore conjecture as a corollary of theorem~\ref{th-II.6} one would have to ``smooth-out'' the de-singularized canonical family of Marques and Neves of any minimal immersion
of area less than $8\pi$ in order to make it as an element of $\mbox{Sweep}_5(\Sigma^g,S^3)$ as we defined it in the previous section.\hfill $\Box$
\end{Rm}

\subsection{Towards Higher Order Sweep-Outs and Higher Levels in the $S^3$ Hierarchy.}
This subsection is more of {\bf speculative nature}.
Denote by ${\mathcal T}$ the space of cycles given by the courant of integration of $T$, the image in $S^3$ by an arbitrary isometry of the oriented Clifford Torus for which we have singled out an orientation of the closed geodesics of minimal lengths compatible with the orientation of the torus. To every such a torus we can assign the oriented core circle which is the unique large circle equidistant to the whole torus contained inside the domain $S^3_{+,T}$ of $S^3$ having $T$ as a boundary  and homologous to one of these oriented geodesics in $S^3_{+,T}$. Vice versa,  any  oriented large circle in $S^3$ being given, there exists a unique embedded $T\in {\mathcal T}$  equipped with an orientation of it's closed geodesics of minimal length whose core circle is given by this circle. Hence
${\mathcal T}$ identifies with the space of  oriented large circles in $S^3$ which itself identifies to the Grassmann manifold $G_{2}({\R}^4)$ of
oriented 2-planes in ${\R}^4$. Recall that
\[
\tilde{G}_{2}({\R}^4)\simeq  S^2\times S^2
\]
where for any $(x,y)\in S^2\times S^2$ the first entry $x$ corresponds to a choice of a complex structure $J_x$ on ${\R}^4$ compatible with the flat metric and the associated oriented 2 plane is given by it's intersection with $S^3$ and corresponds to $\pi_x^{-1}(\{y\})$ where $\pi_x$ is the Hopf fibration compatible with $J_x$ in the sense that $\pi_x$ is holomorphic from $({\R}^4\setminus\{0\},J_x)$ into ${\mathbb C}{\mathbb P}^1$. Since $H^2(S^2\times S^2,{\Z}_2)={\Z}_2\oplus{\Z}_2$ there are two choices for $\om_{\mathcal T}$ and the tree hierarchy is splitting
into two branches. 

\medskip

Hence, denoting by $\Sigma^2$ the abstract closed oriented surface of genus 2, one aims to construct an element in $\mbox{Sweep}_8(\Sigma^2,S^3)$ which is going to be a map in $$\mbox{Lip}(B^3,\mbox{Sweep}_5(\Sigma^2,S^3))$$ whose restriction
to $\p B^3$ is going to be $\ep-$close to a map from $S^2$ into $\mbox{Sweep}_5(T^2,S^3)$ with maximal area $\ep$ close to $2\pi^2$ and such that the topological constraint 
\[
\forall\ z\in \p B^3\quad\quad [\Om_{\vec{\Phi}}\cap \{z\}\times (-1,+1)]\in H_5(\Om_{\vec{\Phi}},\p\,\Om_{\vec{\Phi}},{\Z}_2)\quad\mbox{ is Poincar\'e dual to } \quad\vec{\Phi}^\ast\,\pi_{\mathcal T}^\ast\,\om_{\mathcal T}
\]
is satisfied. The restriction of $\vec{\Phi}$ to $B^3\times\{0\}$ should be given by the connected sum of the two Hopf tori given by the pre-image by  the Hopf fibration of the map $\vec{\Xi}_\ast(0,\cdot)$ showed in figure 6 which has been constructed in the step 2 of the proof of lemma~\ref{lm-II.5e}.

Instead of looking for an almost explicit example as we did in the lower levels of the $S^3$ hierarchy, the existence of  $\mbox{Sweep}_8(\Sigma^2,S^3)$ could be investigated more ``abstractly'' by constructing  first a smooth inclusion map
\[
\iota\ : \ \mbox{Sweep}_5(T^2,S^3)\quad\longrightarrow\quad\mbox{Sweep}_5(\Sigma^2,S^3)
\] 
and by studying the boundary operator $\p_\ast$ in the following associated exact sequence
\[
\begin{CD}
... @>j_\ast>> \pi_3(\mbox{Sweep}_5(\Sigma^2,S^3),\mbox{Sweep}_5(T^2,S^3))@>\p_\ast>> \pi_2(\mbox{Sweep}_5(T^2,S^3))@>\iota_\ast>>...\end{CD}
\]
Then the next step in the hierarchy should produce index 9 surfaces\footnote{Observe that one of the Lawson surfaces, $\xi_{1,2}$, having genus 2 is conjectured by A.Neves in \cite{Nev} to be of index 9.}, still of genus 2,  obtained by everting the sweep-outs $\mbox{Sweep}_8(\Sigma^2,S^3)$ from one orientation to the other. A more systematic study of the successive elements in the $S^3$ hierarchy should be the subject of a forthcoming work.

 \renewcommand{\theequation}{A.\arabic{equation}}
\renewcommand{\theTh}{A.\arabic{Th}}
\renewcommand{\theProp}{A.\arabic{Prop}}
\renewcommand{\theLma}{A.\arabic{Lma}}
\renewcommand{\theCo}{A.\arabic{Co}}
\renewcommand{\theRm}{A.\arabic{Rm}}
\renewcommand{\theequation}{A.\arabic{equation}}
\setcounter{equation}{0} 
\reset
\appendix
\section{Appendix}

We are using the following proposition which is fairly standard in {\it Palais deformation theory} but for which nevertheless we give a proof below.
\begin{Prop}
\label{pr-A.1}
Let $\{\mbox{Sweep}_{N_l}\}_{l\le k}$ be a minmax hierarchy for a $C^2$ Lagrangian $E$ on a Banach manifold ${\mathfrak M}$. Assume moreover that $E$ satisfies the {\it Palais Smale condition}. Then for any $\ep>0$ there exists
$\delta>0$ such that
\[
\forall\, (Y,{\Phi})\in \mbox{Sweep}_{N_k} \quad \max_{y\in Y} E({\Phi}(y))<W_k+\delta
\]
then there exists $y_k^\ep\in Y$ and a critical point ${\Phi}_k$ of $E$ in ${\mathfrak M}$ such that
\[
E({\Phi}_k)=W_k\quad\mbox{ and }\quad d_{\mathbf P}(\vec{\Phi}(y_k^\ep),{\Phi}_k)<\ep\quad.
\]
where $d_{\mathbf P}$ is the usual {\it Palais distance} issued from the {\it Finsler structure} induced by the {\it Banach Manifold structure} on ${\mathfrak M}$\hfill $\Box$
\end{Prop}
\noindent{\bf Proof of proposition~\ref{pr-A.1}.}  Let $\delta>0$ such that
\[
W_{k-1}+\delta<W_{k}\quad.
\]
Consider the pseudo-gradient $X_k$  defined in ${\mathfrak M}$ and locally lipschitz in ${\mathfrak M}\setminus {\mathfrak M}^\ast$ and multiply  by a cut-off function supported in $[W_{k}-\delta,W_{k}+\delta]$ and equal to one on
$[W_{k}-\delta/2,W_{k}+\delta/2]$ in such a way that we have 
\[
\forall \,{\Phi}\in {\mathfrak M} \quad \|X_k({\Phi})\|_{{\Phi}}\le \|DE({\Phi})\|_{{\Phi}}
\]
and
\[
E({\Phi})\in [W_{k}-\delta/2,W_{k}+\delta/2]\quad\Rightarrow\quad \lf<X_k({\Phi}),DE({\Phi})\rg>_{T_{{\Phi}}{\mathfrak M},T^\ast_{{\Phi}}{\mathfrak M}}>  \|DE({\Phi})\|_{{\Phi}}^2
\]
Since the pseudo-gradient is supported in the level sets larger than $W_{k-1}$ it's flow $\phi_t$ generates a family of homeomorphisms
preserving $\mbox{Sweep}_{N_k}$ due to the conditions iii) and v)  in the definition of a hierarchy. For any $\eta\in (0,\delta/2)$ and any $(Y,{\Phi})\in \mbox{Sweep}_{N_k}$ such that
\[
\max_{y\in Y}E({\Phi}(y))\le W_k+\delta/2
\]
we consider the images $(Y,\phi_t({\Phi}))\in \mbox{Sweep}_{N_k}$. Denote by $d_{\mathbf P}$ the {\it Palais distance} associated
to the Finsler structure $\|\cdot\|$ and for which $(\mathfrak M,d_{\mathbf P})$ is complete (see \cite{Riv-columb}). Following lecture 2 of \cite{Riv-columb} we have for all ${y\in Y}$ and  any $t_1<t_2<t_{max}^{y}$
\be
\label{A.I.1}
 d_{\mathbf P}\lf(\phi_{t_1}({\Phi}(y)),\phi_{t_2}({\Phi}(y))\rg)\le \ 2\ \sqrt{t_2-t_1}\ \lf[E(\phi_{t_1}(\vec{\Phi}(y)))-E(\phi_{t_2}(\vec{\Phi}(y)))\rg]^{1/2}
\ee
where $t_{max}^{y}$ is the maximal existence time of $\phi_t({\Phi}(y))$. For a given $y\in Y$, assuming $t_{max}^y<+\infty$, because of the
previous inequality $\phi_t({\Phi}(y))$ realizes a Cauchy sequence for $d_{\mathbf P}$. Since ${\mathfrak M}$ is complete for $d_{\mathbf P}$,
the only possibility is that $\lim_{t\rightarrow t_{max}^y} \phi_t({\Phi}(y))\in {\mathfrak M}^\ast$. 
Assume first that 
\[
\forall\, y\in Y\quad E({\Phi}(y))\ge W_{k}-\delta \quad\Rightarrow \quad \|DE({\Phi}(y))\|_{{\Phi}(y)}\ge\delta^{1/4}
\]
Let $T>0$ and $\al>0$ such that 
\[
\forall \, t< T\quad \forall\, y\in Y\quad E(\phi_t({\Phi}(y)))\ge W_{k}-\delta \quad\Rightarrow \quad \|DE(\phi_t({\Phi}(y)))\|_{\phi_t({\Phi}(y))}\ge\delta^{1/4}
\]
We have for any ${\Phi}\in E^{-1}([W_{k}-\delta/2,W_{k}+\delta/2])$
\[
-\,\lf.\frac{d E(\phi_t({\Phi}))}{dt}\rg|_{t=0}\ge \|DE({\Phi})\|_{{\Phi}}^2
\]
This implies
\[
\max_{y\in Y}E(\phi_T({\Phi}(y)))\le W_k+\delta/2-\sqrt{\delta}\ T
\]
Denote $T_{\max}$ the first time such that either
\[
\exists \,y\in Y\quad\mbox{ s.t. }\quad t_{max}^y=T_{\max}
\]
which implies $DE(\phi_{T_{max}}(\vec{\Phi}(y))=0$ and $E(\phi_{T_{max}}({\Phi}(y)))\ge W_{k}-\delta$ or
\[
\exists \ y\in Y\quad E(\phi_t({\Phi}(y)))\ge W_{k}-\delta\quad \mbox{ and }\quad \|DE(\phi_{T_{max}}({\Phi}(y)))\|_{\phi_{T_{max}}({\Phi}(y))}={\delta}^{1/4}
\]
We clearly have $\sqrt{\delta}\, T_{max}\le\delta/2$. This implies using (\ref{A.I.1})
\[
\max_{y\in Y}d_{\mathbf P}({\Phi}(y),\phi_{T_{max}}({\Phi}(y)))\le {\delta^{1/4}}\ \sqrt{W_{k}/2}
\]
Collecting the various cases and summarizing we have obtained the following :
\be
\label{A.I.2}
\begin{array}{c}
\ds\forall\ \delta>0\quad\mbox{ and }\forall\ (Y,{\Phi})\in  \mbox{Sweep}_{N_k} \quad \mbox{ if }\quad\max_{y\in Y} E({\Phi}(y))<W_k+\delta\\[3mm]
\ds\quad\mbox{ then }\exists\ y\in Y\quad\mbox{ and }{\Phi}_\delta\in {\mathfrak M}\quad\mbox{ s.t. }
 |E({\Phi}_\delta)-W_k|\le\delta\\[3mm]
 \quad\mbox{and }\lf\|DE({\Phi}_\delta)\rg\|\le \delta^{1/4}
 \mbox{ moreover }\quad d_{\mathbf P}({\Phi}_\delta,{\Phi}(y))\le {\delta^{1/4}}\ \sqrt{W_{k}/2}
\end{array}
\ee
We claim now that 
\be
\label{A.I.3}
\begin{array}{l}
\forall\ \ep>0\ \ \exists\,\delta>0\quad\mbox{ s.t. }\quad\forall \ {\Phi}_\delta\in {\mathfrak M}\quad\mbox{ s.t. }
 |E({\Phi}_\delta)-W_k|\le\delta\\[3mm]
 \quad\mbox{and }\lf\|DE({\Phi}_\delta)\rg\|\le \delta^{1/4}\quad\mbox{ then }\quad\exists\ {\Phi}_0\in{\mathfrak M} \\[3mm]
 \quad\mbox{ s. t. }\quad E({\Phi}_0)=W_k\quad,\quad DE({\Phi}_0)=0\quad\mbox{ and }\quad d_{\mathbf P}({\Phi}_\delta,{\Phi}_0)\le \ep
\end{array}
\ee
This fact is a direct consequence of the {\it Palais Smale Condition} .  Indeed if (\ref{A.I.3}) would be false there
would exist $\ep_0>0$, a sequence $\delta_i\rightarrow 0$ and a sequence ${\Phi}_{\delta_i}$ such that
\[
 |E(\vec{\Phi}_\delta)-W_k|\le\delta_i\quad\mbox{ and }\quad\lf\|DE({\Phi}_{\delta_i})\rg\|\le \delta_i^{1/4}
\]
but ${\Phi}_{\delta_i}$ would stay at a {\it Palais distance} larger than $\ep_0$ to any critical point to $E$ at the level set $W_k$ which contradicts {\it (P.S.)}.

\medskip

Combining (\ref{A.I.2}) and (\ref{A.I.3}) we obtain
\[
\begin{array}{c}
\ds\forall\ \ep>0\ \ \exists\,\delta>0\quad\mbox{ s.t. }\forall\ (Y,{\Phi})\in  \mbox{Sweep}_{N_k} \quad \mbox{ if }\quad\max_{y\in Y} E({\Phi}(y))<W_k+\delta\\[3mm]
\mbox{ then }\quad\exists\ \vec{\Phi}_0\in {\mathfrak M}
 \quad\mbox{ s. t. }\quad E({\Phi}_0)=W_k\quad,\quad DE({\Phi}_0)=0\quad\\[3mm]
 \mbox{ and }\quad\exists\ y\in Y\mbox{ s. t. }\quad d_{\mathbf P}({\Phi}(y),{\Phi}_0)\le \ep
\end{array}
\]
 This concludes the proof of proposition~\ref{pr-A.1}.\hfill $\Box$

\medskip
The previous proposition extends to the non {\it Palais-Smale case} in the framework of the {\it viscosity method } to the following theorem.
\begin{Th}
\label{th-A.2}
Let $\beta(0)>0$, let $g\in {\N}$, $M^m$ be a closed sub-manifold of an euclidian space ${\R}^Q$  and ${\mathcal A}$ be a subset of the power set $\wp(\mbox{Imm}(\Sigma^g,M^m))$ of compact subsets of $\mbox{Imm}(\Sigma^g,M^m)$ (equipped with the $W^{2,4}-$topology) which is invariant under the space of homeomorphisms
of $\mbox{Imm}(\Sigma^g,M^m)$ which coincides with the identity on $\mbox{Area}^{-1}([0,\beta({0})-\eta])$ for some $\eta>0$ and assume
\[
\beta(0)=\inf_{A\in{\mathcal A}}\ \sup_{\vec{\Phi}\in A}\mbox{Area}\,(\vec{\Phi})
\]
Then, for any $\ep>0$ there exists  $\tau>0$ and $\delta_\tau>0$ such that, for any $A\in{\mathcal A}$ satisfying
\[
\sup_{\vec{\Phi}\in A}A^\tau(\vec{\Phi})<\beta(\tau)+\delta_\tau\quad,
\]
then 
\[
d_{\mathbf F}(A,{\mathcal C}^0)<\ep\quad,
\]
where ${\mathcal C}^0$ is the non empty space of smooth minimal immersions (possibly branched) of area equal to $\beta(0)$.
\hfill $\Box$
\end{Th}
Before proving theorem~\ref{th-A.2} we shall be proving the following intermediate lemma.
\begin{Lma}
\label{lm-A.3}
Under the same assumptions as the ones of theorem~\ref{th-A.2} introduce for $\sigma>0$
\[
\beta(\sigma)=\inf_{A\in{\mathcal A}}\ \sup_{\vec{\Phi}\in A} A^\sigma(\vec{\Phi})\quad,
\]
where
\[
A^\sigma(\vec{\Phi}):=\mbox{Area}\,(\vec{\Phi})+\sigma^2\, \int_{\Sigma^g}\lf[1+|{\vec{\mathbb I}}_{\vec{\Phi}}|^2\rg]^2\ dvol_{g_{\vec{\Phi}}}\quad,
\]
and  ${\vec{\mathbb I}}_{\vec{\Phi}}$ denotes the second fundamental form of the immersion $\vec{\Phi}$ into $M^m$.

Let $\sigma<\tau$ an assume $\beta(\sigma)<\beta(\tau)$
\[
\forall\, A\in{\mathcal A}\quad\quad\quad\mbox{ s. t. }\quad\max_{\vec{\Phi}\in A}A^\tau(\vec{\Phi})\le \beta(\tau)+(\beta(\tau)-\beta(\sigma))
\]
then
\[
d_{\mathbf P}(A,\tilde{\mathcal C}_{\tau,\sigma})<C\, (\beta(\tau)-\beta(\sigma))^{1/6}\quad,
\]
where
\[
\tilde{\mathcal C}_{\tau,\sigma}:=\lf\{
\begin{array}{c}
\ds\vec{\Phi}\in \mbox{Imm}_0(\Sigma^g,M^m)\quad;\quad A^\tau(\vec{\Phi})\in[2\,\beta(\sigma)-\beta(\tau),2\,\beta(\tau)-\beta(\sigma)]\quad,\\[3mm]
\ds\quad \|DA^\tau(\vec{\Phi})\|_{\vec{\Phi}}<C\, (\beta(\tau)-\beta(\sigma))^{1/3} \quad\\[3mm]
\ds\mbox{ and }\quad (\tau^2-\sigma^2)\,\int_{\Sigma^g}\lf[1+|{\vec{\mathbb I}}_{\vec{\Phi}}|^2\rg]^2\ dvol_{g_{\vec{\Phi}}}\le \, 3\, (\beta(\tau)-\beta(\sigma))
\end{array}
\rg\}\ne\emptyset\quad.
\]
where $C>0$ is a universal constant. 
\hfill $\Box$
\end{Lma}
\noindent{\bf Proof of lemma~\ref{lm-A.3}.} 
Denote
\[
{\mathcal A}_{\tau,\sigma}:=\lf\{A\in {\mathcal A}\quad;\quad\max_{\vec{\Phi}\in A}A^\tau(\vec{\Phi})\le\beta(\tau)+ (\beta(\tau)-\beta(\sigma))\rg\}
\]
Observe that
\[
\begin{array}{l}
\ds\forall \, A\in  {\mathcal A}_{\tau,\sigma}\quad\quad \forall\, \vec{\Phi}\in A\\[3mm]
\ds\quad \beta(\sigma)-  (\beta(\tau)-\beta(\sigma))\le A^\sigma(\vec{\Phi})\quad\Longrightarrow \quad (\tau^2-\sigma^2)\,\int_{\Sigma^g}\lf[1+|{\vec{\mathbb I}}_{\vec{\Phi}}|^2\rg]^2\ dvol_{g_{\vec{\Phi}}}\le \, 3\,(\beta(\tau)-\beta(\sigma))
\end{array}
\]
Let $X^\tau$ be a Lipschitz pseudo-gradient on $\mbox{Imm}_0(\Sigma^g,M^m)$ equipped with the $W^{2,4}-$Finsler structure described in \cite{Riv-columb}. We have on $(\mbox{Imm}_0(\Sigma^g,M^m))^\ast:=\mbox{Imm}_0(\Sigma^g,M^m)\setminus \{\vec{\Phi}\ ;\ DA^\tau(\vec{\Phi})=0\}$
\be
\label{A.I.5}
\lf<X^\tau(\vec{\Phi}),DA^\tau(\vec{\Phi})\rg>\ge \|DA^\tau(\vec{\Phi})\|_{\vec{\Phi}}^2\quad\mbox{ and }\quad \|X^\tau(\vec{\Phi})\|_{\vec{\Phi}}\le 2\ \|DA^\tau(\vec{\Phi})\|_{\vec{\Phi}}
\ee
Introduce 
$\chi\in C^\infty({\R})$ such that $\chi\equiv 1$ on $[1/2,+\infty)$ and $\chi\equiv 0$ on ${\R}_-$ and multiply the pseudo-gradient by the following cut-off
\[
\tilde{X}^\tau_\sigma(\vec{\Phi}):=\chi\lf( \frac{A^\sigma(\vec{\Phi})-2\,\beta(\sigma)+\beta(\tau)}{2\,\beta(\sigma)-\beta(\tau)}\rg)\ X^\tau(\vec{\Phi})
\]
in order for the pseudo-gradient to deform the elements from ${\mathcal A}_\tau$ in the zone exclusively where $A^\sigma(\vec{\Phi})\ge 2\,\beta(\sigma)-\, \beta(\tau)$ and hence in the zone where $\int_{\Sigma^g}\lf[1+|{\vec{\mathbb I}}_{\vec{\Phi}}|^2\rg]^2\ dvol_{g_{\vec{\Phi}}}\le  \, 3\,(\beta(\tau)-\beta(\sigma)) $.

\medskip

 Introduce the unique flow issued from the Lipschitz vector-field on $(\mbox{Imm}_0(\Sigma^g,M^m))^\ast$ given by
\[
\lf\{
\begin{array}{l}
\ds\frac{d \varphi_t(\vec{\Phi})}{dt}=-\,\tilde{X}^\tau_\sigma(\varphi_t(\vec{\Phi}))\\[3mm]
\ds\varphi_0(\vec{\Phi})=\vec{\Phi}
\end{array}
\rg.
\]
Let $A\in {\mathcal A}_\tau$ arbitrary. It is proved in \cite{Riv-columb} that the flow starting from $A$ exists for all time as long as $\|DA^\tau( \varphi_t(A))\|>0$. We have $\forall\,\Phi\in A$
\be
\label{A.I.6}
\lf\{
\begin{array}{l}
\ds A^\sigma(\varphi_t(\vec{\Phi}))\le A^\tau(\varphi_t(\vec{\Phi}))\le  A^\tau(\vec{\Phi})\le2\,\beta(\tau)-\beta(\sigma)\quad,\\[3mm]
\ds A^\sigma(\vec{\Phi})\le 2\,\beta(\sigma)-\beta(\tau)\quad\Longrightarrow\quad\forall\ t>0\quad \varphi_t(\vec{\Phi})=\vec{\Phi}
\end{array}
\rg.
\ee
It is proved in \cite{Riv-minmax} that, for any $\vec{w}\in W^{2,4}\cap W^{1,\infty}(\Gamma\lf(\vec{\Phi}^\ast T M^m\rg))$
\be
\label{A.I.7}
\lf|DA^\sigma(\vec{\Phi})\cdot\vec{w}-DA^\tau(\vec{\Phi})\cdot\vec{w}\rg|\le\ C\, \lf[\beta(\tau)-\beta(\sigma)+ (\tau^2-\sigma^2)^{1/4}\, (\beta(\tau)-\beta(\sigma))^{3/4}\rg]\|\vec{w}\|_{\vec{\Phi}}
\ee
Let $\eta>0$ such that
\be
\label{A.I.8}
C\, \lf[\beta(\tau)-\beta(\sigma)+ (\tau^2-\sigma^2)^{1/4}\, (\beta(\tau)-\beta(\sigma))^{3/4}\rg]\le\eta
\ee
where $C>0$ is the constant in (\ref{A.I.7}) 

Assume 
\[
\forall\, \vec{\Phi}\in A\quad \quad A^\sigma(\vec{\Phi})>2\,\beta(\sigma)-\beta(\tau)\quad\Longrightarrow\quad \|DA^\tau(\vec{\Phi})\|>\eta
\]
Let $T_{\max}>0$ such that
\[
\forall\, t\in [0,T_{\max})\quad,\quad A^\sigma(\varphi_t(\vec{\Phi}))>2\,\beta(\sigma)-\beta(\tau)\quad\Longrightarrow\quad \|DA^\tau(\varphi_t(\vec{\Phi}))\|>\eta
\]
From \cite{Riv-columb} we have combining (\ref{A.I.7}) and (\ref{A.I.8})  that
\be
\label{A.I.9}
\forall\, t\in [0,T_{\max})\quad,\quad A^\sigma(\varphi_t(\vec{\Phi}))> \beta(\sigma)-\,2^{-1}(\beta(\tau)-\beta(\sigma))\quad\Longrightarrow\quad\frac{d A^\sigma(\varphi_t(\vec{\Phi}))}{dt}\le -\frac{\eta^2}{2}
\ee
Since, for all $t>0$ there must always exist $\vec{\Phi}\in A$ such that $A^\sigma(\varphi_t(\vec{\Phi}))\ge\beta(\sigma)$ (by definition of $\beta(\sigma)$) and since 
\be
\label{A.I.9-b}
\forall \, A\in {\mathcal A}_{\tau,\sigma}\quad \forall\ \vec{\Phi}\in A\quad\quad A^\sigma(\vec{\Phi})\le \beta(\sigma)+2\, (\beta(\tau)-\beta(\sigma))
\ee
we have
\be
\label{A.I.11}
\eta^2\ T_{\max}\le 4\, (\beta(\tau)-\beta(\sigma))\quad .
\ee
We chose $\beta(\tau)-\beta(\sigma)=\eta^3$ which implies $T_{\max}\le \eta$. Hence we have
\be
\label{A.I.12}
\begin{array}{l}
\ds\exists\, t\in[0,\eta]\quad\ \vec{\Phi}\in A\quad\mbox{ s. t. }\quad A^\sigma(\vec{\Phi})>2\,\beta(\sigma)-\beta(\tau)\quad,\quad \|D A^\sigma(\varphi_t(\vec{\Phi}))\| <(\beta(\tau)-\beta(\sigma))^{1/3}  \\[3mm]
\quad\mbox{ and }\quad 2\,\beta(\tau)-\beta(\sigma)> A^\tau(\varphi_t(\vec{\Phi})>A^\sigma(\varphi_t(\vec{\Phi}))>2\,\beta(\sigma)-\beta(\tau)
\end{array}
\ee
Using (\ref{A.I.1}) we have
\be
\label{A.I.13}
d_{\mathbf P}(\varphi_t(\vec{\Phi}),\vec{\Phi})\le C\, \sqrt{\eta}\quad\mbox{ and }\quad (\tau^2-\sigma^2)\,\int_{\Sigma^g}\lf[1+|{\vec{\mathbb I}}_{\varphi_t(\vec{\Phi})}|^2\rg]^2\ dvol_{g_{\vec{\Phi}}}\le 3\ (\beta(\tau)-\beta(\sigma))
\ee
This implies lemma~\ref{lm-A.3}. \hfill $\Box$

\medskip

\begin{Lma}
\label{deuxi}
Let $M^m$ be a closed sub-manifold of the euclidian space ${\R}^Q$. For any $W^{2,4}-$immersion $\vec{\Phi}$ of an oriented closed surface $\Sigma$ we denote
\[
F(\vec{\Phi}):=\int_\Sigma\lf(1+|\vec{\mathbb I}_{\vec{\Phi}}|^2_{g_{\vec{\Phi}}}\rg)^2\ dvol_{g_{\vec{\Phi}}}
\]
where $\vec{\mathbb I}_{\vec{\Phi}}$ is the second fundamental for of the immersion into $M^m$. The lagrangian $F$ is $C^2$ and there exists a constant $C$ depending only on $M^m$
such that for any perturbation $\vec{w}$ of the form $\vec{v}\circ\vec{\Phi}$ one has
\be
\label{deuxi-der1}
|DF(\vec{\Phi})(\vec{v}(\vec{\Phi}))|\le C\,\int_{\Sigma}\lf(1+|\vec{\mathbb I}_{\vec{\Phi}}|^2_{g_{\vec{\Phi}}}\rg)\ \lf[ (1+|\vec{\mathbb I}_{\vec{\Phi}}|^2_{g_{\vec{\Phi}}})\ |\p \vec{v}|(\vec{\Phi})+ |\vec{\mathbb I}_{\vec{\Phi}}|_{g_{\vec{\Phi}}}\ |\p^2 \vec{v}|(\vec{\Phi})\rg]\ dvol_{g_{\vec{\Phi}}}
\ee
and
\be
\label{deuxi-der2}
|D^2F(\vec{\Phi})(\vec{v}(\vec{\Phi}),\vec{v}(\vec{\Phi}))|\le C\,\int_{\Sigma}\lf(1+|\vec{\mathbb I}_{\vec{\Phi}}|^2_{g_{\vec{\Phi}}}\rg)\ \lf[ (1+|\vec{\mathbb I}_{\vec{\Phi}}|^2_{g_{\vec{\Phi}}})\ |\p \vec{v}|^2(\vec{\Phi})+\ |\p^2 \vec{v}|^2(\vec{\Phi})\rg]\ dvol_{g_{\vec{\Phi}}}
\ee
\hfill $\Box$
\end{Lma}
\noindent{\bf Proof of lemma~\ref{deuxi}.}

In local coordinates we  denote the second fundamental form
\[
\vec{\mathbb I}_{\vec{\Phi}}=\pi_{\vec{n}}\lf(d^2\vec{\Phi}   \rg)=\pi_{\vec{n}}\lf(\p^2_{x_ix_j}\vec{\Phi}\rg)\ dx_i\otimes dx_j
\]
we have
\be
\label{deuxi-1}
|\vec{\mathbb I}_{\vec{\Phi}}|^2_{g_{\vec{\Phi}}}:=\lf|\pi_{\vec{n}}\lf(d^2\vec{\Phi}   \rg)\rg|^2_{g_{\vec{\Phi}}}=\sum_{i,j,k,l}g^{ik}g^{jl}\pi_{\vec{n}}\p^2_{x_ix_j}\vec{\Phi}\cdot\pi_{\vec{n}}\p^2_{x_kx_l}\vec{\Phi}
\ee
Denote $\pi_T$ the projection onto the tangent plane to the immersion. We have in local coordinates
\be
\label{deuxi-2}
\pi_T(\vec{X})=\sum_{i,j=1}^2 g^{ij}\ \p_{x_i}\vec{\Phi}\cdot\vec{X}\ \p_{x_j}\vec{\Phi}
\ee
Hence
\be
\label{deuxi-3}
\lf.\pi_{\vec{n}}\frac{d\pi_{\vec{n}}}{dt}\rg|_{t=0}(\vec{X})=-\sum_{i,j=1}^2 g^{ij}\ \p_{x_i}\vec{\Phi}\cdot\vec{X}\ \pi_{\vec{n}}\lf(\p_{x_j}\vec{w}\rg)
\ee
We have clearly
\[
\frac{d g_{ij}}{dt}=\p_{x_i}\vec{\Phi}\cdot\p_{x_j}\vec{w}+\p_{x_i}\vec{w}\cdot\p_{x_j}\vec{\Phi}
\]
Hence
\be
\label{deuxi-4}
\frac{d g^{ij}}{dt}=-\,g^{ik} g^{jl}\,\lf[\p_{x_k}\vec{\Phi}\cdot\p_{x_l}\vec{w}+\p_{x_k}\vec{w}\cdot\p_{x_l}\vec{\Phi}\rg]:= -2\,(d\vec{\Phi}\dot{\otimes}_Sd\vec{w})^{ij}
\ee
We have then
\be
\label{deuxi-5}
\lf.\frac{d|\vec{\mathbb I}_{\vec{\Phi}}|^2_{g_{\vec{\Phi}}}}{dt}\rg|_{t=0}=2\, \lf<\pi_{\vec{n}}\lf(d^2\vec{\Phi}   \rg),\pi_{\vec{n}}\lf(D^{g_{\vec{\Phi}}}d\vec{w}\rg)\rg>_{g_{\vec{\Phi}}}-\,4\lf( g\otimes(d\vec{\Phi}\dot{\otimes}_Sd\vec{w})\res \vec{{\mathbb I}}_{\vec{\Phi}}\dot{\otimes}\vec{{\mathbb I}}_{\vec{\Phi}}\rg)
\ee
where $\res$ is the contraction operator between  $4-$contravariant and $4-$covariant tensors and
\be
\label{deuxi-6}
D^{g_{\vec{\Phi}}}d\vec{w}:=\lf[\p^2_{x_ix_j}\vec{w} -\sum_{rs=1}^2g^{rs}\p_{x_r}\vec{\Phi}\cdot \p^2_{x_ix_j}\vec{\Phi}\ \p_{x_s}\vec{w}\rg]\, dx_i\otimes dx_j\quad.
\ee
This gives in particular that
\be
\label{deuxi-7}
\begin{array}{l}
\ds\lf.\frac{d}{dt}\int_{\Sigma}(1+|\vec{\mathbb I}_{\vec{\Phi}}|^2_{g_{\vec{\Phi}}})^2\ dvol_{g_{\vec{\Phi}}}\rg|_{t=0}=DF(\vec{\Phi})(\vec{w})\\[5mm]
\ds\quad=4\, \int_{\Sigma}(1+|\vec{\mathbb I}_{\vec{\Phi}}|^2_{g_{\vec{\Phi}}}) \lf[\lf<\vec{\mathbb I}_{\vec{\Phi}},D^{g_{\vec{\Phi}}}d\vec{w}\rg>_{g_{\vec{\Phi}}}-\,2\lf( g\otimes(d\vec{\Phi}\dot{\otimes}_Sd\vec{w})\rg)\res \lf(\vec{{\mathbb I}}_{\vec{\Phi}}\dot{\otimes}\vec{{\mathbb I}}_{\vec{\Phi}}\rg)\rg]\ dvol_{g_{\vec{\Phi}}}\\[5mm]
\ds\quad+\int_{\Sigma}(1+|\vec{\mathbb I}_{\vec{\Phi}}|^2_{g_{\vec{\Phi}}})^2 \lf<d\vec{\Phi};d\vec{w}\rg>_{g_{\vec{\Phi}}}\ dvol_{g_{\vec{\Phi}}}
\end{array}
\ee
For $\vec{w}:=\vec{v}\lf(\vec{\Phi}\rg)$ we have 
\[
\begin{array}{l}
\ds\p^2_{x_ix_j}\vec{w} -\sum_{rs=1}^2g^{rs}\p_{x_r}\vec{\Phi}\cdot \p^2_{x_ix_j}\vec{\Phi}\ \p_{x_s}\vec{w}=\sum_{\al,\beta=1}^Q\p^2_{z_\al z_\beta}\vec{v}(\vec{\Phi})\, \p_{x_i}\vec{\Phi}^{\,\al}\,\p_{x_j}\vec{\Phi}^{\,\beta}\\[5mm]
\ds+\sum_{\al=1}^Q\p_{z_\al }\vec{v}(\vec{\Phi})\, \lf[\p^2_{x_ix_j}\vec{\Phi}^{\,\al}
-\ds\sum_{rs=1}^2g^{rs}\p_{x_r}\vec{\Phi}\cdot \p^2_{x_ix_j}\vec{\Phi}\  \p_{x_s}\vec{\Phi}^{\,\al}\rg]
\end{array}
\]
We have
\be
\label{deuxi-8}
\pi_{T}\lf(\p^2_{x_ix_j}\vec{\Phi}
\rg)=\sum_{rs=1}^2g^{rs}\,\p^2_{x_ix_j}\vec{\Phi}\cdot\p_{x_r}\vec{\Phi}\  \p_{x_s}\vec{\Phi}\quad.
\ee
Hence
\[
\p^2_{x_ix_j}\vec{\Phi}-\sum_{rs=1}^2g^{rs}\p_{x_r}\vec{\Phi}\cdot \p^2_{x_ix_j}\vec{\Phi}\  \p_{x_s}\vec{\Phi}=\pi_{\vec{n}}(\p^2_{x_ix_j}\vec{\Phi})=\vec{\mathbb I}_{ij}\quad.
\]
This implies that
\be
\label{deuxi-9}
D^{g_{\vec{\Phi}}}d\vec{w}=\sum_{\al,\beta=1}^Q\p^2_{z_\al z_\beta}\vec{v}(\vec{\Phi})\, d\vec{\Phi}^{\,\al}\otimes d \vec{\Phi}^{\,\beta}+\sum_{\al=1}^Q\p_{z_\al }\vec{v}(\vec{\Phi})\, \vec{\mathbb I}_{ij}^{\, \al}
\ee
We deduce
\be
\label{deuxi-10}
|DF(\vec{\Phi})(\vec{v}(\vec{\Phi}))|\le C\,\int_{\Sigma}\lf(1+|\vec{\mathbb I}_{\vec{\Phi}}|^2_{g_{\vec{\Phi}}}\rg)\ \lf[ (1+|\vec{\mathbb I}_{\vec{\Phi}}|^2_{g_{\vec{\Phi}}})\ |\p \vec{v}|(\vec{\Phi})+ |\vec{\mathbb I}_{\vec{\Phi}}|_{g_{\vec{\Phi}}}\ |\p^2 \vec{v}|(\vec{\Phi})\rg]\ dvol_{g_{\vec{\Phi}}}
\ee
We now compute the second derivative
\be
\label{deuxi-11}
\begin{array}{l}
\ds\lf.\frac{d}{dt}DF(\vec{\Phi}_t)(\vec{w})\rg|_{t=0}\\[5mm]
\ds=4\, \int_{\Sigma}(1+|\vec{\mathbb I}_{\vec{\Phi}}|^2_{g_{\vec{\Phi}}}) \lf[\lf<\vec{\mathbb I}_{\vec{\Phi}},D^{g_{\vec{\Phi}}}d\vec{w}\rg>_{g_{\vec{\Phi}}}-\,2\lf( g\otimes(d\vec{\Phi}\dot{\otimes}_Sd\vec{w})\rg)\res \lf(\vec{{\mathbb I}}_{\vec{\Phi}}\dot{\otimes}\vec{{\mathbb I}}_{\vec{\Phi}}\rg)\rg]\ \lf<d\vec{\Phi};d\vec{w}\rg>_{g_{\vec{\Phi}}}\ dvol_{g_{\vec{\Phi}}}\\[5mm]
\ds+\int_{\Sigma}|1+|\vec{\mathbb I}_{\vec{\Phi}}|^2_{g_{\vec{\Phi}}}|^2 \lf|\lf<d\vec{\Phi};d\vec{w}\rg>_{g_{\vec{\Phi}}}\rg|^2+8\, \lf|\lf<\vec{\mathbb I}_{\vec{\Phi}},D^{g_{\vec{\Phi}}}d\vec{w}\rg>_{g_{\vec{\Phi}}}-\,2\lf( g\otimes(d\vec{\Phi}\dot{\otimes}_Sd\vec{w})\rg)\res \lf(\vec{{\mathbb I}}_{\vec{\Phi}}\dot{\otimes}\vec{{\mathbb I}}_{\vec{\Phi}}\rg)\rg|^2\ dvol_{g_{\vec{\Phi}}}\\[5mm]
\ds+\,4\, \int_{\Sigma}(1+|\vec{\mathbb I}_{\vec{\Phi}}|^2_{g_{\vec{\Phi}}}) \frac{d}{dt}\lf[\lf<\vec{\mathbb I}_{\vec{\Phi}_t},D^{g_{\vec{\Phi}_t}}d\vec{w}\rg>_{g_{\vec{\Phi}_t}}-\,2\lf( g_{\vec{\Phi}_t}\otimes(d\vec{\Phi}_t\dot{\otimes}_Sd\vec{w})\rg)\res \lf(\vec{{\mathbb I}}_{\vec{\Phi}_t}\dot{\otimes}\vec{{\mathbb I}}_{\vec{\Phi}_t}\rg)\rg]\ dvol_{g_{\vec{\Phi}}}\\[5mm]
\ds+\int_{\Sigma}(1+|\vec{\mathbb I}_{\vec{\Phi}}|^2_{g_{\vec{\Phi}}})^2 \lf<d\vec{w};d\vec{w}\rg>_{g_{\vec{\Phi}}}\ dvol_{g_{\vec{\Phi}}}-2\,\int_{\Sigma} (1+|\vec{\mathbb I}_{\vec{\Phi}}|^2_{g_{\vec{\Phi}}})^2 \ \lf(d\vec{\Phi}\dot{\otimes}_S d\vec{w}\rg)\res \lf(d\vec{\Phi}\otimes d\vec{w}\rg)\ dvol_{g_{\vec{\Phi}}}
\end{array}
\ee
We have in one hand
\be
\label{deuxi-12}
\pi_{\vec{n}}\frac{d}{dt}\vec{\mathbb I}_{\vec{\Phi}_t}=\pi_{\vec{n}}\lf(D^{g_{\vec{\Phi}}}d\vec{w}\rg)
\ee
in the other hand
\be
\label{deuxi-13}
\begin{array}{l}
\ds\frac{d}{dt}\lf( \sum_{r=1}^2g^{rs}_{\vec{\Phi}_t}\ \p_{x_r}\vec{\Phi}_t\cdot \p^2_{x_ix_j}\vec{\Phi}_t \rg)=\sum_{r=1}^2g^{rs}\ \p_{x_r}\vec{\Phi}\cdot \p^2_{x_ix_j}\vec{w}+g^{rs}\ \p_{x_r}\vec{w}\cdot \p^2_{x_ix_j}\vec{\Phi}+\frac{d g^{rs}}{dt} \p_{x_r}\vec{\Phi}\cdot \p^2_{x_ix_j}\vec{\Phi}\\[5mm]
\ds\quad=\sum_{r=1}^2g^{rs}\ \p_{x_r}\vec{\Phi}\cdot \p^2_{x_ix_j}\vec{w}+g^{rs}\ \p_{x_r}\vec{w}\cdot \pi_{\vec{n}}(\p^2_{x_ix_j}\vec{\Phi})+\frac{d g^{rs}}{dt} \p_{x_r}\vec{\Phi}\cdot \p^2_{x_ix_j}\vec{\Phi}\\[5mm]
\ds\quad+\sum_{r,k,l=1}^2g^{rs}\ g^{kl}\ \p_{x_r}\vec{w}\cdot \p_{x_k}\vec{\Phi}\ \p_{x_l}\vec{\Phi}\cdot\p^2_{x_ix_j}\vec{\Phi}
\end{array}
\ee
we have
\be
\label{deuxi-14}
\begin{array}{l}
\ds\sum_{r=1}^2\frac{d g^{rs}}{dt} \p_{x_r}\vec{\Phi}\cdot \p^2_{x_ix_j}\vec{\Phi}=- \sum_{r,k,l=1}^2 g^{rk}\,g^{sl}\,\p_{x_r}\vec{\Phi}\cdot \p^2_{x_ix_j}\vec{\Phi}\, \lf[\p_{x_k}\vec{w}\cdot\p_{x_l}\vec{\Phi}+ \p_{x_l}\vec{w}\cdot\p_{x_k}\vec{\Phi} \rg]\\[5mm]
\ds=- \sum_{r,k,l=1}^2 g^{lk}\,g^{sr}\,\p_{x_l}\vec{\Phi}\cdot \p^2_{x_ix_j}\vec{\Phi}\ \p_{x_k}\vec{w}\cdot\p_{x_r}\vec{\Phi}     +  g^{lk}\,g^{sr}\,\p_{x_l}\vec{\Phi}\cdot \p^2_{x_ix_j}\vec{\Phi}\ \p_{x_r}\vec{w}\cdot\p_{x_k}\vec{\Phi}
\end{array}
\ee
Combining (\ref{deuxi-13}) and (\ref{deuxi-14}) we obtain
\be
\label{deuxi-15}
\ds\frac{d}{dt}\lf( \sum_{r=1}^2g^{rs}_{\vec{\Phi}_t}\ \p_{x_r}\vec{\Phi}_t\cdot \p^2_{x_ix_j}\vec{\Phi}_t \rg)=\sum_{r=1}^2g^{rs}\ \p_{x_r}\vec{\Phi}\cdot (D^{g_{\vec{\Phi}}}d\vec{w})_{ij}+\sum_{r=1}^2 g^{rs}\, \p_{x_r}\vec{w}\cdot\vec{\mathbb I}_{ij}
\ee
Thus
\be
\label{deuxi-16}
\ds\frac{d}{dt}\lf(D^{g_{\vec{\Phi}_t}}d\vec{w}\rg)=-\sum_{i,j=1}^2\lf[ \sum_{r=1}^2g^{rs}\ \p_{x_r}\vec{\Phi}\cdot (D^{g_{\vec{\Phi}}}d\vec{w})_{ij}\ \p_{x_s}\vec{w}+ g^{rs}\, \p_{x_r}\vec{w}\cdot\vec{\mathbb I}_{ij}\ \p_{x_s}\vec{w}\rg] \ dx_i\otimes dx_j
\ee
Combining (\ref{deuxi-4}), (\ref{deuxi-12}) and (\ref{deuxi-16}) we obtain
\be
\label{deuxi-17}
\begin{array}{l}
\ds \frac{d}{dt}\lf[\lf<\vec{\mathbb I}_{\vec{\Phi}_t},D^{g_{\vec{\Phi}_t}}d\vec{w}\rg>_{g_{\vec{\Phi}_t}}-\,2\lf( g_{\vec{\Phi}_t}\otimes(d\vec{\Phi}_t\dot{\otimes}_Sd\vec{w})\rg)\res \lf(\vec{{\mathbb I}}_{\vec{\Phi}_t}\dot{\otimes}\vec{{\mathbb I}}_{\vec{\Phi}_t}\rg)\rg]\\[5mm]
\ds=\lf|\pi_{\vec{n}}\lf(D^gd\vec{w}\rg)\rg|_{g_{\vec{\Phi}}}^2-\lf<\vec{\mathbb I}\ ;\sum_{i,j=1}^2g^{ij}\p_{x_i}\vec{\Phi}\cdot D^g d\vec{w}\ \p_{x_j}\vec{w}\rg>_{g_{\vec{\Phi}}}+ 4\, 
\lf[(d\vec{\Phi}\otimes_S d\vec{w})\otimes(d\vec{\Phi}\otimes_S d\vec{w})\rg]\res (\vec{\mathbb I}\dot{\otimes}\vec{\mathbb I})\\[5mm]
\ds-4\, \lf(g\otimes (d\vec{\Phi}\otimes_S d\vec{w})\rg)\res \lf(\vec{\mathbb I}\otimes \pi_{\vec{n}}(D^g d\vec{w})+\pi_{\vec{n}}(D^g d\vec{w})\otimes\vec{\mathbb I}\rg)-2\  \lf(g\otimes (d\vec{w}\otimes_S d\vec{w})\rg)\res (\vec{\mathbb I}\dot{\otimes}\vec{\mathbb I}) 
\end{array}
\ee
Combining (\ref{deuxi-11}) and (\ref{deuxi-17}) gives 
\be
\label{deuxi-18}
\begin{array}{l}
\ds D^2F(\vec{\Phi})(\vec{w},\vec{w})=\\[5mm]
\ds4\, \int_{\Sigma}(1+|\vec{\mathbb I}_{\vec{\Phi}}|^2_{g_{\vec{\Phi}}}) \lf[\lf<\vec{\mathbb I}_{\vec{\Phi}},D^{g_{\vec{\Phi}}}d\vec{w}\rg>_{g_{\vec{\Phi}}}-\,2\lf( g_{\vec{\Phi}}\otimes(d\vec{\Phi}\dot{\otimes}_Sd\vec{w})\rg)\res \lf(\vec{{\mathbb I}}_{\vec{\Phi}}\dot{\otimes}\vec{{\mathbb I}}_{\vec{\Phi}}\rg)\rg]\ \lf<d\vec{\Phi};d\vec{w}\rg>_{g_{\vec{\Phi}}}\ dvol_{g_{\vec{\Phi}}}\\[5mm]
\ds+\int_{\Sigma}|1+|\vec{\mathbb I}_{\vec{\Phi}}|^2_{g_{\vec{\Phi}}}|^2 \lf|\lf<d\vec{\Phi};d\vec{w}\rg>_{g_{\vec{\Phi}}}\rg|^2+8\, \lf|\lf<\vec{\mathbb I}_{\vec{\Phi}},D^{g_{\vec{\Phi}}}d\vec{w}\rg>_{g_{\vec{\Phi}}}-\,2\lf( g\otimes(d\vec{\Phi}\dot{\otimes}_Sd\vec{w})\rg)\res \lf(\vec{{\mathbb I}}_{\vec{\Phi}}\dot{\otimes}\vec{{\mathbb I}}_{\vec{\Phi}}\rg)\rg|^2\ dvol_{g_{\vec{\Phi}}}\\[5mm]
\ds+\,4\, \int_{\Sigma}(1+|\vec{\mathbb I}_{\vec{\Phi}}|^2_{g_{\vec{\Phi}}}) \lf[\lf|\pi_{\vec{n}}\lf(D^{g_{\vec{\Phi}}}d\vec{w}\rg)\rg|_{g_{\vec{\Phi}}}^2-\lf<\vec{\mathbb I}_{\vec{\Phi}}\ ;\sum_{i,j=1}^2g^{ij}\p_{x_i}\vec{\Phi}\cdot D^g_{g_{\vec{\Phi}}} d\vec{w}\ \p_{x_j}\vec{w}\rg>_{g_{\vec{\Phi}}}\rg] \ dvol_{g_{\vec{\Phi}}}\\[5mm]
\ds+\,16\, \int_{\Sigma}(1+|\vec{\mathbb I}_{\vec{\Phi}}|^2_{g_{\vec{\Phi}}})   \lf[(d\vec{\Phi}\otimes_S d\vec{w})\otimes(d\vec{\Phi}\otimes_S d\vec{w})\rg]\res (\vec{\mathbb I}_{\vec{\Phi}}\dot{\otimes}\vec{\mathbb I}_{\vec{\Phi}})          \ dvol_{g_{\vec{\Phi}}}\\[5mm]
\ds-\,16\, \int_{\Sigma}(1+|\vec{\mathbb I}_{\vec{\Phi}}|^2_{g_{\vec{\Phi}}})  \  \lf(g_{\vec{\Phi}}\otimes (d\vec{\Phi}\otimes_S d\vec{w})\rg)\res \lf(\vec{\mathbb I}_{\vec{\Phi}}\otimes \pi_{\vec{n}}(D^{g_{\vec{\Phi}}} d\vec{w}) +\pi_{\vec{n}}(D^{g_{\vec{\Phi}}} d\vec{w})\otimes\vec{\mathbb I}_{\vec{\Phi}}\rg)  \ dvol_{g_{\vec{\Phi}}}\\[5mm]
\ds-\, 8\, \int_{\Sigma}(1+|\vec{\mathbb I}_{\vec{\Phi}}|^2_{g_{\vec{\Phi}}})  \lf(g_{\vec{\Phi}}\otimes (d\vec{w}\otimes_S d\vec{w})\rg)\res (\vec{\mathbb I}_{\vec{\Phi}}\dot{\otimes}\vec{\mathbb I}_{\vec{\Phi}})  \ dvol_{g_{\vec{\Phi}}}\\[5mm]
\ds+\int_{\Sigma}(1+|\vec{\mathbb I}_{\vec{\Phi}}|^2_{g_{\vec{\Phi}}})^2 \lf<d\vec{w};d\vec{w}\rg>_{g_{\vec{\Phi}}}\ dvol_{g_{\vec{\Phi}}}-2\,\int_{\Sigma} (1+|\vec{\mathbb I}_{\vec{\Phi}}|^2_{g_{\vec{\Phi}}})^2 \ \lf(d\vec{\Phi}\dot{\otimes}_S d\vec{w}\rg)\res \lf(d\vec{\Phi}\otimes d\vec{w}\rg)\ dvol_{g_{\vec{\Phi}}}
\end{array}
\ee
For $\vec{w}:=\vec{v}(\vec{\Phi})$, using (\ref{deuxi-9}), we deduce
\be
\label{deuxi-19}
|D^2F(\vec{\Phi})(\vec{v}(\vec{\Phi}),\vec{v}(\vec{\Phi}))|\le C\,\int_{\Sigma}\lf(1+|\vec{\mathbb I}_{\vec{\Phi}}|^2_{g_{\vec{\Phi}}}\rg)\ \lf[ (1+|\vec{\mathbb I}_{\vec{\Phi}}|^2_{g_{\vec{\Phi}}})\ |\p \vec{v}|^2(\vec{\Phi})+\ |\p^2 \vec{v}|^2(\vec{\Phi})\rg]\ dvol_{g_{\vec{\Phi}}}
\ee
This concludes the proof of lemma~\ref{deuxi}.\hfill $\Box$

\end{document}